\documentclass[12pt,oneside]{amsbook}
\usepackage{cite}
\usepackage{latexsym}
\usepackage{layout}
\usepackage{amsfonts}
\usepackage{amssymb}
\usepackage{epsfig}

\makeindex

\begin{document}
\title{Notes on the Theory of Algebraic Numbers}
\author{Steve Wright} 
\maketitle
\markboth{}{}

\tableofcontents

\chapter{Motivation for Algebraic Number Theory: Fermat's Last Theorem}

\textbf{Fermat's Last Theorem (FLT)}. \textit{If $n\geq 3$ is an integer then there are no positive integers x, y, z such that $x^n+y^n=z^n$.}

\vspace{0.4cm}
This was first stated by Pierre de Fermat around 1637: in the margin of his copy of Bachet's edition of the complete works of Diophantus, Fermat wrote (in Latin):
\begin{quote}
``It is impossible to separate a cube into two cubes or a bi-quadrate into two bi-quadrates, or in general any power higher than the second into powers of like degree; I have discovered a truly remarkable proof which this margin is too small to contain."
\end{quote}
FLT was proved (finally!) by Andrew Wiles in 1995. One can show FLT if and only if one can prove

\vspace{0.4cm}
\textbf{FLT for primes}. \textit{If p is an odd prime then there are no nonzero integers x, y, z, such that $x^p+y^p+z^p=0$.}

\vspace{0.4cm}

Modern algebraic number theory essentially began in an attack on FLT by the great German number theorist Ernst Eduard Kummer in 1840. In order to explain Kummer's strategy, we need to recall the notion of a unique factorization domain (UFD).

\vspace{0.4cm}
\textbf{Definitions}. Let $D$ be an integral domain, i.e., a commutative ring with identity 1 that contains no zero divisors, i.e., no elements $x, y$ such that $x\not= 0\not= y$ and $xy=0$.

$(i)$ if $a, b\in D$ and if there exists $c\in D$ such that $b=ac$, then $a$ \emph{divides b}, or $a$ is a \emph{factor} of $b$, denoted $a|b$.

$(ii)$ $u\in D$ is a \emph{unit} in $D$ if $u$ has a mulitiplicative inverse in $D$.

If $U(D)$ denotes the set of all units in $D$ then $\{-1, 1\}\subseteq U(D)$ and $U(D)$ is an abelian group under multiplication in $D$; $U(D)$ is the \emph{group of units in D}.

$(iii)$ $a, b\in D$ are \emph{associates} if there exists $u\in U(D)$ such that $a=bu$.

$(iv)$ $0\not= p\in D$ is \emph{prime} (or \emph{irreducible}) if $p\not \in U(D)$ and $p=ab$ for $a, b \in D$ implies that either $a$ or $b$ is a unit.

\vspace{0.4cm}

If $a \in D,\ u\in U(D)$ then $a=(au^{-1})u$, i.e., \emph{every} element of $D$ has a factorization of the form (unit)$\times$(element of $D$). Such factorizations are hence said to be \emph{trivial}. Primes in $D$ are precisely the elements of $D$ with only trivial factorizations. If $\mathbb{Z}$ denotes the ordinary ring of integers then
\[
U(\mathbb{Z})=\{-1, 1\},
\]
and the primes in $\mathbb{Z}$, according to $(iv)$ in the above definitions, are precisely the prime numbers, together with their negatives. Hence a prime in an integral domain is the analog of a prime number in $\mathbb{Z}$. N.B. In order to avoid ambiguity that may arise from the terminology with regard to primes that we have introduced, we will refer to a positive prime number in $\mathbb{Z}$ as a \emph{rational prime}.

\vspace{0.4cm}
\textbf{Definition}. An integral domain $D$ is a \emph{unique factorization domain} (UFD) if

$(i)$ every element in $D\setminus \big(\{0\}\cup U(D)\big)$ can be factored into a product of primes in $D$, and

$(ii)$ If $p_1\cdots p_r$ and $q_1\cdots q_s$ are factorizations into primes of the same element of $D$ then $r=s$ and the $q_j$'s can be reindexed so that $p_i$ and $q_i$ are associates for $i=1,\dots,r$.

\vspace{0.4cm}

\begin{center}
\textit{Kummer's strategy for proving Fermat's Last Theorem}
\end{center}

\vspace{0.4cm}

Assume that $p$ is an odd  rational prime and that there exits nonzero integers $x, y, z$ such that
\begin{equation*}
x^p+y^p=z^p.\tag{$*$}
\]
We want to derive a contradiction from this assumption. In order to do this Kummer split the situation into the following two cases:

\emph{Case} I: $p$ divides none of $x,\ y,\ z$.

\emph{Case }II: $p$ divides at least one of $x,\ y,\ z$.

\noindent We will discuss what Kummer did only for Case I.

\vspace{0.2cm}
\emph{Notation}. In the sequel, if $S$ is a set and $n$ is a positive integer then $S^n$ will denote the Cartesian product of $S$ with itself taken $n$ times, i.e., the set of all $n$-tuples $(s_1,\dots,s_n)$, where $s_i\in S$ for all $i$.

\vspace{0.2cm}
It's easy to derive a contradiction for $p=3$. If neither $x, y,$ nor $z$ is divisible by 3, then $x^3,\ y^3$ and $z^3$ are each $\equiv \pm 1$ mod 9, and so $x^3+y^3 \equiv -2$, 0, or 2 mod 9, whence  $x^3+y^3 \not \equiv z^3$ mod 9, contrary to $(*)$ with $p=3$.

Hence suppose that $p>3$. Let
\[
\omega=\omega_p=e^{2\pi i/p}=\cos \left(\frac{2\pi}{p}\right)+i\sin \left(\frac{2\pi}{p}\right),\]
a \emph{p-th root of unity}. It can be shown (see Proposition 25, \emph{infra}) that the set $\mathbb{Z}[\omega]$ of complex numbers defined by
\[
\mathbb{Z}[\omega]=\left\{\sum_{i=0}^{p-2} a_i\omega^i: (a_0,\dots, a_{p-2})\in \mathbb{Z}^{p-1}\right\}\]
is a \emph{subring} of the set of complex numbers \textbf{C}, i.e., $\mathbb{Z}[\omega]$ is closed under addition, subtraction, and multiplication of complex numbers, and is also clearly an integral domain. Now, suppose that
\begin{quote}
$\mathbb{Z}[\omega]$ \textit{is a} UFD.
\end{quote}
Kummer then proved that 
\begin{equation*}
\textrm{there exists a unit $u\in \mathbb{Z}[\omega]$ and $\alpha\in \mathbb{Z}[\omega]$ such that $x+y\omega=u\alpha^p$}.\tag{$**$}\]
He then used $(*)$, $(**)$, and the assumption that $p$ does not divide $x$ or $y$ (Case I) to show that
\[
x\equiv y\ \textrm{mod}\ p.\]
Applying the same argument using $x^p+(-z)^p=(-y)^p$, he also got that
\[
x\equiv -z\ \textrm{mod}\ p.\]
But then
\[
2x^p\equiv x^p+y^p=z^p\equiv -x^p\ \textrm{mod}\ p,\]
and so
\[
3x^p\equiv 0\ \textrm{mod}\ p,\]
i.e., $p|(3x^p)$, hence $p|x$ or $p|3$. Because $p>3$, it follows that $p|x$, and this contradicts the hypothesis of Case I. Thus Kummer had shown that
\[
\textrm{if $\mathbb{Z}[\omega]$ is a UFD then Case I cannot be true,}\]
i.e., if $(*)$ is true and $\mathbb{Z}[\omega]$ is a UFD then $p$ must divide at least one of $x,\ y,$ or $z$. Kummer was thus led to ask
\[
\textrm{is  $\mathbb{Z}[\omega_p]$ a UFD, for all rational primes $p>3$?}\]
The answer, unfortunately, is no: Kummer was able to prove that $\mathbb{Z}[\omega_{23}]$ is in fact not a UFD. So the next question must be
\[
\textrm{for what $p$ is  $\mathbb{Z}[\omega_p]$ a UFD?}\]
Answer: all $p\leq 19$ and no others! This is very difficult to prove, and was not done until 1971. For some discussion of the ideas that Kummer used to derive a contradiction in Case I when  $\mathbb{Z}[\omega_p]$ is not a unique factorization domain, see Marcus [9], Chapter 1.

Let $\mathbb{Q}$ denote the set of all rational numbers, and let  $\mathbb{Q}[\omega_p]$ denote the set of complex numbers defined by
\[
\mathbb{Q}[\omega_p]=\left\{\sum_{i=0}^{p-2} a_i\omega_p^i: (a_0,\dots, a_{p-2})\in \mathbb{Q}^{p-1}\right\}.\]
Clearly  $\mathbb{Z}[\omega_p] \subseteq \mathbb{Q}[\omega_p]$, and we will eventually prove that $\mathbb{Q}[\omega_p]$ is a \emph{subfield} of \textbf{C}, i.e., $\mathbb{Q}[\omega_p]$ is closed under addition, subtraction, multiplication, and division of complex numbers. It turns out that arithmetic properties of $\mathbb{Z}[\omega_p]$ such as unique factorization and the existence of units with useful algebraic properties are closely tied to algebraic properties of  $\mathbb{Q}[\omega_p]$; indeed, much of Kummer's own work in number theory turns on a deep study of this connection. We are hence arrived at the fundamental questions of algebraic number theory:

(a) What are the \emph{subfields} $F$  of \textbf{C} which have a distinguished \emph{subring} $R$ such that (i) the \emph{arithmetic}, i.e., \emph{ring-theoretic}, structure of $R$ can be used to solve interesting and important problems in number theory and such that (ii) the arithmetic structure of $R$  can be effectively studied by means of the \emph{field-theoretic} structure of $F$?

(b) Given a class of fields $F$ and subrings $R$ of $F$ which answer question (a), what is the mathematical technology which can be used to get the ring-theoretic structure of $R$ from the field-theoretic structure of $F$?

\noindent We will spend our time in these notes getting some good answers to these very important questions.

\chapter{Complex Number Fields}

\textbf{Definition}. A \emph{complex number field} is a nonempty set $F$ of complex numbers such that $F\not= \{0\}$ and $F$ is closed under addition, subtraction, multiplication, and division of complex numbers, i.e., $F$ is a nonzero \emph{subfield} of \textbf{C}.

\vspace{0.4cm}
N.B. If $F$ a complex number field then $ \mathbb{Q} \subseteq F$.

\vspace{0.4cm}
\textbf{Definition}. If $A$ is a commutative ring with identity, a \emph{polynomial over A} is a  polynomial all of whose coefficients are in $A$.
\[
\textrm{Notation: $A[x]=$ the set of all polynomials in the indeterminant $x$ over $A$}.\]

\vspace{0.1cm}
$A[x]$ is a commutative ring with identity under the usual definitions of addition and multiplication of polynomials, and when $A$ is a field $F$, $F[x]$ is a Euclidean domain.

\vspace{0.4cm}
\textbf{Definitions}. Let $F$ be a complex number field. A complex number $\theta$ is \emph{algebraic over F} if there exits $p(x)\in F[x]$ such that $p\not \equiv 0$ and $p(\theta)=0.$

If $\theta$ is algebraic over $F$, let\[
M(\theta)=\{p\in F[x]: p\ \textrm{is \textit{monic} and $p(\theta)=0$}\}\]
(N.B. $M(\theta)\not= \emptyset$). An element of $M(\theta)$ of smallest degree is a \emph{minimal polynomial of} $\theta$ \emph{over} $F$.

\vspace{0.4cm}
\textbf{Proposition 1}. \textit{If $\theta$ is algebraic over $F$ then there is only one element of $M(\theta)$ of smallest degree, i.e., the minimal polynomial of $\theta$ over F is} unique. 

\vspace{0.4cm}
\emph{Proof}. Let $p,q$ be minimal polynomials of $\theta$ over $F$. Because $F[x]$ is a Euclidean domain, there exits $d,r\in F[x]$ such that 
\[
q=dp+r,\ r\equiv 0\ \textrm{or the degree of $r<$ the degree of $p$}.\]
Hence
\[
r(\theta)=q(\theta)-d(\theta)p(\theta)=0.\]
If $r\not \equiv 0$ then, upon dividing $r$ by its leading coefficient, we obtain a monic polynomial over $F$ of lower degree that $p$ and not identically 0 which has $\theta$ as a root, impossible since $p$ is a minimal polynomial of $\theta$ over $F$. Hence $r\equiv 0$ and so $p$ divides $q$ in $F[x]$. Similarly, $q$ divides $p$ in $F[x]$. Hence $p=\alpha q$ for some $\alpha\in F$, and because $p$ and $q$ are both monic, $\alpha=1$, and so $p=q$.$\hspace{15.3cm} \textrm{QED}$

\vspace{0.4cm}
\textbf{Definition}. If $\theta$ be algebraic over $F$ then the \emph{degree of} $\theta$ \emph{over F} is the degree of the minimal polynomial of $\theta$ over $F$.

\vspace{0.4cm}
The proof of Proposition 1 implies the following corollary, which we will use frequently in the sequel.

\vspace{0.4cm}
\textbf{Corollary 2}. \textit{If $\theta$ is algebraic over F, p is the minimal polynomial of $\theta$ over F, and $q\in F[x]$ has $\theta$ as a root, then p divides q in $F[x]$.}

\vspace{0.4cm}
\textbf{Definition}. A polynomial $p$ in $F[x]$ is \emph{irreducible over F} if there do not exist \emph{nonconstant} polynomials $q,\ r\in F[x]$ such that $p=qr$.

\vspace{0.4cm}
\textbf{Proposition 3}. \textit{If $\theta$ is algebraic over $F$ with minimal polynomial p then p is irreducible over F.}

\emph{Proof}. Suppose that $p=qr,\ q,\ r$ both nonconstant polynomials in $F[x]$. Then degree of $q<$ degree of $p$, degree of $r<$ degree of $p$, and $\theta$ is a root of either $q$ or $r$, contrary to the minimality of the degree of $p$. $\hspace{10.2cm} \textrm{QED}$

\vspace{0.4cm}
\textbf{Proposition 4}. \textit{A polynomial p of degree $n>0$ irreducible over a complex number field F has $n$ distinct roots in $\mathbf{C}$.}

\emph{Proof}. We may assume with no loss of generality that $p$ is monic. Suppose that $p$ has a double root $\alpha \in \mathbf{C}$; then\[
p(x)=(x-\alpha)^2q(x),\ q\in \mathbf{C}[x].\]
Hence
\[
p^{\prime}(x)=(x-\alpha)^2q^{\prime}(x)+2(x-\alpha)q(x),\]
(the $\prime$ here denotes differentiation) and so $p$ and $p^{\prime}$ have $\alpha$ as a common root. Now $p^{\prime} \in F[x]$ and $\alpha$ is algebraic over $F$ (because $p\in F[x]$), hence Corollary 2 implies that the minimal polynomial $r$ of $\alpha$ must divide both $p$ and $p^{\prime}$. But $p$ is irreducible, so $p=r$ ($p$ and $r$ are both monic), and so $p$ divides $p^{\prime}$ in $F[x]$, which is impossible because $p^{\prime}\not \equiv 0$ and the degree of $p^{\prime}<$ the degree of $p$. Hence $p$ does not have a double root; because every polynomial over \textbf{C} of degree $n$ has $n$ roots in \textbf{C} counted according to multiplicity, it follows that $p$ has $n$ distinct roots. $\hspace{14.3cm} \textrm{QED}$

\vspace{0.4cm}
\textbf{Definition}. Let $\theta$ be of degree $n$ over $F$. The $n$ (distinct) roots of the minimal polynomial of $\theta$ over $F$ are called the \emph{conjugates of $\theta$ over F}.

\vspace{0.4cm}
\begin{quote}
N.B. If $\theta$ is algebraic over $F$, with minimal polynomial $p$, then the conjugates of $\theta$ over $F$ is the set of all complex numbers algebraic over $F$ which have $p$ as their minimal polynomial over $F$.\end{quote}

\vspace{0.4cm}
\textbf{Theorem 5}. \textit{If F is a complex number field then the set of all complex numbers algebraic over F is a complex number field which contains F}.

\vspace{0.4cm}
The proof of this very important theorem requires some useful mathematical technology involving symmetric polynomials.

\vspace{0.4cm}
\textbf{Definition}. A polynomial $p$ over $F$ in the variables $x_1,\dots,x_m$ is \emph{symmetric} if the following condition holds: if $\sigma$ is a permutation of $\{1,\dots,m\}$ then
\[
p(x_1,\dots,x_m)=p(x_{\sigma(1)},\dots,x_{\sigma(m)}),\]
i.e., $p$ remained unchanged under any permutation of its variables.

\vspace{0.4cm}
\textbf{Definition}. Let $n$ be a positive integer. The polynomials
\[
 \sigma_1=\sum_{i=1}^n x_i,\]

\hspace{7cm}  \vdots 

\hspace{6.4cm} $\sigma_i=\textrm{sum of all products of $i$ different $x_j$},$

\hspace{7cm}  \vdots 
\[
\sigma_n=\prod_{i=1}^n x_i\]
are all symmetric over any complex number field. They are called the \emph{elementary symmetric polynomials in n variables}.

\vspace{0.4cm}
\textbf{Proposition 6}. \textit{If $(\theta_1,\dots,\theta_n)$ is an n-tuple of complex numbers then}
\[
\prod_{i=1}^n (x-\theta_i)=x^n+\sum_{i=1}^n (-1)^i\sigma_i(\theta_1,\dots,\theta_n) x^{n-i}.\]

\vspace{0.4cm}
\emph{Proof}. This follows by a straightforward induction on $n$. $\hspace{5cm} \textrm{QED}$

\vspace{0.4cm}
\textbf{Corollary 7}. \textit{If $\theta$ is algebraic over F of degree n, $\sigma_i$ is an elementary symmetric polynomial in n variables and $\theta_1,\dots,\theta_n$  are the conjugates of $\theta$ over F, then $\sigma_i(\theta_1,\dots,\theta_n)\in F$}.

\vspace{0.4cm}
\emph{Proof}. If $p$ is the minimal polynomial of $\theta$ over $F$ then Proposition 6 implies that
\begin{equation*}
p(x)=\prod_{i=1}^n (x-\theta_i)=x^n+\sum_{i=1}^n (-1)^i\sigma_i(\theta_1,\dots,\theta_n) x^{n-i}.\tag{1}\]
All coefficients of $p$ are in $F$, hence all coefficients of the polynomial on the right-hand side of $(1)$ are also in $F$.$\hspace{12.1cm}\textrm{QED}$
 
 We also need the following result from the classical theory of equations:

\vspace{0.4cm}
\textbf{Lemma 8}. \textit{Let $m,\ n$ be fixed positive integers,  $\tau_1\dots,\tau_m$ (respectively, $\sigma_1,\dots,\sigma_n$) the elementary symmetric polynomials in m (respectively, n) variables. Let $p$ be a polynomial over $F$ in the variables $x_1,\dots,x_m,\ y_1\dots,y_n$ with the following property: if $\sigma$ (respectively, $\tau$) is a permutation of $\{1,\dots,m\}$ (respectively,  $\{1,\dots,n\}$) then
\[
p(x_1,\dots,x_m, y_1\dots,y_n)=p(x_{\sigma(1)},\dots,x_{\sigma(m)}, y_{\tau(1)},\dots,y_{\tau(n)}),\]
i.e., p remains unchanged when the $x_i$'s and the $y_j$'s are permuted among themselves. Then there exits a polynomial $q$ over $F$ in the variables $x_1,\dots,x_m,\ y_1\dots,y_n$ such that}
\begin{eqnarray*}
p(x_1,\dots,&x_m,& y_1,\dots,y_n)\\
&=&q\big(\tau_1(x_1,\dots,x_m),\dots,\tau_m(x_1,\dots,x_m),\sigma_1( y_1,\dots,y_n),\dots,\sigma_n(y_1,\dots,y_n)\big).
\end{eqnarray*}

\vspace{0.4cm}
\emph{Proof}. See Weisner [11], Theorem 49.10.$\hspace{7.5cm} \textrm{QED}$

\vspace{0.4cm}
\emph{Proof of Theorem} 5.

Let $\alpha$ and $\beta$ be algebraic over $F$. We want to show that $\alpha \pm \beta, \alpha \beta,$ and $\alpha/\beta,$ provided that $\beta \not=0$, are all algebraic over $F$. 
We will do this by the explicit construction of  polynomials over $F$ that have these numbers as roots.

Start with $\alpha+\beta$. Let $f$ and $g$ denote the minimal polynomials of, respectively, $\alpha$ and $\beta$, of degree $m$ and $n$, respectively. Let $\alpha_1,\dots,\alpha_m$ and $\beta_1,\dots,\beta_n$ denote the roots of $f$ and $g$ in $\textbf{C}$, with $\alpha_1=\alpha$ and $\beta_1=\beta$. Now consider the polynomial
\begin{equation*}
\prod_{i=1}^m \prod_{j=1}^n (x-\alpha_i-\beta_j)=x^{mn}+\sum_{i=1}^{mn} c_i(\alpha_1,\dots,\alpha_m, \beta_1,\dots,\beta_n)x^{mn-i}, \tag{2}\]
where each coefficient $c_i$ is a polynomial in the $\alpha_i$'s and $\beta_j$'s over $F$ (in fact, over $\mathbb{Z}$). We claim that 
\begin {equation*}
c_i(\alpha_1,\dots,\alpha_m, \beta_1,\dots,\beta_n) \in F,\ i=1,\dots,mn. \tag{3} \] 
If this is true then the polynomial $(2)$ is in $F[x]$ and has $\alpha_1+\beta_1=\alpha+\beta$ as a root, whence $\alpha+\beta$ is algebraic over $F$.

In order to verify (3), observe that the left-hand side of $(2)$ remains unchanged when the $\alpha_i$'s and the $\beta_j$'s are permuted amongst themselves (this simply rearranges the order of the factors in the product), and so the same thing is true for each coefficient $c_i$. It thus follows from Lemma 8 that there exists a polynomial $l_i$ over $F$ in the variables $x_1,\dots,x_m, y_1,\dots,y_n$ such that
\begin{eqnarray*}
c_i(\alpha_1,\dots,&\alpha_m,& \beta_1,\dots,\beta_n)\\
&=&l_i(\tau_1(\alpha_1,\dots,\alpha_m),\dots,\tau_m(\alpha_1,\dots,\alpha_m),\sigma_1( \beta_1,\dots,\beta_n),\dots,\sigma_n( \beta_1,\dots,\beta_n)).
\end{eqnarray*}
where $\tau_1,\dots,\tau_m, \sigma_1,\dots,\sigma_n$ denote, respectively, the elementary symmetric polynomials in $m$ and $n$ variables. It now follows from Corollary 7 that each of the numbers at which $l_i$ is evaluated in this equation is in $F$, and (3) is an immediate consequence of that.

A similar argument shows that $\alpha-\beta$ and $\alpha\beta$ are algebraic over $F$.

Suppose next that $\beta\not=0$ is algebraic over $F$ and let
\[
x^n+\sum_{i=0}^{n-1} a_ix^i\]
be the minimal polynomial of $\beta$ over $F$. Then $1/\beta$ is a root of
\[
1+\sum_{i=0}^{n-1} a_ix^{n-i} \in F[x], \]
and so $1/\beta$ is algebraic over $F$. Then $\alpha/\beta=\alpha \cdot (1/\beta)$ is algebraic over $F$. $\hspace{2.3cm}\ \textrm{QED}$

\vspace{0.4cm}
\textbf{Definitions}. A complex number $\theta$ is \emph{algebraic} if $\theta$ is algebraic over $\mathbb{Q}$. If $\theta$ is not algebraic then $\theta$ is \emph{transcendental}.

\vspace{0.4cm}
The numbers $e$ and $\pi$ are transcendental: for a proof, see Hardy and Wright [6], Theorems 204 and 205.

\vspace{0.4cm}
\textbf{Theorem 9} \textit{$($Gelfand-Schneider, $1934)$. If $\alpha$ and $\beta$ are algebraic, $0\not= \alpha\not= 1$, and $\beta$ is irrational, then $\alpha^{\beta}$ is transcendental}.

\vspace{0.4cm}
\emph{Proof}. See A. Baker [1], Theorem 2.1.$\hspace{8.5cm} \textrm{QED}$

Hilbert's seventh problem asks: is $2^{\sqrt{2}}$ transcendental? Theorem 9 answers: yes it is!

\chapter{Extensions of Complex Number Fields}

\vspace{0.4cm}
\textbf{Definition}. A complex number field $K$ is an \emph{extension} of a complex number field $F$ if $F\subseteq K$.

\vspace{0.4cm}
N.B. Every complex number field is an extension of $\mathbb{Q}$.

If $K$ is an extension of $F$ then $K$ becomes a vector space over $F$ with the vector addition and scalar multiplication defined by the addition and multiplication of elements in $K$.

\vspace{0.4cm}
\textbf{Definitions}. If $K$ is an extension of $F$ then the \emph{degree of K over F}, denoted $[K:F]$, is the dimension of $K$ as a vector space over $F$. $K$ is a \emph{finite extension of F} if $[K:F]$ is finite.

\vspace{0.4cm}
\textbf{Definition}. If $\theta \in \mathbf{C}$ and $F$ is a complex number field then $F(\theta)$ denotes the smallest subfield of \textbf{C} that contains $F$ and $\theta$. $F(\theta)$ is called a \emph{simple extension of F}.

\vspace{0.4cm}
It is easy to see that
\[
F(\theta)=\left\{ \frac{f(\theta)}{g(\theta)}: f,\ g\in F[x],\ g(\theta)\not= 0\right\}.\]
We will now show that if $\theta$ is algebraic over $F$ then the structure of $F(\theta)$ simplifies considerably.

\vspace{0.4cm}
\textbf{Proposition 10}. \textit{If F is a complex number field and $\theta$ is algebraic over F of degree n then $[F(\theta): F]=n$ and $\{1, \theta,\dots, \theta^{n-1}\}$ is a basis of $F(\theta)$ over F}.

\vspace{0.4cm}
\emph{Proof}. The set $\{1, \theta,\dots, \theta^{n-1}\}$ is linearly independent over $F$; otherwise, $\theta$ is the root of a nonzero polynomial over $F$ of degree $<n$, and that is impossible. Now  let $p$ be the minimal polynomial of $\theta$ over $F$,
\[
\alpha=\frac{f(\theta)}{g(\theta)}\in F(\theta).\]
Because $g(\theta)\not= 0=p(\theta)$, $p$ does not divide $g$ in $F[x]$. But $p$ is irreducible over $F$, hence $p$ and $g$ have no non-constant common factor in $F[x]$, i.e., $p$ and $g$ are relatively prime in $F[x]$. As $F[x]$ is a Euclidean domain, it follows that there exist polynomials $s$ and $t$ over $F$ such that
\[
tp+sg=1.\]
Evaluate this equation at $x=\theta$; since $p(\theta)=0$, we get
\[
s(\theta)g(\theta)=1,\ \textrm{i.e.,}\ s(\theta)=\frac{1}{g(\theta)}.\]
Hence
\[
\alpha=\frac{f(\theta)}{g(\theta)}=f(\theta)s(\theta)=h(\theta)\]
for some polynomial $h$ over $F$. Divide $h$ by $p$ in $F[x]$ to obtain $q,\ r\in F[x]$ such that
\[
h=qp+r,\ \textrm{degree of}\ r<\ \textrm{degree of}\ p=n, \]
and then evaluate at $x=\theta$ to obtain $h(\theta)=r(\theta)$. Thus
\[
\alpha=r(\theta)\ \textrm{with}\ r\in F[x],\ \textrm{degree of}\ r\leq n-1,\]
hence $\alpha$ is in the linear span of $\{1, \theta,\dots, \theta^{n-1}\}$ over $F$.$\hspace{6.2cm} \textrm{QED}$

Proposition 10 implies that if $\theta$ is algebraic over $F$ then
\[
[F(\theta):F]=\ \textrm{degree of}\ \theta\ \textrm{over}\ F.\]

\vspace{0.4cm}
\textbf{Definition}. An \emph{algebraic number field}, or, as we will sometimes say more succinctly, a \emph{number field}, is a complex number field that is a \emph{finite} extension of $\mathbb{Q}$.

\vspace{0.4cm}

\begin{center}
\textit{Two important examples}
\end{center}

\vspace{0.4cm}
(1) \textit{Quadratic number fields}

\vspace{0.4cm}

Let $m$ be a square-free integer, i.e., $m$ does not have a nontrivial perfect square as a factor. Then $\sqrt{m}$ is irrational (why?), hence $x^2-m$ is the minimal polynomial of 
 $\sqrt{m}$ over $\mathbb{Q}$. Thus $\mathbb{Q}(\sqrt{m})$ has degree 2 over $\mathbb{Q}$ and $\{1, \sqrt{m} \}$ is a basis of $\mathbb{Q}(\sqrt{m})$ over $\mathbb{Q}$. Hence
 \[
 \mathbb{Q}(\sqrt{m})=\{a+b\sqrt{m}: (a, b)\in \mathbb{Q} \times \mathbb{Q}\}.\]

\vspace{0.4cm}
\textbf{Definitions}. $\mathbb{Q}(\sqrt{m})$ is the \emph{quadratic number field determined by m}. If $m>0$ then  $\mathbb{Q}(\sqrt{m})$ is a \emph{real} quadratic number field and if $m<0$ then  $\mathbb{Q}(\sqrt{m})$ is an \emph{imaginary} quadratic number field.

\vspace{0.4cm}
(2)\textit{Cyclotomic number fields}.

\vspace{0.4cm}
Let $m$ be an integer, $m\geq 2$. Set
\[
\omega_m=e^{2\pi i/m}=\cos \frac{2\pi}{m}+i\sin \frac{2\pi}{m}.\]
The set of complex numbers $\{1,\dots,\omega_m^{m-1}\}$ is an abelian group under multiplication, cyclic of order $m$. An element $\omega^k$ of this group has order $m$, and is hence a generator of this group, if and only if the greatest common divisor of $k$ and $m$ is 1. 

\vspace{0.4cm}
\textbf{Definition}. A generator of the group $\{1,\dots,\omega_m^{m-1}\}$ is called a \emph{primitive $m$-th root of unity}.
 
\vspace{0.4cm}
\emph{Notation}. We will denote the greatest common divisor of the integers $i$ and $j$ by $\gcd(i, j)$.

\vspace{0.4cm}
The number $\omega_m$ is algebraic (it's a root of $x^m-1$); what is it's minimal polynomial (over $\mathbb{Q}$)?

\vspace{0.4cm}
\textbf{Lemma 11}. \textit{The conjugates of $\omega_m$ over $\mathbb{Q}$ are precisely the primitive $m$-th roots of unity}.

\vspace{0.4cm}
\emph{Proof}. See Marcus [9], Chapter 2, Theorem 3 (this is not obvious and requires some work). $\hspace{14.2cm} \textrm{QED}$

\vspace{0.4cm}
\textbf{Definition}. The  \emph{m-th cyclotomic polynomial} is the polynomial
\[
\Phi_m(x)=\prod_{k:\ k \in  \mathbb{Z},\ 1\leq k \leq m,\ \gcd(k, m)=1} (x-\omega_m^k).\]

\vspace{0.4cm}

Lemma 11 implies that
\[
\Phi_m(x)\ \textrm{is the minimal polynomial of $\omega_m$ over $\mathbb{Q}$.}\]

For each positive integer $n$, let 
\[
\varphi(n)=\ \textrm{the cardinality of the set $\{k: k \in  \mathbb{Z},\ 1\leq k \leq n,\ \gcd(k, n)=1\}$}:\]
The function $\varphi$ so defined is called \emph{Euler's totient function}. Proposition 10 implies that
\[
[\mathbb{Q}(\omega_m):\mathbb{Q}]=\ \textrm{the degree of}\ \Phi_m(x)=\varphi(m),\]
hence
\[
\mathbb{Q}(\omega_m)=\left\{ \sum_{i=0}^{\varphi(m)-1} a_i \omega_m^i: (a_0,\dots,a_{\varphi(m)-1}) \in \mathbb{Q}^{\varphi(m)} \right\}.\]

\vspace{0.4cm}
\textbf{Definition}. The  \emph{m-th cyclotomic number field} is the number field $\mathbb{Q}(\omega_m)$.

\vspace{0.4cm}
Now assume that $m=p,$ $p$ an odd prime. Then
\[
\textrm{the degree of $\omega_p$ over $\mathbb{Q}$}=\varphi(p)=p-1.\]
Division of $x-1$ into $x^p-1$ produces the factorization
\begin{equation*}
x^p-1=(x-1)\sum_{i=0}^{p-1} x^i.\tag{*}\]
The numbers $1, \omega_p,\dots,\omega_p^{p-1}$ are precisely the $p$-th roots of unity, and they are all primitive $p$-th roots, hence it follows from $(*)$ and the definition of the $p$-th cyclotomic polynomial $\Phi_p(x)$ that
\[
\Phi_p(x)=\prod_{i=1}^{p-1}(x-\omega_p^i)=\frac{ \prod_{i=0}^{p-1}(x-\omega_p^i)}{x-1}=\frac{x^p-1}{x-1}=\sum_{i=0}^{p-1} x^i.\]
In particular, $\Phi_p(x) \in \mathbb{Z}[x]$. This is no accident; we will see eventually that for all integers $m\geq 2$,  $\Phi_m(x) \in \mathbb{Z}[x]$. 

\chapter{The Primitive Element Theorem}

\vspace{0.4cm}
\textbf{Definition}. If $F$ is a complex number field and $\{\alpha_1,\dots,\alpha_n\} \subseteq \mathbf{C}$ then $F(\alpha_1,\dots,\alpha_n)$ is the smallest subfield of \textbf{C} that contains $F$ and $\{\alpha_1,\dots,\alpha_n\}$.

\vspace{0.4cm}
\textbf{Proposition 12}. If $F[x_1,\dots,x_n]$ denotes the set of all polynomials in the variables $x_1,\dots,x_n$ with coefficients in $F$ then
\begin{equation*}
F(\alpha_1,\dots,\alpha_n)=\left\{\frac{p(\alpha_1,\dots,\alpha_n)}{q(\alpha_1,\dots,\alpha_n)}: p, q\in F[x_1,\dots,x_n],\ q(\alpha_1,\dots,\alpha_n)\not=0\right\};\tag{\textit{i}}\]

\noindent $(ii)$ if each $\alpha_i$ is algebraic over $F$ of degree $d_i$ then
\[
F(\alpha_1,\dots,\alpha_n)=\Big\{p(\alpha_1,\dots,\alpha_n): p\in F[x_1,\dots,x_n],\ \textrm{degree of}\ p\leq \sum_i (d_i-1)\Big\}.\]

\vspace{0.4cm}
\emph{Proof}. Exercise. $\hspace{11.9cm} \textrm{QED}$

\vspace{0.4cm}
\textbf{Definition}. If $K$ and $F$ are complex number fields with $F\subseteq K$ and if there is an element $\theta$ of $K$ such that $K=F(\theta)$ then $\theta$ is called a \emph{primitive element of K over F}.

\vspace{0.4cm}
The following theorem, which will play a crucial role in the sequel, asserts that all finite extensions of a complex number field are simple extensions.

\vspace{0.4cm}
\textbf{The Primitive Element Theorem}. \textit{If K is a finite extension of a complex number field F then every element of K is algebraic over F and there exits a primitive element of K over F. Moreover, the degree over F of any primitive element of K over F is equal to the degree of K over F}.

\vspace{0.4cm}
\emph{Proof}. If $\alpha\in K$ is not algebraic over $F$ then $\{\alpha^i: i=1, 2, 3,\dots \}$ is an infinite subset of $K$ that is linearly independent over $F$, which is not possible because $K$ is a finite extension of $F$.

Let $\{\alpha_1,\dots,\alpha_n\}$ be a basis of $K$ over $F$. Then each $\alpha_i$ is algebraic over $F$ and $K=F(\alpha_1,\dots,\alpha_n)$. We now apply the following lemma to obtain a primitive element of $K$ over $F$. 

\vspace{0.4cm}
\textbf{Lemma 13}. \textit{If $\beta_1,\dots,\beta_k$ are algebraic over $F$ then $F(\beta_1,\dots,\beta_k)$ is a simple extension of $F$}.

\vspace{0.4cm}
\emph{Proof}. Note first that Proposition 12$(ii)$ implies that $F(\beta_1,\dots,\beta_k)$ is a finite extension of $F$, and so by what we have already shown,
\begin{equation*}
\textrm{every element of  $F(\beta_1,\dots,\beta_k)$ is algebraic over $F$.}\tag{1}
\]
We will prove Lemma 13 for $k=2$; once that is done, the general case will follow from $(1)$ and induction on $k$.

Hence, suppose that $\alpha, \beta$ are algebraic over $F$: we must prove that there is a $\theta\in F(\alpha, \beta)$ such that $F(\alpha, \beta)=F(\theta)$. Let $\alpha_1,\dots,\alpha_m,\ \beta_1,\dots,\beta_n$ be the conjugates of $\alpha$ and $\beta$ over $F$, with $\alpha_1=\alpha,\ \beta_1=\beta$. 
Since the conjugates over $F$ are distinct (Propositions 3 and 4), it follows that $\beta\not= \beta_j$ for $j\not= 1$. Hence for each $j\not= 1$ and for each $i$,
\[
\alpha_i+x\beta_j=\alpha+x\beta\]
has exactly one solution $x$ in $F$. Because there are only finitely many of these equations, we can find $c\in F$ such that
\begin{equation*}
\alpha_i+c\beta_j\not =\alpha+c\beta,\ \textrm{for all}\ i\ \textrm{and for all}\ j\not= 1.\tag{2}\]
Let $\theta=\alpha+c\beta$. This is the $\theta$ that we want: we clearly have $F(\theta)\subseteq F(\alpha, \beta)$, and to prove the reverse inclusion, it suffices to show that $\beta\in F(\theta)$, for then $\alpha=\theta-c\beta \in F(\theta)$.

In order to prove that $\beta\in F(\theta)$, consider the minimal polynomials $f$ and $g$ of $\alpha$ and $\beta$ over $F$. Since
\[
f(\theta-c\beta)=f(\alpha)=0,\]
$\beta $ is a common root of $g(x)$ and $f(\theta-cx)$. 
These polynomials have no other common root; if they did, then for some $j\not= 1$, $\theta-c\beta_j$ is a root $\alpha_i$ of $f$, contrary to $(2)$.

The polynomials $g(x)$ and $f(\theta-cx)$ are in $F(\theta)[x]$ with unique common root $\beta$. If\[
h(x)=\ \textrm{the minimal polynomial of $\beta$ over $F(\theta)$,}\]
then Corollary 2 implies that $h(x)$ divides both $g(x)$ and $f(\theta-cx)$ in $F(\theta)[x]$. We claim that the degree of $h(x)$ must be 1; otherwise $g(x)$ and $f(\theta-cx)$ would have more than one common root, namely the at least two distinct roots of $h(x)$! Hence there exits $ \delta \in F(\theta)$ such that
\[
h(x)=x+\delta.\]
But then $0=h(\beta)= \beta+ \delta$, i.e., $\beta=-\delta \in F(\theta)$.$\hspace{6.8cm} \textrm{QED}$

Finally, we note that if $\theta$ is a primitive element of $K$ over $F$ and $m$ is the degree of $\theta$ over $F$ then Proposition 10 implies that $m=[K:F]$.$\hspace{6.8cm} \textrm{QED}$

\vspace{0.4cm}

The proof of Lemma 13 can be used to easily find primitive elements of finite extensions of complex number fields. For example, suppose that we wish to find a primitive element of $\mathbb{Q}(\sqrt{2}, \sqrt [3]{3})$ over $\mathbb{Q}$. Note first that the conjugates of $\sqrt{2}$ over $\mathbb{Q}$ are $\sqrt{2}$ and $-\sqrt{2}$ and the conjugates of $\sqrt [3]{3}$ over $\mathbb{Q}$ are $\sqrt [3]{3},\ \omega \sqrt [3]{3},$ and $\omega^2 \sqrt [3]{3}$, where $\omega=e^{2\pi i/3}.$ Next, in the proof of Lemma 13, we take $\alpha_1,\ \alpha_2,\ \alpha_3$ to be the conjugates of $\sqrt [3]{3}$ and $\beta_1, \beta_2$ to be the conjugates of $\sqrt{2}$, hence we must find $c\in \mathbb{Q}$ such that
\[
\sqrt [3]{3}-c\sqrt{2}\not=\sqrt [3]{3}+c\sqrt{2},\]
\[
 \omega \sqrt [3]{3}-c\sqrt{2}\not=\sqrt [3]{3}+c\sqrt{2},\]
 \[
\omega^2 \sqrt [3]{3} -c\sqrt{2}\not=\sqrt [3]{3}+c\sqrt{2},\]
and $c=1$ will do the job. Hence
\[
\mathbb{Q}(\sqrt{2}, \sqrt [3]{3})=\mathbb{Q}(\sqrt{2}+\sqrt [3]{3}).\]

\vspace{0.4cm}

Exercise: what is $[\mathbb{Q}(\sqrt{2}, \sqrt [3]{3}): \mathbb{Q})]$?

\chapter{Trace, Norm, and Discriminant}

\vspace{0.4cm}
\textbf{Definition}. If $K$ is an extension of the complex number field $F$ then an \emph{embedding of K over F} is a ring monomorphism $\sigma: K\rightarrow \textbf{C}$ which fixes each element of $F$, i.e., $\sigma(\alpha)=\alpha$, for all $\alpha \in F$. 
 
\vspace{0.4cm}
\emph{Notation}. If $n$ is a positive integer then we let
\[
F_n[x]=\textrm{the set of all polynomials in $F[x]$ of degree $\leq n-1$}.\]

 \vspace{0.4cm}
 \textbf{Proposition 14}. \textit{If K is a finite extension of F and $n=[K:F]$ then there are n distinct embeddings of K over F. Moreover, if $\theta$ is a primitive element of K over F, $\theta_1,\dots,\theta_n$ are the conjugates of $\theta$ over F, and if $\sigma_i:K \rightarrow \mathbf{C}$ is defined by
 \[
 \sigma_i(\alpha)=q(\theta_i),\ \alpha\in K,\ i=1,\dots,n,\]
 where $q(x)$ is the unique polynomial in $F_n[x]$ such that $\alpha=q(\theta)$, then $\sigma_1,\dots,\sigma_n$ are the embeddings of K over F}.

 \vspace{0.4cm}
 \emph{Proof}. We prove first that each $\sigma_i$ is an embedding of $K$ over $F$. Let $\alpha, \beta \in K$. Because $\theta$ is a primitive element of $K$ over $F$, there exist polynomials $q, r\in F_n[x]$ such that $\alpha=q(\theta), \beta=r(\theta)$. If $\alpha \in F$ then $q(x)\equiv x$ hence $\sigma_i(\alpha)=q(\alpha)=\alpha$. Also, $\alpha\pm \beta=(q\pm r)(\theta)$ and $q\pm r\in F_n[x]$, hence,
 \[
 \sigma(\alpha\pm \beta)=(q \pm r)(\theta_i)=q(\theta_i)\pm r(\theta_i)=\sigma_i(\alpha)\pm \sigma_i(\beta).\]
Let $s\in F_n[x]$ satisfy $\alpha \beta=s(\theta)$. Then $s-qr\in F[x]$ and $(s-qr)(\theta)=s(\theta)-q(\theta)r(\theta)=0$, so Corollary 2 implies that the minimal polynomial $p(x)$ of  $\theta$ over $F$ divides $s(x)-q(x)r(x)$. But $p(\theta_i)=0$ for all $i=1,\dots,n$, hence $(s-qr)(\theta_i)=0$, $i=1,\dots,n$, hence
\[
\sigma_i(\alpha \beta)=s(\theta_i)=q(\theta_i)r(\theta_i)=\sigma_i(\alpha)\sigma_i(\beta).\]
To prove that $\sigma_i$ is injective, assume that $0=\sigma_i(\alpha)=q(\theta_i)$. Then $q(x)\equiv 0$; otherwise, $\theta_i$ is the root of a nonzero polynomial over $F$ of degree $<n$, which is impossible because the degree of $\theta_i$ over $F$ is $n$. Hence $\alpha=q(\theta)=0$. Finally, note that because the $\theta_i$'s are distinct, it follows that if $i\not= j$ then $\sigma_i(\theta)=\theta_i\not=\theta_j=\sigma_j(\theta)$, i.e., the maps $\sigma_1,\dots,\sigma_n$ are all distinct.
  
Suppose next that $\sigma:K\rightarrow \textbf{C}$ is an embedding of $K$ over $F$. We must prove that there is an $i\in \{1,\dots\,n\}$ such that $\sigma=\sigma_i$. In order to do that, observe that
\[
p\big(\sigma(\theta)\big)=\sigma\big(p(\theta)\big)=\sigma(0)=0,\]
i.e., $\sigma(\theta)$ is a root of $p(x)$, hence for some $i\in \{1,\dots\,n\}$, $\sigma(\theta)=\theta_i$. Then for $\alpha=q(\theta)\in K, q\in F_n[x]$, we have
\[
\sigma(\alpha)=\sigma \big(q(\theta)\big)=q \big(\sigma(\theta)\big)=q(\theta_i)=\sigma_i(\alpha),\]
and so $\sigma=\sigma_i$. $\hspace{13cm} \textrm{QED}$

\vspace{0.4cm}
\textbf{Definitions}. Let $K$ be a finite extension of $F$ of degree $n$, $\sigma_1,\dots,\sigma_n$ the embeddings of $K$ over $F$. If $\alpha \in K$ then the \emph{trace of $\alpha$ over F} is
\[
T(\alpha)=\sum_i \sigma_i(\alpha),\]
and the \emph{norm of $\alpha$ over F} is
\[
N(\alpha)=\prod_i \sigma_i(\alpha).\]

\vspace{0.4cm}
The definitions of the trace $T$ and norm $N$ imply that
\[
T(\alpha+\beta)=T(\alpha)+T(\beta),\ N(\alpha \beta)=N(\alpha)N(\beta),\ \forall \alpha, \beta \in K.\]
Because $\sigma_i$ fixes each element of $F$, it follows that
\[
T(\alpha \beta)=\alpha T(\beta),\ N(\alpha \beta)=\alpha^nN(\beta),\ \forall \alpha \in F,\ \forall \beta \in K;\]
in particular (taking $\beta=1$),
\[
T(\alpha)=n \alpha,\ N(\alpha)=\alpha^n,\ \forall \alpha \in F.\]

Proposition 14 implies that if $\theta$ is a primitive element of $K$ over $F$, $\theta_1,\dots,\theta_n$ the conjugates of $\theta$ over $F$, $\alpha \in K$ with $\alpha=q(\theta)$ for some $q\in F_n[x]$, then
\[
T(\alpha)=\sum_i q(\theta_i),\ N(\alpha)=\prod_i q(\theta_i).\]
Moreover, these formulas for $T$ and $N$ do not depend on the primitive element used to define them, since a different choice of primitive element simply permutes the numbers $q(\theta_1),\dots,q(\theta_n)$.

\vspace{0.4cm}
\textbf{Proposition 15}. \textit{T and N map K into F.}

\vspace{0.4cm}
\emph{Proof}. Let $\theta$ be a primitive element of $K$ over $F$,\[
p(x)=x^n+\sum_{i=0}^{n-1} a_ix^i\]
the minimal polynomial of $\theta$ over $F$. Consider the companion matrix
\[
M=\left[ \begin{array}{cccccc}
0 & 0 & 0 & \dots & 0 & -a_0\\
1 & 0 & 0 & \dots & 0 & -a_1\\
0 & 1 & 0 & \dots & 0 & -a_2\\
\vdots & \vdots & \vdots & \dots & \vdots & \vdots\\
0 & 0 & 0 & \dots & 1 & -a_{n-1} \end{array} \right] \]

\vspace{0.2cm}
\noindent of $p(x)$. The characteristic polynomial of $M$ is $p(x)$, hence the eigenvalues of $M$ are the conjugates $\theta_1,\dots,\theta_n$ of $\theta$ over $F$. Hence $M$ is an $n\times n$ matrix with $n$ distinct eigenvalues, and so $M$ is similar to the diagonal matrix

\[
\left[ \begin{array}{cccc}
\theta_1 & 0  & \dots & 0\\
0 & \theta_2 & \dots & 0 \\
\vdots & \vdots & \dots & \vdots\\
0 & 0 & \dots & \theta_n\\
 \end{array} \right] \]

\vspace{0.2cm}

\noindent (Hungerford [8], Theorem VII.5.5, exercise 8, Section VII.5). 

Now let $\alpha \in K$, with $q\in F_n[x]$ such that $\alpha=q(\theta)$. Then

\vspace{0.02cm}
\[
q(M)\ \textrm{is similar to}\ \left[ \begin{array}{cccc}
q(\theta_1) & 0  & \dots & 0\\
0 & q(\theta_2) & \dots & 0 \\
\vdots & \vdots & \dots & \vdots\\
0 & 0 & \dots & q(\theta_n)\\
 \end{array} \right], \]
and so 
\begin{equation*}
T(\alpha)=\sum_i q(\theta_i)=\textrm{trace of}\ q(M),\tag{$*$}\]
\begin{equation*}
N(\alpha)=\prod_i q(\theta_i)=\textrm{ determinant of}\ q(M). \tag{$**$}\]
As all of the entries of $M$ are in $F$, so also are all of the entries of $q(M)$, and so the trace and the determinant of $q(M)$ are in $F$. That $T(\alpha)$ and $N(\alpha)$ are in $F$ is now a consequence of $(*)$ and $(**)$. $\hspace{13cm} \textrm{QED}$

N.B. One can show that if $m=[K: F(\alpha)]$ and if 
\[
x^k+\sum_{i=0}^{k-1} c_ix^i\]
is the minimal polynomial of $\alpha$ over $F$ then
\[
T(\alpha)=-mc_{k-1},\ N(\alpha)=\big((-1)^kc_0\big)^m\]
(Marcus [9], Chapter 2, Theorem 2.4, Corollary 2.1).

N.B. Let $\{\alpha_1,\dots,\alpha_n\}$ be a basis of $K$ over $F$, $K$ considered in the usual way as a vector space over $F$. For $\alpha \in K$, let $[a_{ij}]$ be the $n\times n$ matrix over $F$ defined by
\[
\alpha \alpha_i=\sum_j a_{ij}\alpha_j;\]
$[a_{ij}]$ is the representing matrix with respect to $\{\alpha_1,\dots,\alpha_n\}$ of the linear transformation of $K$ that is defined by multiplication by $\alpha$ on $K$. Then one can prove that
\[
T(\alpha)=\textrm{the trace of $[a_{ij}]$},\ N(\alpha)=\textrm{the determinant of $[a_{ij}]$}\]
(Marcus [9], Chapter 2, exercise 17).

\vspace{0.2cm}

\emph{Notation}: if $M$ is a matrix we will denote the determinant of $M$ by $\det M$.

\vspace{0.4cm}
\textbf{Definitions}. $(i)$ Let $K$ be a finite extension of $K$ over $F$, let $\theta$ be a primitive element of $K$ over $F$, let $\theta_1,\dots\theta_n$ be the conjugates of $\theta$ over $F$, and let $\alpha=q(\theta) \in K$, where $q\in F_n[x]$. The \emph{conjugates of $\alpha$ over $K$} are the complex numbers defined by
\[
q(\theta_i),\ i=1,\dots,n.\]
We will denote the conjugates of $\alpha$ over $K$ by $\alpha^{(1)},\dots,\alpha^{(n)}$, i.e., $\alpha^{(i)}=q(\theta_i),\ i=1,\dots,n.$

$(ii)$ Let $(\alpha_1,\dots,\alpha_n)$ be an $n$-tuple in $K^n$ and let $\alpha_i^{(1)},\dots,\alpha_i^{(n)}$ be the conjugates of $\alpha_i$ over $K$. The \emph{discriminant of }$(\alpha_1,\dots,\alpha_n)$, denoted by $\Delta(\alpha_1,\dots,\alpha_n)$, is

\vspace{0.2cm}
\[
\left( \det \left[ \begin{array}{cccc}
\alpha_1^{(1)} &  \alpha_1^{(2)} & \dots &  \alpha_1^{(n)}\\
\alpha_2^{(1)} &  \alpha_2^{(2)} & \dots &  \alpha_2^{(n)}\\
\vdots & \vdots & \dots & \vdots\\
\alpha_n^{(1)} &  \alpha_n^{(2)} & \dots &  \alpha_n^{(n)}\\
 \end{array} \right] \right)^2.\]

\vspace{0.4cm}
N.B. (1) It is an immediate consequence of Proposition 14 that taking the $i$-th conjugate over $K$ is an $F$-linear homomorphism, i.e., if $a, b\in F$ and $\alpha, \beta \in K$ then $(a\alpha+b\beta)^{(i)}=a\cdot \alpha^{(i)}+b\cdot \beta^{(i)}$ and $(\alpha \beta)^{(i)}=\alpha^{(i)} \beta^{(i)}$.

(2) The value of $\Delta(\alpha_1,\dots,\alpha_n)$ depends on neither the primitive element of $K$ over $F$ used to define it nor on the way the $\alpha_1,\dots,\alpha_n$ are ordered in the the $n$-tuple $(\alpha_1,\dots,\alpha_n)$. A different primitive element simply permutes the columns of the matrix the square of whose determinant defines $\Delta(\alpha_1,\dots,\alpha_n)$, and a different ordering of $\alpha_1,\dots,\alpha_n$ simply interchanges certain rows of that matrix, hence the determinants differ only by a sign, and so the square of their values are the same. We also note that if two coordinates of $(\alpha_1,\dots,\alpha_n)$ are the same then $\Delta(\alpha_1,\dots,\alpha_n)=0$.

\vspace{0.4cm}
\textbf{Proposition 16}. \textit{If $T$ is the trace of $K$ over $F$ and $(\alpha_1,\dots,\alpha_n)$ is an $n$-tuple in $K^n$ then}
\[
\Delta(\alpha_1,\dots,\alpha_n)=\det\ [T(\alpha_i \alpha_j)].\]

\vspace{0.4cm}
\emph{Proof}. If $\sigma_1,\dots,\sigma_n$ are the embeddings of $K$ over $F$ then
\[
\alpha_i^{(j)}=\sigma_j(\alpha_i).\]
If we multiply the matrix $[\sigma_j(\alpha_i)]$ by its transpose $[\sigma_i(\alpha_j)]$ then the $(i, j)$ entry of the product is
\[
\sum_k \sigma_k(\alpha_i) \sigma_k(\alpha_j)=\sum_k \sigma_k(\alpha_i \alpha_j)=T(\alpha_i \alpha_j),\]
and so $[T(\alpha_i \alpha_j)]$ is the matrix product of $[\alpha_i^{(j)}]$ and the transpose of $[\alpha_i^{(j)}]$. It follows that
\[
\det\ [T(\alpha_i \alpha_j)]=\det [\alpha_i^{(j)}]\cdot \det (\textrm{the transpose of}\ [\alpha_i^{(j)}])=\big(\det\  [\alpha_i^{(j)}]\big)^2=\Delta(\alpha_1,\dots,\alpha_n).\]
$\hspace{15.6cm} \textrm{QED}$

Propositions 15 and 16 immediately imply

\vspace{0.4cm}
\textbf{Proposition 17}. \textit{If $(\alpha_1,\dots,\alpha_n) \in K^n$ then $\Delta(\alpha_1,\dots,\alpha_n) \in F$}.

\vspace{0.4cm}
\textbf{Corollary 18}. \textit{The subset $\{\alpha_1,\dots,\alpha_n\}$ of K is a basis of K over F if and only if $\Delta(\alpha_1,\dots,\alpha_n)\not= 0$}.

\vspace{0.4cm}
\emph{Proof}. $(\Rightarrow)$. Suppose that $\Delta(\alpha_1,\dots,\alpha_n)=0$. Then Proposition 16 implies that the system of linear equations
\[
\sum_i x_iT(\alpha_i \alpha_j)=0,\ j=1,\dots,n,\]
has a nontrivial solution $x_i=a_i\in F,\ i=1,\dots,n.$ Let $\alpha=\sum_i a_i\alpha_i \not= 0$. Then $T(\alpha\alpha_j)=0,$ for all $j$. We are assuming that $\{\alpha_1,\dots,\alpha_n\}$ is a basis, and so this implies that $T(\alpha \beta)=0$, for all $\beta \in K$. Taking $\beta=1/\alpha$, we deduce that $n=T(1)=0$, a clear contradiction.

$(\Leftarrow)$. Suppose that $\{\alpha_1,\dots,\alpha_n\}$ is linearly dependent over $F$. Then there is a nonzero $n$-tuple $(a_1,\dots,a_n)\in F^n$ such that $\sum_i a_i\alpha_i = 0$. Multiply this equation by $\alpha_j$ and take the trace to obtain
\[
\sum_i a_iT(\alpha_i \alpha_j)=0,\ j=1,\dots,n.\]
The matrix $[T(\alpha_i \alpha_j)]$ is hence singular, and so, by Proposition 16, 
\[
\Delta(\alpha_1,\dots,\alpha_n)=\det [T(\alpha_i \alpha_j)]=0.\]
$\hspace{15.6cm} \textrm{QED}$

\emph{Notation}: If $\theta$ is a primitive element of $K$ over $F$, we let $D(\theta)$ denote the discriminant $\Delta(1,\dots,\theta^{n-1})$.

\vspace{0.4cm}
\textbf{Proposition 19}. \textit{$(i)$ If $\theta$ is a primitive element of $K$ over $F$ of degree n and $\theta_1,\dots,\theta_n$ are the conjugates of $\theta$ over F, then}
\[
D(\theta)=\prod_{1\leq i<j\leq n} (\theta_i-\theta_j)^2.\]

$(ii)$ \textit{If $(\alpha_1,\dots,\alpha_n)$ and $(\beta_1,\dots,\beta_n)$ are in $F^n$ and if }
\[
\alpha_j=\sum_ic_{ij}\beta_i,\ j=1,\dots,n,\]
\textit{where $c_{ij} \in F$, for all i and j, then}
\[
\Delta(\alpha_1,\dots,\alpha_n)=\big(\det [c_{ij}]\big)^2 \Delta(\beta_1,\dots,\beta_n).\]

\vspace{0.4cm}
\emph{Proof}. $(i)$ The definition of the discriminant implies that

\vspace{0.02cm}
\begin{eqnarray*}
\sqrt{D(\theta)}&=&\det \left[ \begin{array}{cccc}
1 &  1 & \dots &  1\\
\theta^{(1)} &  \theta^{(2)} & \dots &  \theta^{(n)}\\
(\theta^2)^{(1)} &  (\theta^2)^{(2)} & \dots &  (\theta^2)^{(n)}\\
\vdots & \vdots & \dots & \vdots\\
(\theta^{n-1})^{(1)} &  (\theta^{n-1})^{(2)} & \dots &  (\theta^{n-1})^{(n)}\\
 \end{array} \right] \\
&=&\det \left[ \begin{array}{cccc}
1 &  1 & \dots &  1\\
\theta_1 &  \theta_2 & \dots &  \theta_n\\
\theta_1^2 &  \theta_2^2 & \dots &  \theta_n^2\\
\vdots & \vdots & \dots & \vdots\\
\theta_1^{n-1} &  \theta_2^{n-1} & \dots &  \theta_n^{n-1}\\
 \end{array} \right] \\
&=& \textrm{Vandermonde determinant of}\ \theta_1,\dots,\theta_n\\
&=& \prod_{1\leq i<j\leq n} (\theta_i-\theta_j)\ .
\end{eqnarray*}

\vspace{0.2cm}
$(ii)$ Because $c_{ki} \in F$, for all $i$ and $k$, it follows that
\[
\alpha_i^{(j)}=\sum_{k=1}^n c_{ki} \beta_k^{(j)},\]
i.e.,
\[
[\alpha_i^{(j)}]=\big(\textrm{transpose of}\ [c_{ij}] \big) [\beta_i^{(j)}],\]
hence taking the determinant of this equation and squaring yields\[
\Delta(\alpha_1,\dots,\alpha_n)=\big(\det [c_{ij}]\big)^2 \Delta(\beta_1,\dots,\beta_n)\ .\]
$\hspace{15.6cm}\textrm{QED}$
\vspace{0.02cm}
\begin{center}
\textit{Examples}
\end{center}

(1) \emph{Quadratic number fields}

\vspace{0.2cm}
If $m$ is a square-free integer then $\sqrt m$ is a primitive element of $\mathbb{Q}(\sqrt m)$ and $\sqrt m$, $-\sqrt m$ are the conjugates of $\sqrt m$ over $\mathbb{Q}$. Hence the discriminant $D(\sqrt m)$ is

\vspace{0.2cm}
\[
\left( \det \left[ \begin{array}{cc}
1 &  1\\
\sqrt m & -\sqrt m\\
 \end{array} \right] \right)^2=4m\ .\]

\vspace{0.2cm}
(2) \emph{Cyclotomic number fields}

\vspace{0.4cm}
\textbf{Proposition 20}. \textit{If p is an odd prime and $\omega=e^{2\pi i/p}$ then the discriminant $D(\omega)$, taken with respect to $\mathbb{Q}(\omega)$, is}
\[
(-1)^{(p-1)/2}p^{p-2}.\]

\vspace{0.4cm}
\emph{Proof}. $\omega$ is a primitive element of $\mathbb{Q}(\omega)$, of degree $p-1$, with conjugates $\omega,\omega^2,\dots,\omega^{p-1}$ over $\mathbb{Q}(\omega)$ (Lemma 11). Hence Proposition 19$(i)$ implies that
\begin{equation*}
D(\omega)=\prod_{1\leq i<j\leq p-1} (\omega^i-\omega^j)^2. \tag{1}\]

We proceed to calculate the product on the right-hand side of (1). Begin with the factorization
\begin{equation*}
\frac{x^p-1}{x-1}=\sum_{i=0}^{p-1} x^i=\prod_{i=1}^{p-1} (x-\omega^i).\tag{2}\]
Differentiate this equation with respect to $x$ and then set $x=\omega^j$, Since $\omega^p=1$, we obtain
\begin{equation*}
-\frac{p\omega^{p-j}}{1-\omega^j}=\prod_{\substack{i=1 \\ i\not= j}}^{p-1} (\omega^j-\omega^i).\tag{3}\]
Now let $x=1$ and then let $x=0$ in (2) to get, respectively,
\begin{equation*}
p=\prod_{i=1}^{p-1} (1-\omega^i)\ ; \tag{4}\]
\begin{equation*}
\prod_{i=1}^{p-1}\omega^i=(-1)^{p-1}=1. \tag{5}\]
Hence (5) implies that
\begin{equation*}
\prod_{j=1}^{p-1}\omega^{p-j}=\omega^{p-1}\cdots \omega^1=1.\tag{6}\]
(3), (4), (6) imply that
\begin{equation*}
\prod_{j=1}^{p-1}\ \prod_{\substack{i=1 \\ i\not= j}}^{p-1} (\omega^j-\omega^i)=p^{p-1}\frac{\displaystyle {\prod_{j=1}^{p-1}\omega^{p-j}}}{\displaystyle{\prod_{j=1}^{p-1} (1-\omega^j)}}=\frac{p^{p-1}}{p}=p^{p-2}.\tag{7}\]
In the product on the left-hand side of (7), $i<j$ for half of the factors and $j<i$ for the other half. There are $(p-1)(p-2)$ factors in all. Hence the product on the left-hand side of (7) is
\begin{eqnarray*}
p^{p-2}&=&(-1)^{(p-1)(p-2)/2} \prod_{1\leq i<j\leq p-1} (\omega^i-\omega^j)^2\\
&=&(-1)^{(p-1)(p-2)/2} D(\omega),\ \textrm{by} (1)\\
&=&(-1)^{(p-1)/2} D(\omega),\ \textrm{because $p$ is odd},
\end{eqnarray*}
hence
\[
D(\omega)=(-1)^{(p-1)/2}p^{p-2}.\]
$\hspace{15.6cm}\textrm{QED}$

\chapter{Algebraic Integers and Number Rings}

\textbf{Definition}. A complex number $\alpha$ algebraic over $\mathbb{Q}$ is an \emph{algebraic integer} if all coefficients of the minimal polynomial of  $\alpha$ over $\mathbb{Q}$ are integers.

\vspace{0.4cm}
\emph{Notation}. We will let $\mathcal{A}$ denote the set of all algebraic integers.

\vspace{0.4cm}
\textbf{Theorem 21}.\textit{ The set of all algebraic integers $\mathcal{A}$ is a subring of $\mathbf{C}$ which contains the set $\mathbb{Z}$ of all integers}.

\vspace{0.4cm}
\emph{Proof}. That $\mathbb{Z}\subseteq \mathcal{A}$ is clear. Let $\alpha, \beta \in \mathcal{A}$. We must prove that $\alpha \pm \beta$ and $\alpha \beta$ are in $\mathcal{A}$.

Let $\alpha_1=\alpha,\dots,\alpha_m$ and $\beta_1=\beta,\dots,\beta_n$ be the conjugates over $\mathbb{Q}$ of $\alpha$ and $\beta$. Observe first that Proposition 6 and the fact that the minimal polynomials of $\alpha$ and $\beta$ have integer coefficients imply that every elementary symmetric polynomial in $m$ (respectively, $n$) variables evaluated at $\alpha_1,\dots,\alpha_m$ (respectively, $\beta_1,\dots,\beta_n$) is an integer. Observe next that it follows from Weisner [11], Theorem 49.9, that Lemma 8 remains true if the field $F$ in its statement is replaced by $\mathbb{Z}$. The proof of Theorem 5 now applies with the appropriate straightforward modifications to show that the polynomials
\[
\prod_{i, j}\big(x-(\alpha_i \pm \beta_j)\big),\ \prod_{i, j}(x-\alpha_i \beta_j)\]
are monic polynomials with integer coefficients that have  $\alpha \pm \beta$ and $\alpha \beta$ as roots. We now invoke the following lemma to conclude that $\alpha \pm \beta$ and $\alpha \beta$ are in $\mathcal{A}$. $\hspace{3cm}\textrm{QED}$

\vspace{0.4cm}
\textbf{Lemma 22}. \textit{If $\alpha$ is a root of a monic polynomial with integer coefficients then $\alpha$ is an algebraic integer.}

\vspace{0.4cm}
\emph{Proof}. Suppose that $f \in \mathbb{Z}[x]$ is monic and $f(\alpha)=0$. If $p(x)$ is the minimal polynomial of $\alpha$ over $\mathbb{Q}$, we must prove that $p(x) \in \mathbb{Z}[x]$.

Corollary 2 implies that there is a $q \in \mathbb{Q}[x]$ such that $f=pq$. We now find a rational number $a/b$ and polynomials $u, v \in \mathbb{Z}[x]$ such that
\[
f=\frac{a}{b} uv\]
and
\begin{quote}
$u$ (respectively, $v$) is a constant multiple of $p$ (respectively, $q$) and $u$ (respectively, $v$) has all of its coefficients relatively prime, i.e., the greatest common divisor of all of the coefficients is 1.
\end{quote}

\emph{Claim $($Gauss' Lemma$)$}: all coefficients of $uv$ are relatively prime.

Assume this claim for now. Because $bf=auv$, $f$ is monic, and $u,v \in \mathbb{Z}[x]$, it follows that $a$ divides $b$ in $\mathbb{Z}$, say $b=am$ for some $m \in \mathbb{Z}$. Hence
\[
mf=uv.\]
But $f \in \mathbb{Z}[x]$, and so $m$ is a common factor of all of the coefficients of $uv$, and so by the claim, $m=\pm 1$, hence
\[
f=\pm uv.\]
$f$ monic hence implies that the leading coefficient of $u$ is $\pm 1$. But $u$ is a constant multiple of $p$ and $p$ is monic, hence $p=\pm u\in \mathbb{Z}[x]$.

\emph{Proof of the claim}. Suppose that the coefficients of $uv$ have a common prime factor $r$. Let $\mathbb{Z}_r$ denote the field of ordinary residue classes mod $r$. If $s \in \mathbb{Z}[x]$ and if we let $\bar{s}$ denote the polynomial in $\mathbb{Z}_r[x]$ obtained from $s$ by reducing the coefficients of $s$ mod $r$, then $s \rightarrow \bar{s}$ defines a homomorphism of $\mathbb{Z}[x]$ onto $\mathbb{Z}_r[x]$. Because $r$ divides all of the coefficients of $uv$, it hence follows that
\[
0=\overline{uv}=\bar{u}\bar{v}\ \textrm{in}\  \mathbb{Z}_r[x].\]
Because $\mathbb{Z}_r$ is a field, $\mathbb{Z}_r[x]$ is an integral domain, in particular $\mathbb{Z}_r[x]$ has no zero divisors,  hence we conclude from this equation that either $\bar{u}$ or $\bar{v}$ is 0 in $\mathbb{Z}_r[x]$, i.e., either all of the coefficients of $u$ or of $v$ are divisible by $r$. This contradicts the fact that the coefficients of $u$ (respectively, $v$) are relatively prime. $\hspace{10cm}\textrm{QED}$

\vspace{0.4cm}
\textbf{Definitions}. If $F$ is a complex number field then the subring $\mathcal{A}\cap F$ of $F$ is called the \emph{ring of integers in F}. If $F$ is a number field, i.e., a finite extension of $\mathbb{Q}$, then the ring of integers in $F$ is called a \emph{number ring}.

\vspace{0.4cm}
N.B. From now on, when we speak of an \emph{integer} in a complex number field $F$, we will always mean an element of $\mathcal{A}\cap F$, i.e., an algebraic integer in $F$. An element of $\mathbb{Z}$ will be called a \emph{rational} integer, in order to avoid conflict with this terminology.

\vspace{0.2cm}
\begin{center}
\textit{Examples}
\end{center}

(1) \textit {Quadratic number fields}

\vspace{0.4cm}
\textbf{Proposition 23}. $(i)$ \textit{If m is a square-free rational integer and $m\not \equiv 1$ \textnormal{mod 4} then}
\[
\mathcal{A}\cap \mathbb{Q}(\sqrt m)=\{i+j\sqrt m: (i, j) \in \mathbb{Z} \times \mathbb{Z} \}:=\mathbb{Z}+\sqrt m\ \mathbb{Z};\]

$(ii)$ \textit{If $m\equiv 1$ \textnormal{mod 4} then}
\[
\mathcal{A}\cap \mathbb{Q}(\sqrt m)=\left\{i+j\left(\frac{1+\sqrt m}{2} \right): (i, j) \in \mathbb{Z} \times \mathbb{Z} \right\}:=\mathbb{Z}+\left(\frac{1+\sqrt m}{2} \right) \mathbb{Z}.\]
 
\vspace{0.4cm}
\emph{Proof}. Let $R=\mathcal{A}\cap \mathbb{Q}(\sqrt m)$. $\sqrt m$ is an algebraic integer, hence Theorem 21 implies that 
\[
\mathbb{Z}+\sqrt m\ \mathbb{Z} \subseteq R.\]

Let $\alpha \in R\setminus (\mathbb{Z}+\sqrt m\ \mathbb{Z})$. We can write
\[
\alpha=\frac{i+j\sqrt m}{k},\ i, j, k \in \mathbb{Z},\ k>0,\ \gcd(i, j, k)=1.\]

\emph{Claim}: $k=2$.

In order to verify this claim, note first that the degree of $\alpha$ over $\mathbb{Q}$ is either 1 or 2: if the degree of $\alpha=1$ then $\alpha \in \mathbb{Z}$, contrary to its choice. Hence $\alpha$ is the root of a quadratic monic polynomial with coefficients in $\mathbb{Z}$ i.e.,
\[
0=(i+j\sqrt m)^2+bk(i+j\sqrt m)+ck^2,\ \textrm{for some}\ b, c \in \mathbb{Z},\]
i.e.,
\[
0=i^2+j^2m+bki+ck^2+j(2i+bk)\sqrt m,\]
hence
\begin{equation*}
0=i^2+j^2m+bki+ck^2, \tag{1}\]
\begin{equation*}
0=j(2i+bk).\tag{2}\]
Now $j\not=0$; otherwise, $\alpha \in \mathcal{A} \cap \mathbb{Q}=\mathbb{Z}$, \emph{contra hypothesis}, hence by (2),
\begin{equation*}
bk=-2i, \tag{3}\]
and so (1) becomes
\begin{equation*}
j^2m-i^2+ck^2=0.\tag{4}\]

Suppose $p$ is a common prime factor of $i$ and $ k$. (4) implies that $p^2|(j^2m)$ and as $m$ is square-free, we conclude that $p|j$, and this contradicts the fact that  $\gcd(i, j, k)=1$. Hence $\gcd(i, k)=1$. Then by (3), $i|bk$, and so $i|b$ i.e., $b=li$ for some $l \in \mathbb{Z}$. (3) now implies that $lk=-2$, hence $k=1$ or 2. From the choice of $\alpha$, it follows that $k\not= 1$. Hence $k=2$.

We now have that $\alpha=\displaystyle{\frac{i+j\sqrt m}{2}}$ is of degree 2 and is a root of
\[
p(x)=x^2-ix+\frac{i^2-j^2m}{4},\]
hence $p(x)$ is the minimal polynomial of $\alpha$ over $\mathbb{Q}$. Because $\alpha \in R$, it follows that $\displaystyle{\frac{i^2-j^2m}{4}} \in \mathbb{Z}$, hence $i^2\equiv j^2m\ \textrm{mod}\ 4.$ Now $1=\gcd(i, k)=\gcd(i, 2)$, hence $i$ is odd, and so $i^2\equiv 1$ mod 4. Thus
\begin{equation*}
j^2m\equiv 1\ \textrm{mod}\ 4.\tag{5}\]

We assert that $m\equiv 1$ mod 4. In order to see this, observe first that $m$ square-free implies that $m\equiv$ 1, 2, or 3 mod 4. Suppose that $m\equiv 2$ mod 4. Then (5) implies that $2j^2\equiv$ 1 mod 4, hence $2j^2$ is odd, a contradiction. Suppose $m\equiv 3$ mod 4. Then by (5) again, $3j^2 \equiv$ 1 mod 4, and so $j$ is odd, hence $j^2\equiv$ 1 mod 4, whence 3 $\equiv$ 1 mod 4, also not possible. Our assertion follows.

We conclude that if $m\not \equiv$ 1 mod 4 then there are no elements $\alpha \in R\setminus (\mathbb{Z}+\sqrt m\ \mathbb{Z})$, i.e., $R=\mathbb{Z}+\sqrt m\ \mathbb{Z}$. This proves $(i)$.

The fact that  $m\equiv 1$ mod 4 together with (5) implies that $j^2\equiv 1$ mod 4, and so $j$ must be odd. Hence
\[
\alpha=\frac{i+j\sqrt m}{2},\ \textrm{with $i$ and $j$ both odd}.\]
 We have shown that whenever $m\equiv 1$ mod 4,
 \[
R\setminus (\mathbb{Z}+\sqrt m\ \mathbb{Z})\subseteq \left\{\frac{i+j\sqrt m}{2}: (i, j)\in \mathbb{Z} \times \mathbb{Z}, \textrm{ $i$ and $j$ both odd} \right\}.\]
On the other hand, every element $\displaystyle{\frac{i+j\sqrt m}{2}}$ of the set on the right-hand side of this inclusion is in $\mathbb{Q}(\sqrt m)$ and, as per our previous reasoning, is of degree 2 over $\mathbb{Q}$ and has $p(x)$ above as its minimal polynomial over $\mathbb{Q}$. Because $i$ and $j$ are odd and $m\equiv 1$ mod 4, we have that $i^2\equiv j^2m$ mod 4, and so $\displaystyle{\frac{i^2-j^2m}{4}} \in \mathbb{Z}$. Hence $p(x)\in \mathbb{Z}[x]$, whence $\displaystyle{\frac{i+j\sqrt m}{2}}$ is an algebraic integer, and so is in $R\setminus (\mathbb{Z}+\sqrt m\ \mathbb{Z})$. We conclude that if  $m\equiv 1$ mod 4 then
\begin{eqnarray*}
R&=&\mathbb{Z}+\sqrt m\ \mathbb{Z} \cup  \left\{\frac{i+j\sqrt m}{2}: (i, j)\in \mathbb{Z} \times \mathbb{Z}, \textrm{ $i$ and $j$ both odd} \right\}\\
&=&\mathbb{Z}+\left(\frac{1+\sqrt m}{2} \right) \mathbb{Z}.
\end{eqnarray*}
This proves $(ii)$. $\hspace{12.5cm}\textrm{QED}$

In particular, setting $m=-1$ in Proposition 23, we deduce
 
\vspace{0.4cm}
\textbf{Corollary 24}. \textit{If $i=\sqrt{-1}$ then $\mathcal{A} \cap \mathbb{Q}(i)=\mathbb{Z}+i\mathbb{Z}$.}

\vspace{0.4cm}
\textbf{Definition}. $\mathbb{Z}+i\mathbb{Z}$ is called the set of \emph{Gaussian integers}.

\vspace{0.4cm}
(2) \emph{Cyclotomic number fields}

\vspace{0.4cm}
Let $p$ be an odd rational prime, $\omega=e^{2\pi i/p}$.

\vspace{0.4cm}
\textbf{Proposition 25}.
\[
\mathcal{A} \cap \mathbb{Q}(\omega)=\left\{\sum_{k=0}^{p-2} a_k\omega^k: (a_0,\dots,a_{p-2}) \in \mathbb{Z}^{p-1} \right\}.\]

\vspace{0.4cm}
The proof of Proposition 25 lies somewhat deeper than the proof of Proposition 23, requiring some additional mathematical tools for the requisite analysis of the algebraic integers that are involved in the argument. We begin the acquisition of those tools in the remaining part of this chapter and continue doing that in Chapter 7. As we will eventually see, Proposition 25 will be an immediate consequence of Proposition 33, proved in Chapter 7.

\vspace{0.4cm}
\textbf{Proposition 26}. \textit{ Let F be a number field, $R=\mathcal{A}\cap F,\ n=[F:\mathbb{Q}]$.}

$(i)$ \textit{If T and N are the trace and norm mappings of F over $\mathbb{Q}$ then $T(R) \subseteq \mathbb{Z}$ and $N(R) \subseteq \mathbb{Z}$;}

$(ii)$ \textit{ if $(a_1,\dots,a_n) \in R^n$ then $\Delta(a_1,\dots,a_n) \in \mathbb{Z}$}.

\vspace{0.4cm}
\emph{Proof}. Proposition 16 and Theorem 21 imply that $(ii)$ is a consequence of $(i)$, so we need only verify $(i)$, and to do that, it suffices by Proposition 15 to prove that $T(R) \subseteq \mathcal{A}$ and $N(R) \subseteq \mathcal{A}$.

Toward that end, let $\alpha \in F,\ m=$ the degree of $\alpha$ over $\mathbb{Q}$, and let $\alpha_1,\dots, \alpha_m$ and $\alpha^{(1)},\dots, \alpha^{(n)}$ be, respectively, the conjugates of $\alpha$ over $\mathbb{Q}$ and $F$.

\emph{Claim} 1: $m|n$ (in $\mathbb{Z}$) and the list of numbers $\alpha^{(1)},\dots, \alpha^{(n)}$ is obtained, not necessarily in the indicated order, by repeating each number $\alpha_i$ $n/m$ times.

Assume this for now. Then for $\alpha \in R$, we have that $\{\alpha_1,\dots,\alpha_m\}\subseteq \mathcal{A}$, and so Claim 1 and Theorem 21 imply that
\[
T(\alpha)=\sum_1^n \alpha^{(i)}=\frac{n}{m} \sum_1^m \alpha_i \in \mathcal{A},\]
\[
N(\alpha)=\prod_1^n \alpha^{(i)}=\prod_1^m \alpha_i^{n/m}  \in \mathcal{A}.\]

\emph{Proof of Claim} 1. Let $\theta$ be a primitive element of $F$ over $\mathbb{Q}$, $\theta_1,\dots,\theta_n$ the conjugates of $\theta$ over $\mathbb{Q}$. Let $\alpha \in F$ and let $r(\theta_1),\dots,r(\theta_n)$ be the conjugates of $\alpha$ over $F$ (recall that $r(x)$ is the polynomial in $\mathbb{Q}_n[x]$ such that $\alpha=r(\theta)$). Consider the polynomial
\[
f(x)=\prod_1^n \big(x-r(\theta_i)\big);\]
$f(x)$ is called the \emph{field polynomial of $\alpha$ over F}.

\emph{Claim} 2: $f \in \mathbb{Q}[x]$.

To see this, note first that Proposition 6 implies that the coefficients of $f$ are 
\[
(-1)^i\sigma_i\big(r(\theta_1),\dots,r(\theta_n)\big),\] 
where the $\sigma_i$'s are the elementary symmetric polynomials in $n$ variables. As $r \in \mathbb{Q}_n[x]$, each of these coefficients is a symmetric polynomial in $\theta_1,\dots,\theta_n$ over $\mathbb{Q}$. Now apply the following lemma, an immediate consequence of Weisner [11], Theorem 49.10:

\vspace{0.4cm}
\textbf{Lemma 27}. \textit{Every symmetric polynomial in $n$ variables over $\mathbb{Q}$ is a polynomial in the elementary symmetric polynomials $\sigma_1,\dots,\sigma_n$ over  $\mathbb{Q}$}.

\vspace{0.4cm}
Hence there exits a $p_i\in \mathbb{Q}[x_1,\dots,x_n]$ such that
\[
\sigma_i\big(r(\theta_1),\dots,r(\theta_n)\big)=p_i\big(\sigma_1(\theta_1,\dots,\theta_n),\dots,\sigma_n(\theta_1,\dots,\theta_n)\big),\ \textrm{for all}\ i.\]
But $\sigma_i(\theta_1,\dots,\theta_n)\in \mathbb{Q}$, for all $i$ (Corollary 7), hence
\[
(-1)^i\sigma_i\big(r(\theta_1),\dots,r(\theta_n)\big)\in \mathbb{Q},\ \textrm{for all}\ i.\]
This verifies Claim 2.

We have that $f(\alpha)=0$, and so if $p$ is the minimal polynomial of $\alpha$ over $\mathbb{Q}$ then Claim 2 implies that for some nonnegative rational integer $s$,
\begin{equation*}
f=hp^s,\ h\in \mathbb{Q}[x],\ \textrm{with $h$ and $p$ relatively prime in $\mathbb{Q}[x]$}.\tag{6}\]

\emph{Claim} 3: $h(x)\equiv 1$.

Claim 3 and (6) imply that $f=p^s$, i.e., the field polynomial is a power of the minimal polynomial, and Claim 1 is an immediate consequence of that.

\emph{Proof of Claim} 3. Note that $h$ constant implies that $h(x)\equiv 1$ because $f$ and $p$ are both monic, so we need only show that $h$ is constant.

Suppose $h$ is not constant. Then $h$ has one of the $r(\theta_i)$ as a root. Hence $h\circ r$  has $\theta_i$ as a root, hence if $q$ is the minimal polynomial of $\theta_i=$ the minimal polynomial of $\theta$, then $q$ divides $h\circ r$ , and so $0=h\big(r(\theta)\big)=h(\alpha)$. But $h$ and $p$ are relatively prime in $\mathbb{Q}[x]$, hence there exist $u, v\in \mathbb{Q}[x]$ such that $1=hu+pv$, and evaluating this equation at $x=\alpha$ gives
\[
1=h(\alpha)u(\alpha)+p(\alpha)v(\alpha)=0.\]
$\hspace{15.6cm} \textrm{QED}$

\vspace{0.4cm}
\textbf{Lemma} 28. \textit{If $\alpha$ is algebraic over $\mathbb{Q}$ then there is a nonzero $n\in \mathbb{Z}$ such that $n\alpha$ is an algebraic integer.}

\vspace{0.4cm}
\emph{Proof}. Exercise. $\hspace {12cm}\textrm{QED}$

\vspace{0.4cm}
If $F$ is a number field and $\theta$ is a primitive element of $F$ over $\mathbb{Q}$ then Lemma 28 implies that we can always assume that $\theta \in \mathcal{A}$.

\chapter{Integral Bases}

\textbf{Definition}. If $F$ is a number field then the set $\{\alpha_1,\dots,\alpha_s\}$ of integers in $F$ is an \emph{integral basis of F} if every integer $\alpha$ in $F$ can be expressed \emph{uniquely} in the form
\[
\alpha=\sum_{i=1}^s a_i \alpha_i,\ \textrm{where}\ (a_1,\dots,a_s) \in \mathbb{Z}^s.\]

\vspace{0.4cm}
N. B. The concept of an integral basis is the key idea that will be used to prove Proposition 25.

\vspace{0.4cm}
\textbf{Proposition 29}. \textit{Every integral basis of F is a vector-space basis of F over $\mathbb{Q}$}.

\vspace{0.4cm}
\emph{Proof}. Let $\{\alpha_1,\dots,\alpha_s\}$ be an integral basis of $F$, $\alpha \in F$. As $\alpha$ is algebraic over $\mathbb{Q}$, Lemma 28 implies that there exists $0 \not=n\in \mathbb{Z}$ such that $n\alpha$ is an integer in $F$. Hence there exits $ (a_1,\dots,a_s) \in \mathbb{Z}^s$ such that $n\alpha=\sum_{i=1}^s a_i \alpha_i$, and so
\[
\alpha=\sum_i \frac{a_i}{n}\ \alpha_i \in \textrm{span of $\{\alpha_1,\dots,\alpha_s\}$ over $\mathbb{Q}$}.\]

Suppose that $0=\sum_i q_i \alpha_i$ for some $(q_1,\dots,q_s)\in \mathbb{Q}^s$. Multiply this equation by the greatest common denominator $D$ of the $q_i$'s to obtain 
\[
0=\sum_i a_i \alpha_i,\ (a_1,\dots,a_s) \in \mathbb{Z}^s,\]
hence by the uniqueness of this representation, $a_i=0$ for all $i$, and so $q_i=a_i/D=0$, for all $i$.$\hspace{15.4cm}\textrm{QED}$

\vspace{0.4cm}
\textbf{Corollary 30}. \textit{Every integral basis of F has cardinality $[F:\mathbb{Q}]$}.

\vspace{0.4cm}
\textbf{Theorem 31}. \textit{Every number field F has an integral basis}.

\vspace{0.4cm}
\emph{Proof}. Let $\theta$ be an integer in $F$ which is a primitive element of $F$ over $\mathbb{Q}$, and let $n=[F:\mathbb{Q}]$. Consider all of the bases $\{\omega_1,\dots,\omega_n\}$ of $F$ over $\mathbb{Q}$ all of whose elements are integers in $F$ ($\{1,\dots,\theta^{n-1}\}$ is one such basis). Corollary 18 and Proposition 26$(ii)$ imply that
\[
|\Delta(\omega_1,\dots,\omega_n)|\]
is a positive rational integer. Let  $\{\omega_1,\dots,\omega_n\}$  be a basis of this type for which the value $d$ of $|\Delta(\omega_1,\dots,\omega_n)|$ is a minimum.

\emph{Claim}: $\{\omega_1,\dots,\omega_n\}$ is an integral basis of $F$. 

Suppose not. Since $\{\omega_1,\dots,\omega_n\}$ is a basis, there is an integer $\omega$ in $F$ such that $\omega=\sum_iq_i\omega_i$, for some $(q_1,\dots,q_s) \in \mathbb{Q}^s$, but not all $q_i$'s are in $\mathbb{Z}$. We may assume that $q_1 \not \in \mathbb{Z}$. Write
\[
q_1=a+r,\ a\in \mathbb{Z},\ 0<r<1.\]
Define
\[
\omega_1^{*}=\omega-a\omega_1=(q_1-a)\omega_1+\sum_2^n q_i\omega_i,\]
\[
\omega_i^{*}=\omega_i,\ i=2,\dots,n.\]
(We assume here that $n>1$: if $n=1$ then $F=\mathbb{Q}$ and an integral basis in this case is clearly the singleton set $\{1\}$.) Then
\[
\det \left[ \begin{array}{cccc}
q_1-a &  q_2 & \dots &  q_n\\
0 &  1 & \dots &  0\\
0 &  0 & \dots &  0\\
\vdots & \vdots & \dots & \vdots\\
0 &  0 & \dots &  1\\
 \end{array} \right]=q_1-a=r \not= 0, \]
 hence $\{\omega_1^{*},\dots,\omega_n^{*}\}$ is a basis of $F$ over $\mathbb{Q}$ consisting of integers in $F$. Proposition 19$(ii)$ implies that
 \[
|\Delta(\omega_1^{*},\dots,\omega_n^{*})|=r^2|\Delta(\omega_1,\dots,\omega_n)|< |\Delta(\omega_1,\dots,\omega_n)|=d,
\]
contrary to the minimality of $d$.$\hspace{10cm}\textrm{QED}$

\vspace{0.4cm}
\textbf{Proposition 32}. \textit{All integral bases of F have the same discriminant}.

\vspace{0.4cm}
\emph{Proof}. Let $\{\alpha_1,\dots,\alpha_n\}, \{\beta_1,\dots,\beta_n\}$ be integral bases of $F$. 
Then
\[
\alpha_j=\sum_i c_{ij} \beta_i,\ \textrm{where $c_{ij} \in \mathbb{Z}$ for all $i$ and $j$}.\]
Proposition 19$(ii)$ implies that
\begin{equation*}
\Delta(\alpha_1,\dots,\alpha_n)=\big(\det [c_{ij}]\big)^2\Delta(\beta_1,\dots,\beta_n).\tag{$*$}\]
The left-hand side of $(*)$ and the factors of the right-hand side of $(*)$ are in $\mathbb{Z}$, hence
\[
\Delta(\beta_1,\dots,\beta_n)\ \textrm{divides $\Delta(\alpha_1,\dots,\alpha_n)$ in $\mathbb{Z}$}.\]
Similarly,
\[
\Delta(\alpha_1,\dots,\alpha_n)\ \textrm{divides $\Delta(\beta_1,\dots,\beta_n)$ in $\mathbb{Z}$}.\]
Hence
\[
\Delta(\alpha_1,\dots,\alpha_n)=\pm \Delta(\beta_1,\dots,\beta_n),\]
and so (*) implies that the $+$ sign must hold here.$\hspace{7cm}\textrm{QED}$

N.B. The proof of Theorem 31 and Proposition 32 imply that the integral bases of $F$ are precisely the $\mathbb{Q}$-bases of $F$ in $\mathcal{A} \cap F$ the absolute value of whose discriminant is minimal. For this reason, an integral basis is sometimes referred to as a \emph{minimal basis of F}.

\vspace{0.4cm}
\textbf{Definition}. The common value of the discriminants of the integral bases of $F$ is called the \emph{discriminant of F}.

\vspace{0.4cm}
\begin{center}
\textit{Examples}
\end{center}

\vspace{0.4cm}
(1) \textit{Quadratic Number Fields}

\vspace{0.4cm}
Let $m$ be a square-free rational integer. Proposition 23 implies that if $m\not \equiv 1$ mod 4 (respectively, $m\equiv 1$ mod 4) then $\{1, \sqrt m\}$ (respectively, $\{1, \frac{1}{2}(1+\sqrt m)\}$) is an integral basis of $\mathbb{Q}(\sqrt m)$. Hence the discriminant of $\mathbb{Q}(\sqrt m)$ is

\[
\left( \det \left[ \begin{array}{cc}
1 &  1\\
\sqrt m & -\sqrt m\\
 \end{array} \right] \right)^2=4m,\ \textrm{if $m\not \equiv 1$ mod 4},\]

or 
\vspace{0.4cm} 
 \[
 \left( \det \left[ \begin{array}{cc}
1 &  1\\
\displaystyle{\frac{1+\sqrt m}{2}} & \displaystyle{\frac{1-\sqrt m}{2}}\\
 \end{array} \right] \right)^2=m,\ \textrm{if $m\equiv 1$ mod 4}.\]

\vspace{0.4cm}
(2) \textit{Cyclotomic Number Fields}

\vspace{0.4cm}
Let $p$ be an odd rational prime, $\omega=e^{2\pi i/p}.$

\vspace{0.4cm}
\textbf{Proposition 33}. \textit{$\{1, \omega,\dots,\omega^{p-2}\}$ is an integral basis of $\mathbb{Q}(\omega)$}.

\vspace{0.4cm}
N.B. Proposition 33 implies Proposition 25.

\vspace{0.4cm}
Let $\lambda=1-\omega$. The proof of Proposition 33 requires

\vspace{0.4cm}
\textbf{Lemma 34}. \textit{$\{1, \lambda,\dots,\lambda^{p-2}\}$ is an integral basis of $\mathbb{Q}(\omega)$}.

\vspace{0.4cm}
\emph{Proof} (of Lemma 34). Let $\{\alpha_1,\dots,\alpha_{p-1}\}$ be an integral basis of $\mathbb{Q}(\omega)$. 
Then
\begin{equation*}
\lambda^j=\sum_{i=1}^{p-1} c_{ij} \alpha_i,\ j=0, 1,\dots,p-2,\ \textrm{where $c_{ij} \in \mathbb{Z}$ for all $i$ and $j$}.\tag{1}\]
Proposition 19$(ii)$ implies that
\begin{equation*}
\Delta(1, \lambda,\dots,\lambda^{p-2})=\big(\det [c_{ij}]\big)^2\Delta(\alpha_1,\dots,\alpha_{p-1}).\tag{2}\]
From the equations
\begin{eqnarray*}
\lambda&=&1-\omega\\
\lambda^2&=&1-2\omega+\omega^2\\
\lambda^3&=&1-3\omega+3\omega^2-\omega^3,
\end{eqnarray*}
$\hspace{6.7cm} \vdots$
\begin{eqnarray*}
\omega&=&1-\lambda\\
\omega^2&=&1-2\lambda+\lambda^2\\
\omega^3&=&1-3\lambda+3\lambda^2-\lambda^3,
\end{eqnarray*}
$\hspace{6.7cm} \vdots$

\noindent we deduce that
\[
\Delta(1, \lambda,\dots,\lambda^{p-2})=\big(\det [a_{ij}]\big)^2\Delta(1, \omega,\dots,\omega^{p-2}),\]
\[
\Delta(1, \omega,\dots,\omega^{p-2})=\big(\det [a_{ij}]\big)^2\Delta(1, \lambda,\dots,\lambda^{p-2}),\]
where $[a_{ij}]$ is the matrix of binomial coefficients which occurs on the right-hand sides of the above expansions of the powers of $\lambda$ and $\omega$. Since $\Delta(1, \omega,\dots,\omega^{p-2})\not= 0$, these equations imply that $\big(\det [a_{ij}]\big)^2=1$, hence
\[
\Delta(1, \omega,\dots,\omega^{p-2})=\Delta(1, \lambda,\dots,\lambda^{p-2}),\]
and so (2) and Proposition 20 imply that
\begin{equation*}
\pm p^{p-2}=\big(\det [c_{ij}]\big)^2\Delta(\alpha_1,\dots,\alpha_{p-1}).\tag{3}\]
But $\Delta(\alpha_1,\dots,\alpha_{p-1})$ and the $c_{ij}$'s are all in $\mathbb{Z}$, hence (3) implies that
\begin{equation*}
\det [c_{ij}]=\pm p^r,\ \textrm{for some nonnegative $r \in \mathbb{Z}$.}\tag{4}\]
If we now solve the equations (1) for the $\alpha_i$'s using Cramer's rule, we find from (4) that each $\alpha_i$ has the form
\[
p^{-r}\sum_{i=0}^{p-2} a_i \lambda^i,\textrm{for some $(a_0,\dots,a_{p-2})\in \mathbb{Z}^{p-1}$.}\]
Because $\{\alpha_1,\dots,\alpha_{p-1}\}$ is an integral basis of $\mathbb{Q}(\omega)$, it follows that \emph{every} integer in $\mathbb{Q}(\omega)$ is of this form.

We can now prove that $\{1, \lambda,\dots,\lambda^{p-2}\}$ is an integral basis of $\mathbb{Q}(\omega)$. The argument splits naturally into the two cases $r=0$ and $r\geq 1$.

Suppose first that $r=0$. Then
\[
\mathcal{A} \cap \mathbb{Q}(\omega)=\sum_{i=0}^{p-2} \lambda^i \mathbb{Z}.\]
Since $\Delta(1, \lambda,\dots,\lambda^{p-2})=\Delta(1, \omega,\dots,\omega^{p-2})\not=0$, Corollary 18 implies that $\{1, \lambda,\dots,\lambda^{p-2}\}$ is linearly independent over $\mathbb{Q}$, hence $\{1, \lambda,\dots,\lambda^{p-2}\}$ is an integral basis.

Assume next that $r\geq1$. Suppose that $\{1, \lambda,\dots,\lambda^{p-2}\}$ is not an integral basis of $\mathbb{Q}(\omega)$. Then there is an integer in $\mathbb{Q}(\omega)$ of the form
\[
p^{-r}\sum_{i=0}^{p-2} a_i \lambda^i,\textrm{for some $(a_0,\dots,a_{p-2})\in \mathbb{Z}^{p-1}$}\]
such that $p^r$ does not divide every coefficient $a_i$ in $\mathbb{Z}$. If $p^s$ is the highest power of $p$ that divides all of the $a_i$'s then $s\leq r-1$, hence
\[
\frac{1}{p}\sum_{i=0}^{p-2} b_i \lambda^i=\frac{1}{p}\sum_{i=0}^{p-2} \frac{a_i}{p^s} \lambda^i=p^{r-s-1}\frac{\sum_{i=0}^{p-2} a_i \lambda^i}{p^r}\]
is an integer in $\mathbb{Q}(\omega)$ with $b_i\in \mathbb{Z}$, for all $i$, and $p$ does not divide every coefficient $b_i$. Let 
\[
b_m=\textrm{ the $b_i$ with the smallest $i$ such that $p$ does not divide $b_i$.}\]
Then
\[
\frac{1}{p}\sum_{i=m}^{p-2} b_i \lambda^i \in \mathcal{A} \cap \mathbb{Q}(\omega).\]
From the proof of Proposition 20,
\begin{eqnarray*}
p&=&(1-\omega)(1-\omega^2)\cdots(1-\omega^{p-1})\\
&=&(1-\omega)(1-\omega)\cdots(1-\omega)\alpha\\
&=&\lambda^{p-1}\alpha\\
&=&\lambda^{m+1}\alpha^{\prime},\ \textrm{where $\alpha$ and $\alpha^{\prime}$ are integers in $\mathbb{Q}(\omega)$},
\end{eqnarray*}
hence 
\[
\frac{1}{\lambda^{m+1}}\sum_{i=m}^{p-2} b_i \lambda^i \in \mathcal{A} \cap \mathbb{Q}(\omega).\]
Now $\lambda^{m+1}$ cancels into all of the terms of this sum save for the first one, hence
\begin{equation*}
\frac{b_m}{\lambda} \in \mathcal{A} \cap \mathbb{Q}(\omega).\tag{5}
\end{equation*}

We will prove next that $b_m/\lambda \not \in \mathcal{A}$. This contradicts (5) and will finish the proof of Lemma 34. Set 
\[
b_m=b,\]

\[
x=\frac{b}{\lambda}=\frac{b}{1-\omega}.\]
Then
\[
\omega=1-\frac{b}{x},\]
and so
\[
1=\left(1-\frac{b}{x} \right)^p,\ \textrm{i.e., $x^p=(x-b)^p$},\]
hence $b_m/\lambda$ is a root of
\[
g(x)=px^{p-1}+p(\dots)-b^p,\ \textrm{where $p$ does not divide $b$}\]
(We have used here the fact that $p$ divides each binomial coefficient $\Big( \begin{array}{c}
p\\
k\\ 
 \end{array} \Big),\ k=1,\dots,p-1$).

\emph{Claim}: $g(x)$ is irreducible over $\mathbb{Q}$.

By virtue of this claim, $g(x)/p$ is the minimal polynomial of $b/\lambda$ over $\mathbb{Q}$; since the constant term of  $g(x)/p$ is $b^p/p$ and $p$ does not divide $b$,  $g(x)/p \not \in \mathbb{Z}[x]$, hence  $b/\lambda \not \in \mathcal{A}$.

\emph{Proof of the claim}. We let
\[
f(x)=x^{p-1}g\left(\frac{1}{x}\right)=b^px^{p-1}+p(\dots)+p.\]
Now apply to $f(x)$

\vspace{0.4cm}
\textbf{Proposition 35} (\textit{Eisenstein's irreducibility criterion}). \textit{If $p$ is a rational prime,} 
\[
h(x)=\sum_{i=0}^n a_ix^i \in \mathbb{Z}[x],\]
\textit{$p$ does not divide $a_n$, $p^2$ does not divide $a_0$,and $p$ divides $a_i,\ i=1,\dots,n-1$, then $h(x)$ is irreducible over $\mathbb{Q}$}.

\vspace{0.4cm}
\emph{Proof}. See Hungerford [8], Theorem III.6.15.$\hspace{7cm}\textrm{QED}$

Proposition 35 implies that $f(x)$ is irreducible over  $\mathbb{Q}$, hence $g(x)$ is irreducible over  $\mathbb{Q}$ (if $g=uv$ is a nontrivial factorization of $g$ over $\mathbb{Q}$ then
\[
f(x)=x^{\mu}u\left(\frac{1}{x}\right)x^{\nu}v\left(\frac{1}{x}\right),\ \textrm{where $\mu=$ degree of $u$ and $\nu=$ degree of $v$,}\]
is a nontrivial factorization of $f$ over $\mathbb{Q}$). $\hspace{8.5cm}\textrm{QED}$

\vspace{0.4cm}
\emph{Proof of Proposition} 33. Lemma 34 and the equations
\begin{eqnarray*}
\lambda&=&1-\omega\\
\lambda^2&=&1-2\omega+\omega^2\\
\lambda^3&=&1-3\omega+3\omega^2-\omega^3,
\end{eqnarray*}
$\hspace{6.7cm} \vdots$

\noindent imply that
\[
\mathcal{A}\cap \mathbb{Q}(\omega)=\sum_{i=0}^{p-2}\ \omega^i\mathbb{Z}.\]
Because $\{1, \omega,\dots,\omega^{p-2}\}$ is linearly independent over $\mathbb{Q}$, it is hence an integral basis of $\mathbb{Q}(\omega)$.

$\hspace{15cm}\textrm{QED}$

Propositions 20 and 33 imply

\vspace{0.4cm}
\textbf{Proposition 36}. \textit{If $p$ is an odd rational prime and $\omega=e^{2\pi i/p}$ then the discriminant of $\mathbb{Q}(\omega)$ is} $(-1)^{(p-1)/2}p^{p-2}$.

\vspace{0.4cm}
\textbf{Proposition 37}. \textit{If $m\geq 3$ is a rational integer, $\omega_m=e^{2\pi i/m}$, and $\varphi$ is Euler's totient, then $\{1, \omega_m,\dots,\omega_m^{\varphi(m)-1}\}$ is an integral basis of  $\mathbb{Q}(\omega_m)$.}

\vspace{0.4cm}
\emph{Proof}. Marcus [9], Corollary 2, p. 35. $\hspace{8.2cm}\textrm{QED}$

\vspace{0.4cm}
N.B. Although it is true for quadratic and cyclotomic number fields, it is not always the case that a number field of degree $n$ over $\mathbb{Q}$ has an integral basis of the form $\{1, \theta,\dots,\theta^{n-1}\}$ for a primitive element $\theta$ over $\mathbb{Q}$; see Marcus [9], exercise 30, Chapter 2.

\chapter{The Problem of Unique Factorization in a Number Ring}

We saw in the discussion of Fermat's Last Theorem that if $p$ is an odd rational prime and $\omega=e^{2\pi i/p}$ then the question of when $\mathbb{Z}[\omega]$ is a unique factorization domain (UFD) was of interest. Hence we will now study the following question:
\begin{quote}
if $F$ is a number field with ring of integers $R=\mathcal{A}\cap F$, when is $R$ a UFD?
\end{quote}

\vspace{0.4cm}
\textbf{Proposition 38}. \textit{$\alpha \in R$ is a unit of $R$ if and only if the norm $N(\alpha)$ of $\alpha$ is $\pm 1$}.

\vspace{0.4cm}
\emph{Proof}. Let $n=[F:\mathbb{Q}],\ \alpha_1,\dots,\alpha_n$ denote the conjugates of $\alpha$ over $F$. $\alpha$ is a unit of $R$ if and only if $1/\alpha \in R$, hence whenever $\alpha$ is a unit in $R$ then

\[
1=N(1)=N(\alpha)N\left(\frac{1}{\alpha}\right).\]

\vspace{0.2cm}
\noindent Proposition 26$(i)$ implies that $N(\alpha)$ and $N\left(1/\alpha \right)$ are both in $\mathbb{Z}$ and so $N(\alpha)=\pm 1$. Conversely, if
\[
 \pm 1=N(\alpha)=\alpha_1,\cdots,\alpha_n,\]

\vspace{0.2cm}
\noindent then, as $\alpha_i\in \mathcal{A}$ for all $i$, it follows that
\[
\frac{1}{\alpha}=\pm \alpha_2\cdots \alpha_n\in \mathcal{A}.\]

\vspace{0.2cm}
\noindent But $1/\alpha \in F$, since $F$ is a field, hence $1/\alpha\in R$. $\hspace{7cm}\textrm{QED}$

\vspace{0.4cm}
\textbf{Theorem 39}. \textit{Every nonzero nonunit in $R$ is the product of prime elements of $R$}.

\vspace{0.4cm}
\emph{Proof}. If $\alpha \in R\setminus\big(\{0\} \cup U(R)\big)$ is not prime then $\alpha=\beta \gamma$ for some $\beta, \gamma \in  R\setminus\big(\{0\} \cup U(R)\big)$. Now replace $\alpha$ in this factorization procedure by $\beta$ and $\gamma$ in turn, and continue in this way.

\textit{Claim}: this factorization procedure terminates after finitely many steps, and in the final product $\alpha=\gamma_1\cdots\gamma_s$, $\gamma_i$ is prime for all $i$.

Suppose the factorization procedure never terminates. Then for arbitrarily large $n$, $\alpha=\gamma_1\cdots\gamma_n$ and $\gamma_i \not \in \{0\} \cup U(R)$ for all $i$, and so Proposition 38 implies that $N(\gamma_i)$ is a rational integer of absolute value larger than 1, for all $i$. Hence
\[
|N(\alpha)|=\prod_{i=1}^n |N(\gamma_i)|\]
is arbitrarily large, which is obviously impossible. In the final product $\alpha=\gamma_1\cdots\gamma_s$, each $\gamma_i$ is prime; otherwise the factorization procedure will continue beyond its termination point.

$\hspace{15cm}\textrm{QED}$

N.B. Theorem 39 implies that in any number ring, prime factorization always holds. As we will now prove, what can fail is \emph{uniqueness} of the prime factorization.

\vspace{0.4cm}
\emph{Unique prime factorization can fail.}

We show that uniqueness of prime factorization fails in the ring of algebraic integers in $\mathbb{Q}(\sqrt{-5})$. Observe first that by virtue of Proposition 23,
\[
\mathcal{A}\cap \mathbb{Q}(\sqrt{-5})=\mathbb{Z}+\sqrt{-5}\ \mathbb{Z}.\]

\vspace{0.4cm}
\textbf{Proposition 40}. $(i)$ $U(\mathbb{Z}+\sqrt{-5}\ \mathbb{Z})=\{-1, 1\}$.

$(ii)$ 3, 7, and $1\pm 2\sqrt{-5}$ are primes in $\mathbb{Z}+\sqrt{-5}\ \mathbb{Z}.$

\vspace{0.4cm}
Since
\[
21=3\cdot7=(1+2\sqrt{-5})(1-2\sqrt{-5}),\]

\noindent Proposition 40 implies that the prime factorization of 21 in $\mathbb{Z}+\sqrt{-5}\ \mathbb{Z}$ is not unique up to order of the factors and multiplication by units, hence 
$\mathbb{Z}+\sqrt{-5}\ \mathbb{Z}$ is not a UFD.

Set $R=\mathbb{Z}+\sqrt{-5}\ \mathbb{Z}$.

\emph{Proof of Proposition} 40. $(i)$ If $\alpha=a+b\sqrt{-5}\in R$ then
\[
N(\alpha)=(a+b\sqrt{-5})(a-b\sqrt{-5})=a^2+5b^2 \geq 0\]
and so
\[
\alpha \in U(R)\ \textrm{iff}\ N(\alpha)=1\ \textrm{iff}\ a^2+5b^2=1\ \textrm{iff}\ a=\pm1\ \textrm{and $b=0$}.\]

$(ii)$ 3 \emph{is prime in R}.

Suppose $3=\alpha \beta,\  \alpha, \beta \in R\setminus U(R)$. Then $9=N(3)=N(\alpha)N(\beta)$, with $N(\alpha)\not= \pm1 \not= N(\beta)$. Hence $N(\alpha)=N(\beta)=3$ and so for some $a, b \in \mathbb{Z}$, $a^2+5b^2=3$, which is clearly not possible.

7 \emph{is prime in R}.

Arguing as before from $7=\alpha \beta,\  \alpha, \beta \in R\setminus U(R)$, we obtain  $a^2+5b^2=7$ for some $a, b \in \mathbb{Z}$, clearly impossible.

$1\pm 2\sqrt{-5}$ \emph{are primes in R}.

$1\pm 2\sqrt{-5}=\alpha \beta,\  \alpha, \beta \in R\setminus U(R)$ implies that
\[
N(\alpha)N(\beta)=N(1\pm 2\sqrt{-5})=(1+2\sqrt{-5})(1-2\sqrt{-5})=3\cdot7,\] 
with  $N(\alpha)\not= \pm1 \not= N(\beta)$, hence $N(\alpha)=$ 3 or 7, neither of which can be true, as shown above.$\hspace{14.5cm}\textrm{QED}$

\vspace{0.4cm}
\emph{Unique prime factorization can hold.}

\vspace{0.4cm}
Let $i=\sqrt{-1}$. We will prove

\vspace{0.4cm}
\textbf{Proposition 41}. \textit{Unique prime factorization holds in the ring $\mathbb{Z}+i \mathbb{Z}$ of Gaussian integers.}

\vspace{0.4cm}
\emph{Proof}. Let $R=\mathbb{Z}+i \mathbb{Z}$. This will be proven by adapting to $R$ the standard argument which verifies unique prime factorization in $\mathbb{Z}$. Thus, we need the following two lemmas:

\vspace{0.4cm}
\textbf{Lemma 42} \textit{$($Division algorithm in R$)$. If $\alpha, \beta \in R$ with $\beta\not= 0$, then there exists $\pi, \rho \in R$ such that}
\[
\alpha=\pi \beta+\rho,\ \textrm{with}\ N(\rho)<N(\beta).\]

\vspace{0.4cm}
\textbf{Lemma 43} \textit{$($Fundamental property of primes in R$)$. If $\pi$ is a prime in $R$ and $\pi$ divides the product $\alpha \beta$ in R, then either  $\pi$ divides $\alpha$ or $\pi$ divides $\beta$ in R}.

\vspace{0.4cm}
We will assume Lemmas 42 and 43 for now and use them to prove Proposition 41 like so: suppose that
\[
\pi_1\cdots\pi_r=\sigma_1\cdots\sigma_s\]
are products of primes in $R$, with $s\leq r$, say. Lemma 43 implies that $\sigma_1$ must divide one of the $\pi_i$'s, say $\pi_1$. Since $\pi_1$ is prime, there is a unit $u_1$ in $R$ such that $\pi_1=u_1\sigma_1$. Then
\[
\pi_2\cdots\pi_r=u_1\sigma_2\cdots\sigma_s.\]
Suppose that $s<r$. Then a continuation of this procedure stops at a sequence of units $u_1,\dots,u_s$ in $R$ for which
\[
\pi_{s+1}\cdots\pi_r=u_1\cdots u_s.\]
By virtue of Proposition 38, the $\pi_i$'s here are all units, which contradicts the definition of prime element. Hence $s=r$, and upon reindexing of $\pi_1,\dots,\pi_r$, we obtain
$\pi_i=u_i\sigma_i$, with $u_i$ a unit, for all $i$.$\hspace{13.3cm}\textrm{QED}$

\emph{Proof of Lemma} 42. There exist $q_1, q_2 \in \mathbb{Q}$ such that $\alpha/\beta=q_1+q_2i$. Choose $s_1, s_2 \in \mathbb{Z}$ such that 
\[
|q_j-s_j|\leq \frac{1}{2}\ ,\ j=1, 2.\]

Let $\pi=s_1+s_2i$, $\rho=\alpha-\pi\beta$. 
We will show that $N(\rho)<N(\beta)$. To that end, we calculate that
\begin{eqnarray*}
|\rho|&=&|\beta|\sqrt{(q_1-s_1)^2+(q_2-s_2)^2}\\
&\leq&|\beta|\sqrt{\frac{1}{2^2}+\frac{1}{2^2}}\\
&<&|\beta|,
\end{eqnarray*}

\noindent hence $N(\rho)=|\rho|^2<|\beta|^2=N(\beta)$.$\hspace{9.4cm}\textrm{QED}$

\vspace{0.3cm}
\emph{Proof of Lemma} 43. If $\pi$ divides $\alpha$ in $R$, done, so suppose that $\pi$ does not divide $\alpha$ in $R$. Lemma 42 implies that there exist $\delta, \rho \in R$ such that $\alpha=\delta \pi+\rho,\ 0<N(\rho)<N(\pi)$. Set
\[
T=\{\alpha \xi+\pi \eta \not=0: (\xi, \eta)\in R\times R\}.\]
N.B. $\rho=\alpha-\pi \delta \in T$. The norm of each element of $T$ is a positive rational integer, so choose $\gamma=\alpha\xi_0+\pi \eta_0$ of minimum norm in $T$. Then
\begin{equation*}
N(\gamma)\leq N(\rho)<N(\pi).\tag{1}\]

\emph{Claim}: $\gamma$ is a unit.

In order to see this, choose, by Lemma 42, $\theta, \zeta \in R$ such that
\[
\pi=\theta \gamma+\zeta,\ N(\zeta)<N(\gamma).\]
We have that $\zeta=\alpha(-\theta\xi_0)+\pi(1-\theta\eta_0)$. If $\zeta\not=0$ then $\zeta \in T$. But $ N(\zeta)<N(\gamma)$, contrary to the minimality of $N(\gamma)$. Hence $\zeta=0$, and so $\pi=\theta \gamma$. But $\pi$ is prime, hence $\theta$ or $\gamma$ is a unit. If $\theta$ is a unit then $N(\theta)=1$ and so $N(\pi)=N(\theta)N(\gamma)=N(\gamma)$, contrary to (1). Hence $\theta$ is not a unit, and so $\gamma$ is a unit.

Now
\begin{equation*} 
\gamma\beta=\alpha\beta+\pi\beta\eta_0.\tag{2}\]

\noindent Because $\pi$ divides $\alpha\beta$ in $R$, (2) implies that $\pi$ divides $\gamma\beta$ in $R$, i.e., $\gamma\beta=\tau\pi$ for some $\tau \in R$, hence $\beta=\tau\pi/\gamma$. But $\tau/\gamma\in R$ ($1/\gamma \in R$, because $\gamma$ is a unit!), and so $\pi$ divides $\beta$ in $R$. 

$\hspace{15cm}\textrm{QED}$

Propositions 40 and 41 hence raise the following very important question:
\begin{quote}
how can you tell when a number ring has unique prime factorization?
\end{quote}
The answer to this requires that we investigate the structure of \emph{ideals} in a number ring.

\chapter{Ideals in a Number Ring}
\textbf{Definitions}. Let $A$ be a commutative ring with identity. An \emph{ideal of A} is a subring $I$ of $A$ such that $a\in A$, $b\in I$ implies $ab\in I$. $I$ is a \emph{principal ideal of A} if there exists $a\in A$ such that $I=\{ab: b\in A\}$; in this case we say that \emph{I is generated by a} and we denote that by $I=(a)$. $A$ is a \emph{principal-ideal domain} (PID) if $A$ is an integral domain and every ideal of $A$ is principal.  

A basic result of ring theory is

 \vspace{0.4cm}
 \textbf{Proposition 44}. \textit{Every} PID \textit{is a} UFD.

 \vspace{0.4cm}
 \emph{Proof}. Hungerford [8], Theorem III.3.7.$\hspace{7.7cm}\textrm{QED}$

\vspace{0.3cm}
The converse of Proposition 44 is far from true in general, but, remarkably enough, the converse \emph{is} true for number rings.

 \vspace{0.4cm}
\textbf{Theorem 45}. \textit{If F is a number field and $R=\mathcal{A}\cap F$ then R is a} UFD \textit{if and only if R is a} PID.

 \vspace{0.4cm}
By virtue of Theorem 45, we can determine when a number ring $R$ is a UFD by determining when $R$ is a PID. This latter problem is much more tractable than the former one, but it is not easy; as we will see, determining when a number ring is a PID is in general \emph{very difficult}, even in the case of quadratic number fields.

We saw in Chapter 8 that unique factorization can fail to hold in a number ring. We will now develop a theory of ideals in a number ring $R$ which will show that a close analog of unique prime factorization of elements of $R$ in fact holds for \emph{ideals of R}. This structure theory of ideals will also eventually lead to a proof of Theorem 45.

For the remainder of this chapter, let $F$ be a fixed number field of degree $n$ over $\mathbb{Q}$, and let $R=\mathcal{A}\cap F$.

 \vspace{0.01cm}
 \textbf{Definitions}. Let $I$ be an ideal of $R$, $S=\{\alpha_1,\dots,\alpha_s\}\subseteq I$. $I$ is \emph{generated by $S$}, denoted by $I=(\alpha_1,\dots\alpha_s)$, if
 \[
 I=\left\{\sum_i r_i\alpha_i: (r_1,\dots,r_s)\in R^s\right\}.\]
In this case, the elements of $S$ are called the \emph{generators of R}. $S$ is a \emph{basis of I} if every element $x\in I$ can be expressed \emph{uniquely} in the form
\begin{equation*}
x=\sum_i z_i\alpha_i,\ \textrm{for some}\ (z_1,\dots,z_s) \in \mathbb{Z}^s.\tag{*}\]

 \vspace{0.4cm}
 N.B. A generating set of an ideal need not be a basis. The ideal in $\mathbb{Z}+i\mathbb{Z}$ generated by $\{2\}$ is
 \[
 \{2a+2bi: (a, b)\in \mathbb{Z} \times \mathbb{Z}\},\]
 A basis for this ideal is $\{2, 2i\}$, \emph{not} $\{2\}$.

\vspace{0.4cm}
\textbf{Proposition 46}. \textit{Every ideal $I\not=\{0\}$ of $R$ has a basis and every basis of $I$ has cardinality $n=[F:\mathbb{Q}]$}.

\vspace{0.4cm}
\emph{Proof}. Existence of a basis of $I$ follows from a straightforward modification of the proof of Theorem 31. Consider all sets $\{\alpha_1,\dots,\alpha_n\}$ of integers in $I$ which are a vector-space basis of $F$ over $\mathbb{Q}$ (if $\{\beta_1,\dots,\beta_n\}$ is an integral basis of $F$ and $0\not=\alpha \in I$ then $\{\alpha\beta_1,\dots,\alpha\beta_n\}$ is one such set). $|\Delta(\alpha_1,\dots,\alpha_n)|$ is a positive rational integer, and if $\{\alpha_1,\dots,\alpha_n\}$ is chosen with the absolute value of its discriminant minimal, then the proof of Theorem 31 verifies that $\{\alpha_1,\dots,\alpha_n\}$ is a basis of $I$.

Let $\{\alpha_1,\dots,\alpha_s\}$ be a basis of $I$: we must prove that $s=n$. Uniqueness of the representation of the elements of $I$ in the form (*) implies that  $\{\alpha_1,\dots,\alpha_s\}$ is linearly independent over $\mathbb{Q}$ (see the proof of Proposition 29), hence $s\leq n$. Suppose that $s<n$. Let $\{\beta_1,\dots,\beta_n\}$ be an integral basis of $F$ and  let $0\not=\alpha \in I$. Then $\{\alpha\beta_1,\dots,\alpha\beta_n\}$ is linearly independent over $\mathbb{Q}$. Since $\alpha\beta_i\in I$ for all $i$, it follows that
\[
\alpha\beta_j=\sum_{i=1}^n a_{ij} \alpha_i,\ \textrm{with $a_{ij}\in \mathbb{Z}$ for all $i, j$},\]
where we set $\alpha_i=0$ for $s<i\leq n$. Proposition 19$(ii)$ implies that
\[
\Delta(\alpha\beta_1,\dots,\alpha\beta_n)=\big(\det [a_{ij}]\big)^2\Delta(\alpha_1,\dots,\alpha_n).\]
However, $\Delta(\alpha_1,\dots,\alpha_n)=\Delta(\alpha_1,\dots,\alpha_s,0,\dots,0)=0$, hence $\Delta(\alpha\beta_1,\dots,\alpha\beta_n)=0$, which contradicts Corollary 18. Thus $s=n$.$\hspace{9cm}\textrm{QED}$

\vspace{0.4cm}
 \textbf{Definition}. Let $A$ be a commutative ring with identity, $I, J$ ideals of $A$. The \emph{product IJ of I and J} is the ideal of $A$ generated by the set $\{xy: (x, y)\in I\times J\}$, i.e., the smallest ideal of $A$, relative to subset inclusion, which contains this set.

\newpage
 \emph{Easily proven facts}
 
\vspace{0.03cm} 
$(i)$ $IJ$ consists precisely of all sums of the form $\sum_i x_iy_i$, where $x_i\in I$, $y_i\in J$.

$(ii)$ If $I=(x_1,\dots,x_r)$, $J=(y_1,\dots,y_s)$, then
\[
IJ=(x_iy_j: (i, j)\in \{1,\dots,r\}\times \{1,\dots,s\}),\]
and the product does not depend on the choice of generating sets for $I$ and $J$.

$(iii)$ The ideal product is commutative and associative.

\noindent These facts will be used freely and without reference in all of what follows.

\vspace{0.4cm}
 \textbf{Definitions}. If $A$ is a commutative ring with identity then an ideal $P$ of $A$ is \emph{prime} if $\{0\}\not=P\not=A$ and if $a, b \in A$ and $ab \in P$ then $a\in P$ or $b \in P$. An ideal $M$ of $A$ is \emph{maximal} if $\{0\}\not=M\not=A$ and if $I$ is an ideal of $A$ and $M \subseteq I$ then $M=I$ or $I=A$.

\vspace{0.4cm}
If $P$ is a prime ideal of $R$ then we will eventually show (Proposition 54 below) that if $P$ is the product of ideals $I$ and $J$ of $R$ then $\big\{I, J \big\}=\big\{P, R\big\}$, i.e., the only ideal factors of $P$ are $P$ and $R$. This says that the prime ideals of $R$ behave with respect to ideal factorization exactly like the prime elements of $R$ behave with respect to factorization of the elements of $R$.

Another basic result of ring theory is

\vspace{0.4cm}
\textbf{Proposition 47}. \textit{Every maximal ideal in a commutative ring with identity is prime}.

\vspace{0.4cm}
\emph{Proof}. Hungerford [8], Theorem III.2.19.$\hspace{7.8cm}\textrm{QED}$

\vspace{0.4cm}
\textbf{Proposition 48}. \textit{ If P is a prime ideal in a commutative ring $A$ with identity and $I$ and $J$ are ideals of $A$ such that $IJ \subseteq P$ then $I\subseteq P$ or $J\subseteq P$}.

\vspace{0.4cm}
\emph{Proof}. Exercise.$\hspace{12.1cm}\textrm{QED}$

\vspace{0.02cm}
The converse of Proposition 47 is false in general, but the converse is true for number rings. This is a very important fact!

\vspace{0.4cm}
\textbf{Proposition 49}. \textit{An ideal I of R is prime if and only if it is maximal}.

\vspace{0.4cm}
\emph{Proof}. Let $I$ be a prime ideal of $R$: we need to prove that $I$ is a maximal ideal, i.e., we take an ideal $J$ of $R$ which properly contains $I$ and show that $J=R$.

Toward that end, let $\{\alpha_1,\dots,\alpha_n\}$ be an integral basis of $F$, and let $0 \not= \beta \in I$. If
\[
x^m+\sum_{i=0}^{m-1} z_i x^i \]
is the minimal polynomial of $\beta$ over $\mathbb{Q}$ then $z_0 \not=0$ (otherwise, $\beta$ is the root of a nonzero polynomial over $\mathbb{Q}$ of degree less that $m$) and
\[
z_0=-\beta^m-\sum_1^{m-1} z_i\beta^i \in I,\]
hence $\pm z_0 \in I$, and so $I$ contains a positive integer $a$. 

\emph{ Claim}: each element of $R$ can be expressed in the form
\[
a\gamma+\sum_1^n r_i \alpha_i,\] 
where $\gamma \in R,\  r_i \in \{0,\dots, a-1\}, i=1,\dots,n.$

Assume this for now, and let $\alpha \in J \setminus I$. Then for each $k \in \{1, 2, 3,\dots \}$,
\[
\alpha^k=a\gamma_k+\sum_1^n r_{ik} \alpha_i,\ \gamma_k \in R,\  r_{ik} \in \{0,\dots, a-1\},\ i=1,\dots,n,\]
hence the sequence $(\alpha^k-a\gamma_k: k=1, 2, 3,\dots)$ has only finitely many values; consequently there exist positive integers $l<k$ such that
\[
\alpha^l-a \gamma_l=\alpha^k-a\gamma_k.\]
Hence
\[
\alpha^l(\alpha^{k-l}-1)=\alpha^k-\alpha^l=a(\gamma_k-\gamma_l) \in I\  (a \in I\ !).\]
Because $I$ is prime, either $\alpha^l \in I$ or $\alpha^{k-l}-1 \in I$. However, $\alpha^l \not \in I$ because $\alpha \not \in I$ and $I$ is prime. Hence
\[
 \alpha^{k-l}-1 \in I \subseteq J.\]
 But $k-l>0$ and $\alpha \in J$ (by the choice of $\alpha$), and so $-1 \in J$. As $J$ is an ideal, this implies that $J=R$ .

Our claim must now be verified. Let $\alpha \in R$, and find $z_i \in \mathbb{Z}$ such that 
\[
\alpha=\sum_{i=1}^n  z_i\alpha_i.\]
The division algorithm in $\mathbb{Z}$ implies that there exist $m_i \in \mathbb{Z}$, $r_i \in\{0,\dots, a-1\},\ i=1,\dots,n$, such that $z_i=m_ia+r_i,\ i=1,\dots,n$. Thus
\[
\alpha=a\sum_i m_i\alpha_i+\sum_ir_i\alpha_i=a\gamma+\sum_ir_i\alpha_i,\]
with $\gamma \in R$. $\hspace{13.2cm}\textrm{QED}$

\vspace{0.02cm}
Proposition 49 is a key tool needed in the proof of the next theorem, and that theorem is one of the most fundamental results of algebraic number theory. It was  first proved by Richard Dedekind in 1895.

\vspace{0.4cm}
\textbf{Theorem 50}. \textit{$($Fundamental Theorem of Ideal Theory$)$. If I is a nonzeo, proper ideal of R then I is a product of prime ideals, and this factorization of I is unique up to the order of the factors}.

\vspace{0.4cm}
By virtue of Theorem 50, even though unique prime factorization of the elements of $R$ may fail to hold, unique prime-ideal factorization of the \emph{ideals} of $R$ always holds.

The proof of Theorem 50 requires three lemmas. We will proceed by first proving the lemmas and then proving Theorem 50. 

Let $I$ be an ideal of $R$, $\{0\}\not= I \not= R$.

\vspace{0.4cm}
\textbf{Lemma 51}. \textit{There exist prime ideals $P_1,\dots,P_s$ of $R$ such that $I\subseteq P_i$, for all $i$ and $P_1\cdots P_s\subseteq I$.}

\vspace{0.4cm}
\emph{Proof}. If $I$ is prime, done, with $s=1$, hence suppose that $I$ is not prime. Then there exists a product $\beta \gamma$ of elements of $R$ which is in $I$ and $\beta\not \in I,\ \gamma \not \in I$. Let $\{\alpha_1,\dots,\alpha_n\} $ be a basis of $I$, and set
\[
J=(\alpha_1,\dots,\alpha_n,\beta),\ K=(\alpha_1,\dots,\alpha_n,\gamma).\]
Then
\[
JK\subseteq I,\ I\subsetneqq J,\ I\subsetneqq K.\]
If $J, K$ are both prime, done, with $s=2$. Otherwise apply this procedure to each nonprime ideal that occurs, and continue in this way as long as the procedure produces nonprime ideals. Note that after each step of the procedure,

$(i)$ the product of all the ideals obtained in that step is contained in $I$,

$(ii)$ $I$ is contained in each ideal obtained in that step, and

$(iii)$ each ideal obtained in that step is properly contained in an ideal from the immediately preceding step.

\emph{Claim}: this procedure terminates after finitely many steps.

If this is true then each ideal obtained in the final step is prime; otherwise the procedure would continue by applying it to a nonprime ideal. If $P_1,\dots,P_s$ are the prime ideals obtained in the final step then this sequence of ideals satisfies Lemma 51 by virtue of $(i)$ and $(ii)$ above. 

\emph{Proof of the claim}. Suppose this is false. The $(ii)$ and $(iii)$ above imply that the procedure produces a sequence of ideals $J_0, J_1,\dots, J_n,\dots$ such that $J_0=I$ and $J_i\subsetneqq J_{i+1}$, for all $i$. We will now prove that $I$ is contained in only finitely many ideals, hence no such sequence of ideals is possible.

The proof of Proposition 49 implies that $I$ contains a positive rational integer $a$. We show: $a$ belongs to only finitely many ideals.

Suppose that $J$ is an ideal, with basis $\{\beta_1,\dots,\beta_n\}$, and $a \in J$. Then we also have that
\[
J=(\beta_1,\dots,\beta_n, a).\]
By the claim in the proof of Proposition 49, for each $i$, there is $\gamma_i, \delta_i\in R$ such that $\beta_i=a\gamma_i+\delta_i$, and $\delta_i$ can take on only at most $an$ values. But then
\[
J=(a\gamma_1+\delta_1,\dots,a\gamma_n+\delta_s)=(\delta_1,\dots,\delta_n, a).\]
Because each $\delta_i$ assumes at most $an$ values, it follows that $J$ is one of only at most $an^2$ ideals. $\hspace{14.2cm}\textrm{QED}$

The statement of the next lemma requires the following definition:

\vspace{0.4cm}
\textbf{Definition}. If $J$ is an ideal of $R$ then 
\[
J^{-1}=\{\alpha \in F: \alpha \beta \in R,\ \textrm{ for all}\ \beta \in J\}.\]

\vspace{0.4cm}
\textbf{Lemma 52}. \textit{If P is a prime ideal of R then $P^{-1}$ contains an element of $F\setminus R$.}

\vspace{0.4cm}
\emph{Proof}. Let $x\in P$. Lemma 51 implies that $(x)$ contains a product $P_1\cdots P_s$ of prime ideals. Choose a product with the smallest number $s$ of factors.

Suppose that $s=1$. Then $P_1\subseteq (x)\subseteq P$. $P_1$ maximal (Proposition 49) implies that $P=P_1=(x)$. Hence $1/x\in P^{-1}$. Also, $1/x\not \in R$; otherwise, $1=x\cdot 1/x \in P$, contrary to the fact that $P$ is proper.

Suppose that $s>1$. Then $P_1\cdots P_s \subseteq (x) \subseteq P$, so Proposition 48 implies that $P$ contains a $P_i$, say $P_1$. $P$ maximal implies that $P=P_1$. $P_2\cdots P_s \nsubseteq (x)$ by minimality of $s$, hence there exits $\alpha \in P_2\cdots P_s$ such that $\alpha \not \in (x)$, and so $\alpha/x\not \in R$.

\emph{Claim}: $\alpha/x \in P^{-1}$. 

Let $\beta \in P$. We must prove: $\beta (\alpha/x)\in R$. To do that, observe that
\[
(\alpha)P \subseteq P_2\cdots P_sP=P_1\cdots P_s \subseteq (x),\]
and so there is a $\gamma \in R$ such that $\alpha \beta=x\gamma$, i.e., $\beta(\alpha/x)=\gamma$.$\hspace{4.6cm}\textrm{QED}$

The next lemma is the key technical tool that allows us to prove Theorem 50; it will be used to factor an ideal into a product of prime ideals and to show that this factorization is unique up to the order of the factors. In order to state it, we need to extend the definition of products of ideals to products of  arbitrary subsets of $R$ like so:

\vspace{0.4cm}
\textbf{Definition}. If $S$ and $T$ are subsets of $R$ then the \emph{product ST of S and T} is the set
\[
ST=\left\{\sum_i s_it_i: (s_i, t_i)\in S \times T\right\}.\]

\vspace{0.4cm}
This product is clearly commutative and associative, and it agrees with the product defined before when $S$ and $T$ are ideals of $R$.

\vspace{0.4cm}
\textbf{Lemma 53}. \textit{If P is a prime ideal of R and I is an ideal of R then $P^{-1}PI=I$}.

\vspace{0.4cm}
\emph{Proof}. It suffices to show that $P^{-1}P=(1)$. $J=P^{-1}P$ is an ideal of $R$ (prove!). As $1\in P^{-1}$, it follows that $P\subseteq J$ and so $P$ maximal implies that $P=J$ or $J=(1)$.

Suppose that $J=P$. Let $\{\alpha_1,\dots,\alpha_n\}$ be a basis of $P$, and use Lemma 52 to find $\gamma \in P^{-1},\ \gamma \not \in R$. Then $\gamma \alpha_i \in P$, for all $i$, and so
\[
\gamma \alpha_i=\sum_j a_{ij}\alpha_j,\ \textrm{where $a_{ij} \in \mathbb{Z}$ for all $i$, $j$}.\]
As a consequence of these equations, $\gamma$ is an eigenvalue of the matrix $[a_{ij}]$, hence it is a root of the characteristic polynomial of $[a_{ij}]$. Since this characteristic polynomial is a monic polynomial in $\mathbb{Z}[x]$, Lemma 22 implies that $\gamma$ is an algebraic integer, contrary to its choice. Hence $P\not=J$, and so $J=(1)$.$\hspace{11.5cm}\textrm{QED}$

\vspace{0.4cm}
\textbf{Proposition 54}. \textit{The only ideal factors that a prime ideal $P$ has are $P$ and $(1)$}.

\vspace{0.4cm}
\emph{Proof}. Suppose that $P=IJ$, $I$ and $J$ ideals of $R$. Proposition 48 implies that we may assume with no loss of generality that $I \subseteq P$. If $I=(1)$ then $P=J$, done, hence assume that $I$ is proper. Because $P=IJ \subseteq I$, it follows that $P=I$. Lemma 53 implies that
\[
(1)=P^{-1}P=P^{-1}IJ=P^{-1}PJ=J.\]

$\hspace{15cm}\textrm{QED}$

Theorem  50 is a consequence of the next two propositions.

\vspace{0.4cm}
\textbf{Proposition 55}. \textit{Every nonzero proper ideal of R is a product of prime ideals}.

\vspace{0.4cm}
\emph{Proof}. Lemma 51 implies that every nonzero proper ideal of $R$ contains a product $P_1\cdots P_r$ of prime ideals, where we choose a product with the smallest number $r$ of factors. The argument now proceeds by induction on $r$.

Let $\{0\}\not= I \not= R$ be an ideal with $r=1$, i.e., $I$ contains a prime ideal $P$. $P$ maximal implies that $I=P$, done.

Assume now that $r>1$ and every nonzero, proper ideal that contains a product of fewer than $r$ prime ideals is a product of prime ideals. 

Let $\{0\}\not= I \not= R$ be an ideal that contains a product $P_1\cdots P_r$ of prime ideals, with $r$ the smallest number of prime ideals with this property. Lemma 51 implies that
$I \subseteq Q$, $Q$ a prime ideal. Hence  $P_1\cdots P_r\subseteq Q$, and so Lemma 48 implies that $Q$ contains a $P_i$, say $P_1$.  $P_1$ maximal implies that $Q=P_1$. Hence $I \subseteq P_1$. Then $IP_1^{-1}$ is an ideal of $R$; $I\subseteq IP_1^{-1}$ ($1\in P^{-1}$), and so  $IP_1^{-1} \not= \{0\}$. $IP_1^{-1} \not=R$; otherwise, $P_1\subseteq I$, hence $I=P_1$, contrary to the fact that $r>1$. Lemma 53 implies that
\[
P_2\cdots P_r \subseteq P_1^{-1}P_1\cdots P_r \subseteq IP_1^{-1},\]
hence by the induction hypothesis, $IP_1^{-1}$ is a product $P_1^{\prime}\cdots P_k^{\prime}$ of prime ideals, and so by Lemma 53 again,
\[
I=(IP_1^{-1})P_1=P_1^{\prime}\cdots P_k^{\prime}P_1\]
is a product of prime ideals.$\hspace{10.7cm}\textrm{QED}$

\vspace{0.4cm}
\textbf{Proposition 56}. \textit{Factorization as a product of prime ideals is unique up to the order of the factors}.

\emph{Proof}. Suppose that $P_1\cdots P_r=Q_1\cdots Q_s$ are products of prime ideals, with $r\leq s$, say. $Q_1\cdots Q_s\subseteq Q_1$, hence $P_1\cdots P_r\subseteq Q_1$  and so Lemma 48 and the maximality of the $P_i$'s imply, after reindexing one of the $P_i$'s, that $Q_1=P_1$. Then Lemma 53 implies that
\[
P_2\cdots P_r=P_1^{-1}P_1\cdots P_r=Q_1^{-1}Q_1\cdots Q_s=Q_2\cdots Q_s.\]
Continuing in this way, we deduce, upon reindexing of the  $P_i$'s, that $P_i=Q_i$, $i=1,\dots r$, and also, if $r<s$, that
\[
(1)=Q_{r+1}\cdots Q_s.\]
But this equation implies that $R=(1)\subseteq Q_{r+1}$, which is impossible as $Q_{r+1}$ is a proper ideal. Hence $r=s$.$\hspace{12.1cm}\textrm{QED}$

\vspace{0.4cm}
\textbf{Definition}. If $I$ is a nonzero, proper ideal of $R$ and $P$ is a prime-ideal factor of $I$, then the highest power of $P$ that occurs in the prime-ideal factorization of $I$ is called the \emph{multiplicity of P in I}

\vspace{0.4cm}
N.B. If $I$ is a nonzero, proper ideal of $R$ and $I$ is contained in a prime ideal $Q$ of $R$ then $Q$ must be one of the factors of $I$ which occur in the prime-ideal factorization of $I$, and each of these factors clearly contains $I$, i.e., the set $\{P_1,\dots,P_k\}$ of prime-deal factors of $I$ is precisely the set of prime ideals containing $I$. If $m_i$ is the multiplicity of $P_i$ in $I$ then we can factor $I$ as
\[
I=P_1^{m_1}\cdots P_k^{m_k},\]
which is an exact analog for ideals of the prime factorization of the positive rational integers.
\chapter{Some Structure Theory for Ideals in a Number Ring}

In this chapter, we will illustrate how the Fundamental Theorem of Ideal Theory is used to obtain some important features of the structure of ideals in a number ring. 

We start with an elegant proof of Theorem 45: for a number ring $R$, UFD implies PID.

\vspace{0.02cm}
\emph{Proof of Theorem} 45.

Assume that $R$ is a UFD. Suppose that every prime ideal of $R$ is principal. If $\{0\}\not= I \not= R$ is an ideal of $R$ with prime factorization $\prod_i P_i$ then choose $\alpha_i \in R$ such that $P_i=(\alpha_i)$, for all $i$, to obtain
\[
I=\prod_i P_i=\prod_i (\alpha_i)=\Big(\prod_i \alpha_i\Big).\]
Hence we need only prove that each prime ideal $P$ of $R$ is principal.

Let $0\not=z \in P$. Then $z$ is not a unit since $P$ is proper. Let $z=\prod_i \rho_i$ be the factorization of $z$ into prime elements of $R$. Then $\prod_i(\rho_i)=(z) \subseteq P$, hence by proposition 48, $(\rho)\subseteq P$ for some prime element $\rho$ of $R$. $R$ a UFD implies that $\rho$ has the following property: if $\alpha, \beta \in R$ and $\rho$ divides $\alpha\beta$ in $R$ then $\rho$ divides either $\alpha$ or $\beta$ in $R$ (prove!). But this says that $(\rho)$ is a prime ideal, hence maximal, in $R$. Thus $P=(\rho)$ is principal.$\hspace{13.9cm}\textrm{QED}$

\vspace{0.2cm}
Proposition 40 implies that the ring $\mathbb{Z}+\sqrt{-5}\ \mathbb{Z}$ is not a UFD, hence Theorem 45 implies that $\mathbb{Z}+\sqrt{-5}\ \mathbb{Z}$ has a non-principal ideal. In fact, the ideal $(3, 1+2\sqrt{-5})$ is not principal (prove!).

Let $R$ be a fixed number ring.

\vspace{0.4cm}
\textbf{Proposition 57}. \textit{If I and J are ideals of R with $I\not= \{0\}$ then $J\subseteq I$ if and only if there exists an ideal $K$ of $R$ such that $J=IK$. Moreover, $K$ is uniquely determined by this equation.}

\vspace{0.4cm}
\emph{Proof}. It is clear that $J=IK$ for some ideal $K$ of $R$ implies that $J\subseteq I$. For the converse, we assume that $J\subseteq I$ and, without loss of generality, that $I \not= (1)$ and $\{0\}\not= J\not= (1)$. Let $I=P_1\cdots P_r$ and $J=Q_1\dots Q_s$ be the prime (ideal) factorizations of $I$ and $J$. N.B. The factors occurring in these factorizations may not be all distinct. The proof of Proposition 56 implies that every prime factor $P_i$ of $I$ occurs as a prime factor $Q_j$ of $J$, hence by reindexing the $P_i$'s, we have $r\leq s$ and $P_i=Q_i$, $i=1,\dots,r$.

If $r=s$ then $I=J$, so take $K=(1)$. If $r<s$, take $K=Q_{r+1}\cdots Q_s$; then
\[
J=P_1\cdots P_rQ_{r+1}\cdots Q_s=KI.\]

To prove uniqueness, suppose that $K$ and $L$ are ideals with $KI=LI$. Then using Lemma 53, we cancel off all of the prime factors of $I$ from this equation to obtain $K=L$.$\hspace{1cm}\textrm{QED}$

\vspace{0.4cm}
\textbf{Definitions}. If $I$ and $J$ are ideals of $R$ with $I\not= \{0\}$ and $J\subseteq I$ then we will say that $I$ \emph{divides J}. If $I$ divides $J$ then the ideal $K$ in the conclusions of Proposition 57 is called the \emph{quotient of J by I} and is denoted by $J/I$ or $\displaystyle{\frac{J}{I}}$.

\vspace{0.4cm}
The following important corollary is an immediate consequence of Propositions 48 and 57; it states that, with respect to ideal factorization, the prime ideals possess the `` fundamental property of primes".

\vspace{0.4cm}
\textbf{Corollary 58}. \textit{If a prime ideal divides an ideal product then it must divide at least one of the factors of the product. }

\vspace{0.4cm} 
\textbf{Definition}. Let $I$ and $J$ be ideals of $R$. A \emph{greatest common divisor of I and J}, denoted by $\gcd(I, J)$, is an ideal $K$ of $R$ with the following properties: $K$ divides both $I$ and $J$, and if $L$ is an ideal of $R$ dividing both $I$ and $J$ then $L$ divides $K$.

\vspace{0.4cm}
\textbf{Proposition 59}. \textit{$(i)$ If $I\not=\{0\} \not=J$ then
the greatest common divisor of I and J is unique}.

\textit{$(ii)$ Suppose that $\{0\}\not= I\not= (1)$ and $\{0\}\not= J\not= (1)$. If $I$ and $J$ have no common prime-ideal factors, then $\gcd(I, J)=(1)$. If $P_1,\dots,P_r$ are the distinct prime-ideal factors that are the common factors of \emph{both} I and J and if $e_i$ is the highest power of $P_i$ which divides both I and J then}
\[
\gcd(I, J)=\prod_{i=1}^r P_i^{e_i}.\]

\textit{$(iii)$ If $I=(\alpha_1,\dots,\alpha_r)$ and $J=(\beta_1,\dots,\beta_s)$ then}
\[
\gcd(I, J)=(\{\alpha_1,\dots,\alpha_r\} \cup \{\beta_1,\dots,\beta_s\}).\]

\vspace{0.4cm}
\emph{Proof}. $(i)$ and $(iii)$ follow straightforwardly from Proposition 57 and $(ii)$ follows straightforwardly from the Fundamental Theorem of Ideal Theory. $\hspace{5.1cm}\textrm{QED}$ 

N.B. If $I$ and $J$ are nonzero ideals of $R$ then Proposition 46 and Proposition 59$(i), (iii)$ imply that the greatest common divisor of $I$ and $J$ exists and is unique.

\vspace{0.4cm}
\textbf{Corollary 60}. \textit{If I and J are ideals of R and $P_1,\dots,P_r$ are the distinct prime factors of I then the $\gcd(I, J)=(1)$ if and only if $\gcd(P_i, J)=(1)$, for all $i$.}

\vspace{0.4cm}
A corollary of Proposition 46 is that every ideal in a number ring is finitely generated, i.e., every ideal is generated by a finite set. We will now use Propositions 57 and 59 and Corollaries 58 and 60 to show that in fact every ideal in a number ring is generated by at most only two elements; moreover, one of the generators can be taken to be any nonzero element of the ideal.

\vspace{0.4cm}
\textbf{Lemma 61}. \textit{If $I\not=\{0\} \not=J$ are ideals then there exists $0\not=\alpha\in I$ such that $\gcd\big((\alpha)/I, J\big)= (1)$.}

\vspace{0.4cm}
N.B. $\alpha \in I$ implies that $(\alpha)\subseteq I$, and so the ideal quotient $(\alpha)/I$ exists and is unique by Proposition 57.

\vspace{0.4cm}
\emph{Proof}. If $J=(1)$ then any  $0\not=\alpha\in I$ will work, since
\[
\gcd\big((\alpha)/I, J\big)=\gcd\big((\alpha)/I, (1)\big)=(1).\]
Similarly, if $I=(\gamma)$ is principal then
\[
\gcd\big((\gamma)/I, J\big)=\gcd\big((1), J\big)=(1).\]
Hence assume that $J\not=(1)$, $I$ is not principal, and let $P_1,\dots,P_r$ be the distinct prime factors of $J$.

Suppose that $r=1$; then set $P=P_1$. Corollary 60 implies that we must find $0\not=\alpha\in I$ such that $\gcd\big((\alpha)/I, P\big)= (1)$. Choose $\alpha \in I\setminus IP$. Such an $\alpha$ exits; if not then $I=IP$, hence the uniqueness in Proposition 57 implies that $P=(1)$, not possible since $P$ is proper.

This $\alpha$ works. Let $L=(\alpha)/I$. Then $\{0\}\not=L$, and because $I$ is not principal, $L\not=(1)$. Suppose that $\gcd(L, P)\not=(1)$. Proposition 59$(ii)$ implies that $P|L$, i.e., $L=MP$ for some ideal $M$, and so $(\alpha)=IL=MIP$, and so $\alpha \in IP$, contrary to its choice. Hence $\gcd\big((\alpha)/I, P\big)= (1)$.

Next, suppose that $r>1$. Corollary 60 implies that we must find $0\not= \alpha \in I$ such that
\begin{equation*}
\gcd \big((\alpha)/I, P_k\big)=(1),\ \textrm{for all $k$}.\tag{$*$}
\]

For each $k=1,\dots,r$, consider the ideals
\[
I_k=\frac{IP_1\cdots P_r}{P_k}\ \textrm{and $P_k$}.\]
From the $r=1$ case, we find $\alpha_k \in I_k\setminus I_kP_k$. Let $\alpha=\sum_k \alpha_k$. Since $I|I_k$, we have that $I_k\subseteq I$, for all $k$, and so $\alpha_k\in I$, for all $k$, hence $\alpha \in I$.

\emph{Claim}: $\alpha\not \in IP_k$, for all $k$.

If this is so, then $\alpha\not=0$, and the argument from the $r=1$ case implies that $(*)$ is true.

\emph{Proof of the claim}. For $j\not=k$,
\[
(\alpha_j)\subseteq I_j=\frac{IP_1\cdots P_r}{P_j}=IP_k\frac{P_1\cdots P_r}{P_kP_j}\subseteq IP_k,\]
and so $\alpha_j\in IP_k$, for all $j\not=k$.

Suppose now by way of contradiction that $\alpha \in IP_k$. Then
\[
\alpha_k=\alpha-\sum_{j\not=k} \alpha_j \in IP_k,\]
contrary to the choice of $\alpha_k$.$\hspace{10.4cm}\textrm{QED}$ 

\vspace{0.4cm}
\textbf{Theorem 62}. \textit{$($Ideal Generation Theorem$)$ If I is a nonzero ideal of R and $0\not= \beta \in I$, then there exist $\alpha\in I$ such that $I=(\alpha, \beta)$.}

\vspace{0.4cm}
\emph{Proof}. As Theorem 62 is trivial for $I=(1)$, we can assume that $I\not=(1)$. Lemma 61 implies that there exists $0\not= \alpha \in I$ such that
\[
\gcd\big((\alpha)/I, (\beta)/I)\big)= (1).\]

\emph{Claim}: $\gcd\big((\alpha), (\beta)\big)=I$.

If this is true then Proposition 59$(iii)$ implies that $I=(\alpha, \beta)$.

\emph{Proof of the claim}. We will show that

$(i)$ the set of prime factors of $I$ =  the set of common prime factors of $(\alpha)$ and $(\beta)$;

$(ii)$ If $P$ is a prime factor of $I$ and $e, f,$ and $g$ is, respectively, the multiplicity of $P$ in $I, (\alpha)$, and $(\beta)$, then $e=\min \{f, g\}$.

\noindent The claim will then follow from $(i), (ii)$, and Proposition 59 $(ii)$.

In order to verify $(i)$, note first that $(\alpha)\subseteq I, (\beta)\subseteq I$, hence $I$ is a factor of both $(\alpha)$ and $(\beta)$, hence every prime factor of $I$ is a common prime factor of $(\alpha)$ and $(\beta)$. Let $Q$ be a common prime factor of $(\alpha)$ and $(\beta)$. Then $Q$ divides the ideal products
\[
\big((\alpha)/I\big)I\ \textrm{and}\ \big((\beta)/I\big)I.\]

Suppose that $Q$ is not a factor of $I$. Corollary 58 implies that $Q$ divides both $(\alpha)/I$ and $(\beta)/I$, which is not possible because $\gcd\big((\alpha)/I, (\beta)/I)\big)= (1).$ Hence $Q$ is a factor of $I$.

As for $(ii)$, let $P, e, f, g$ be as in $(ii)$, and set $h=\min \{f, g\}$. Because $P^e$ divides $I$, it also divides  $(\alpha)$ and $(\beta)$, hence $e\leq h$.

Suppose that $e<h$. Setting $P^{-e}=(P^{-1})^e$, we have that
\[
P^{-e}(\alpha)=\big((\alpha)/I\big)P^{-e}I.\]
Since $P^h$ divides $(\alpha)$, it follows that $P^{h-e}$, and hence $P$, divides the left-hand side, and thus the right-hand side, of this equation. Because $e$ is the highest power of $P$ which divides $I$, it follows that $P^{-e}I$ is an ideal which is \emph{not} divisible by $P$. Hence Corollary 58 implies that $P$ divides $(\alpha)/I$. But by the same reasoning, $P$ also divides $(\beta)/I$, and that is impossible since $\gcd\big((\alpha)/I, (\beta)/I)\big)= (1).$ Hence $e=h$.$\hspace{7.5cm}\textrm{QED}$ 
\chapter{An Abstract Characterization of Ideal Theory in a Number Ring} 

If $D$ is an integral domain then $D$ is contained in a field $F$ such that for all $x\in F$, there exits $a, b\in D, b\not=0$, such that $x=ab^{-1}$. The field $F$ is uniquely determined by this property and is called the \emph{ field of fractions of D} (Hungerford [8], Theorem III.4.3, Corollary III.4.6).

\vspace{0.4cm}
\textbf{Definition}.  An integral domain $D$ is \emph{integrally closed} if the following condition is satisfied: if $F$ is the field of fractions of $D$ and if $\alpha \in F$ is the root of a monic polynomial in $D[x]$, then $\alpha \in D$, i.e., the only roots in $F$ of monic polynomials with coefficients in $D$ are the elements of $D$.

\vspace{0.4cm}
\textbf{Definition}. An integral domain $D$ is a \emph{Dedekind domain} if $D$ satisfies the following conditions:

$(i)$ every ideal of $D$ is finitely generated;

$(ii)$ every prime ideal of $D$ is maximal;

$(iii)$ $D$ is integrally closed.

\vspace{0.4cm}
\textbf{Proposition 63}. \textit{Every number ring is a Dedekind domain.}

\vspace{0.4cm}\emph{Proof}. Let $F$ be a number field, $R=\mathcal{A} \cap F$. Proposition 46 (respectively, Proposition 49) implies that $(i)$ (respectively, $(ii)$) in the definition of Dedekind domain holds for $R$.

We need to verify $(iii)$ in the definition of Dedekind domain for $R$. Let $\alpha \in F$. Lemma 28 implies that there exits $0\not=n \in \mathbb{Z}$ such that $n\alpha \in R$, hence $\alpha=n\alpha/n=a/b$ with $a, b \in R$, and so $F$ is the field of fractions of $R$.

Let $\alpha \in F$ be a root of a monic polynomial $\mu(x)=x^m+\sum_{i=0}^{m-1} a_i x^i$ in $R[x]$. We must prove: $\alpha \in R$.

\emph{Claim}: $\mathbb{Z}[a_0,\dots,a_{m-1}, \alpha]$ is finitely generated over $\mathbb{Z}$, i.e., there exist nonzero elements $p_1,\dots,p_k$ of  $\mathbb{Z}[a_0,\dots,a_{m-1}, \alpha]$  such that
\[
\mathbb{Z}[a_0,\dots,a_{m-1}, \alpha]=\sum_{i=1}^k p_i\mathbb{Z}.\]
If this claim is true then
\[
\alpha p_i=\sum_{j=1}^k z_{ij} p_j,\ \textrm{with $z_{ij} \in \mathbb{Z}$, for all $i, j$}.\]
Hence $\alpha$ is a root of the characteristic polynomial of the matrix $[z_{ij}]$ and this polynomial is a monic polynomial over $\mathbb{Z}$, whence Lemma 22 implies that $\alpha \in \mathcal{A} \cap F=R$.

\emph{Proof of the claim}. We assert first that $\mathbb{Z}[a_0,\dots,a_{m-1}]$ is finitely generated over $\mathbb{Z}$. To see this, first consider the ring $\mathbb{Z}[a_0]$. We have $a_0\in \mathcal{A}$ and so $a_0$ is a root of a monic polynomial $p$ over $\mathbb{Z}$, of degree $l$, say. If $s\in \mathbb{Z}[x]$ then apply the division algorithm in $\mathbb{Z}[x]$ to find $q, r \in \mathbb{Z}[x]$ such that
\[
s=pq+r,\ 0\leq\ \textrm{degree of $r<l$,}\]
hence $s(a_0)=r(a_0)$, and so $\mathbb{Z}[a_0]$ is generated over $\mathbb{Z}$ by $\{1,\dots,a_0^{l-1}\}$.

Now, assume inductively that $m>2$ and $A=\mathbb{Z}[a_0,\dots,a_{m-2}]$ is generated over $\mathbb{Z}$ by $p_1,\dots,p_k$. Let $p$ now denote the minimal polynomial over $\mathbb{Q}$ of $a_{m-1}$, of degree $t$, say. Set
\[
S=\big\{a_{m-1}^ip_j: (i, j)\in \{0,\dots,t-1\} \times \{1,\dots,k\}\big\},\]
and let $w\in \mathbb{Z}[a_0,\dots,a_{m-1}]$. Then there is an $s \in A[x]$ such that $w=s(a_{m-1})$. Use the division algorithm in $A[x]$ to find $q, r \in A[x]$ such that $s=pq+r$, $0\leq$ degree of $r<t$. Then $w=r(a_{m-1})$, and since every coefficient of $r(x)$ is in $A$, it follows that $w$ is in the set generated over $\mathbb{Z}$ by $S$. Hence $\mathbb{Z}[a_0,\dots,a_{m-1}]$ is finitely generated over $\mathbb{Z}$.

If $\beta \in \mathbb{Z}[a_0,\dots,a_{m-1}, \alpha]$ then $\beta=s(\alpha)$ for some polynomial $s(x)$ in $\mathbb{Z}[a_0,\dots,a_{m-1}][x]$. The polynomial $\mu(x)$ is monic of degree $m$ over $\mathbb{Z}[a_0,\dots,a_{m-1}]$ with root $\alpha$, and so by dividing $s(x)$ by $\mu(x)$ using the division algorithm in  $\mathbb{Z}[a_0,\dots,a_{m-1}][x]$, we find as before a polynomial $r$ over $\mathbb{Z}[a_0,\dots,a_{m-1}]$ of degree $<m$ such that $\beta=r(\alpha)$. If we now take a set $p_1^{\prime},\dots, p_w^{\prime}$ of generators of $\mathbb{Z}[a_0,\dots,a_{m-1}]$ over $\mathbb{Z}$, it follows that $\beta$ is in the set generated over $\mathbb{Z}$ by $\big \{\alpha^ip_j^{\prime}: (i, j)\in \{0,\dots,m-1\} \times \{1,\dots,w\}\big\}$.$\hspace{7.7cm}\textrm{QED}$ 

\vspace{0.2cm}
N.B. The converse of Proposition 63 is far from true. There are Dedekind domains of positive characteristic; no such Dedekind domain can be isomorphic to a number ring, as all number rings have characteristic 0.

\vspace{0.4cm}
\textbf{Theorem 64} \textit{$($Fundamental Theorem of Ideal Theory for Dedekind Domains$)$. Every nonzero proper ideal in a Dedekind domain is a unique product of prime ideals.}

\vspace{0.4cm}
\emph{Proof}. See Marcus [9],  Chapter 3, Theorem 16.$\hspace{6.7cm}\textrm{QED}$ 

\vspace{0.2cm}
Theorem 45 and all of the results from Proposition 57 through Theorem 62 above remain valid for arbitrary Dedekind domains. Moreover, all of these results can be proved for general Dedekind domains in exactly the same way as they were proved for number rings, since the only things required for those proofs are the properties coming from the definition of Dedekind domain  and Theorem 64. A Dedekind domain residing in its field of fractions is hence the abstract analog of the number ring $\mathcal{A} \cap F$ residing in a number field $F$. Many modern accounts of algebraic number theory, including Marcus' treatment in [9], take this more general approach.
\chapter{Ideal-Class Group and the Class Number}

We now develop a way to give a precise, quantitative measure of how far a number ring is from being a UFD.

Let $F$ be a number field with number ring $R=\mathcal{A} \cap F$, fixed for the rest of this chapter.

\vspace{0.4cm}
\textbf{Definition}. If $I$ and $J$ are ideals of $R$ then $I$ is \emph{equivalent to J}, denoted by $I\sim J$, if there exist nonzero elements $\alpha, \beta \in R$ such that
\[
(\alpha)I=(\beta)J.\]

\vspace{0.4cm}
\textbf{Proposition 65}. \textit{$\sim$ is an equivalence relation on the set of all ideals of R.}

\vspace{0.4cm}
\emph{Proof}. Exercise.$\hspace{12cm}\textrm{QED}$ 

\vspace{0.4cm}
\textbf{Definitions}. If $I$ is an ideal of $R$ then $[I]$ denotes the equivalence class of $\sim$ which contains $I$. $[I]$ is called the \emph{ideal class containing $($or determined by$)$ I}. The ideal class $[(1)]$ containing $(1)$ is called the \emph{principal class}.

\vspace{0.4cm}
\textbf{Lemma 66}. \textit{The principal class is the set of all nonzero principal ideals of R.}

\vspace{0.4cm}
\emph{Proof}. If $I\in [(1)]$ then there exist nonzero $\alpha, \beta \in R$ such that $(\alpha)I=(\beta)$, hence $I=(\beta/\alpha)R$, and so $\beta/\alpha \in I$, whence $\{0\}\not= I= (\beta/\alpha)$ is principal. If $I=(\alpha)$ is principal, with $\alpha\not=0$, then $(1)I=(\alpha)(1)$, hence $I\sim (1)$.$\hspace{8.2cm}\textrm{QED}$ 

\vspace{0.4cm}
\textbf{Definition}. If $[I]$ and $[J]$ are ideal classes of $R$ then the \emph{product $[I][J]$ of $[I]$ and $[J]$} is the ideal class $[IJ]$ containing $IJ$.

\vspace{0.4cm}
\textbf{Proposition 67}. \textit{$(i)$ The ideal-class product is well-defined, i.e., if I, J, K, L are ideals of R with $I\sim J$ and $K\sim L$ then $IK\sim JL$.}

\textit{$(ii)$ The ideal-class product is commutative and associative.}

\textit{$(iii)$ For all ideals I of R, $[I][(1)]=[I]$.}

\vspace{0.4cm}
\emph{Proof}. Exercise.$\hspace{12cm}\textrm{QED}$ 

\vspace{0.4cm}
\textbf{Proposition 68}. \textit{If $I\not= \{0\}$ is an ideal of R then there exists an ideal $J\not=\{0\}$ of R such that IJ is a principal ideal.}

\vspace{0.4cm}
\emph{Proof}. Let $0\not=\alpha \in I$ and set
\[
J=\{\beta \in R: (\beta)I\subseteq (\alpha)\}.\]
$J$ is an ideal of $R$, nonzero since $\alpha \in J$.

Now consider the ideal $K=(1/\alpha)IJ$ of $R$. If $K=(1)$ then $IJ=(\alpha)$, done, so assume by way of contradiction that $K\not= (1)$. Then we can find a prime ideal $P$ containing $K$. Lemma 52 implies that
\begin{equation*}
\textrm{there exits a}\ \gamma \in P^{-1}\cap (F\setminus R).\tag{1}\]
Then by definition of $P^{-1}$ and the fact that $K\subseteq P$, it follows that
\begin{equation*}
(\gamma)K\subseteq R.\tag{2}\]

\emph{Claim}: $(\gamma)J\subseteq J.$

To see this, let $\beta \in J, \delta \in I$; we must prove that $\gamma\beta\delta \in (\alpha)$, i.e., $(\gamma\beta\delta)/\alpha \in R$. But this is an immediate consequence of the definition of $K$ and (2).

Let $\{\alpha_1,\dots,\alpha_n\}$ be a basis of $J$. The claim implies that $\gamma\alpha_i \in J$, for all $i$. Now follow the argument in the proof of Lemma 51 to conclude that $\gamma \in R$, which contradicts (1). QED

After observing that $[(0)]=\big\{\{0\}\big\}$, we let 
\[
\mathcal{C}(R)=(\textrm{the set of all ideal classes of $R$})\setminus \big\{\{0\}\big\}.\]

Let $[I] \in \mathcal{C}(R)$. Then $I\not= \{0\}$, and so Proposition 68 implies that there exits an ideal $J\not= \{0\}$ and $0\not= \alpha \in R$ such that $IJ=(\alpha)$. Hence Lemma 66 implies that $[I][J]=[(1)]$. This, together with Proposition 67, shows that $\mathcal{C}(R)$ becomes an abelian group when endowed with the ideal-class product, whose identity element is the principal class.

\vspace{0.4cm}
\textbf{Definitions}. $\mathcal{C}(R)$ is called the \emph{ideal-class group of R} and its order is the \emph{class number of R}.

\vspace{0.4cm}
\textbf{Theorem 69}. \textit{$($Finiteness of the class number$)$ The class number of R is finite}.

\vspace{0.4cm}
\emph{Proof}. The proof of this very important theorem is based on the following lemma:

\vspace{0.4cm}
\textbf{Lemma 70}. \textit{Let $N:F\rightarrow \mathbb{Q}$ be the norm mapping. There exits a positive $M\in \mathbb{Z}$ such that for each $\alpha, \beta \in R$ with $\beta\not=0$, there exits $t\in \mathbb{Z}, 1\leq t\leq M$, and $\omega \in R$ such that}
\[
|N(t\alpha-\omega\beta)|<|N(\beta)|.\]

\vspace{0.4cm}
Assume this for now, and let $M$ be the rational integer provided by Lemma 70. Recall from the proof of the claim in the proof of Lemma 51 that $M$! is contained in only finitely many ideals $I_1,\dots,I_k$ of $R$. Let $I\not=\{0\}$ be an ideal of $R$. We will prove that for some $r$, $I\sim I_r$, hence the class number of $R$ is at most $k$.

For each nonzero element $\alpha$ of $I$, $|N(\alpha)|$ is a positive rational integer, so choose $0\not= \beta \in I$ such that $|N(\beta)|$ is a minimum. Lemma 70 implies that for each $\alpha \in I$, there exits a $t\in \mathbb{Z}$ with $1\leq t \leq M$ and $\omega \in R$ such that
\[
|N(t\alpha-\omega\beta)|<|N(\beta)|.\]
Because $t\alpha-\omega\beta \in I$ and $|N(\beta)|$ is a minimum, this inequality implies that $t\alpha-\omega\beta=0$. Hence
\begin{equation*}
(M!)I\subseteq (\beta).\tag{3}\]
Let $J=(M!/\beta)I$. (3) implies that $J$ is an ideal of $R$, and the definition of $J$ implies that
\begin{equation*}
(M!)I=(\beta)J.\tag{4}\]
Because $\beta \in I$, (4) implies that $M!\beta\in (\beta)J$, hence $M!\in J$. Thus $J=I_r$ for some $r$ and (4) now implies that $I\sim J=I_r$.$\hspace{10.4cm} \textrm{QED}$

\emph{Proof of Lemma} 70. The proof we give here follows an ingenious geometric argument due to A. Hurwitz in 1895. It suffices to prove that there exists a positive $M \in \mathbb{Z}$ such that for each $\gamma \in F$, there is a $t\in \mathbb{Z}$, $1\leq t\leq M$, and an $\omega\in R$ such that
\[
|N(t\gamma-\omega)|<1.\]
The conclusion of Lemma 70 follows from this by taking $\gamma=\alpha/\beta$ and using the multiplicativity of $N$.

Let $\{\alpha_1,\dots,\alpha_n\}$ be an integral basis of $F$. Let $\delta \in F$; then 
\[
\delta=\sum_i \delta_i\alpha_i,\ \textrm{with $\delta_i\in \mathbb{Q}$ for all $i$}.\]
We now estimate $N(\delta)$ like so: with $n=[F:\mathbb{Q}]$,
\begin{eqnarray*}
|N(\delta)|&=&\left|\prod_j \left(\sum_i \delta_i \alpha_i^{(j)}\right)\right|\\
&\leq&\prod_j\left(\max_i |\delta_i| \sum_i |\alpha_i^{(j)}|\right)\\
&\leq&\left(\max_i |\delta_i|\right)^n\prod_j \left(\sum_i |\alpha_i^{(j)}|\right)\\
&=&C\left(\max_i |\delta_i|\right)^n,\ C=\prod_j \left(\sum_i |\alpha_i^{(j)}|\right).
\end{eqnarray*}
Choose $m\in \mathbb{Z}$ such that $m>\sqrt[n]{C}$ and then let $M=m^n$.

This is an $M$ that works. In order to see that, let $E^n$ denote Euclidean $n$-space and let
\[
[0, 1]^n=\{(x_1,\dots,x_n)\in E^n:0\leq x_i \leq 1,\ \textrm{for all $i$}\}\]
denote the unit cube in $E^n$.
Define a map $\phi: F\rightarrow [0, 1]^n$ as follows: if $\gamma \in F$, with $\gamma=\sum_i\gamma_i\alpha_i$, write
\[
\gamma_i=a_i+b_i,\ \textrm{where $a_i\in \mathbb{Z}, 0\leq b_i<1$, for all $i$},\]
and then set
\[
\phi(\gamma)=(b_1,\dots,b_n) \in [0, 1]^n.\]

\vspace{0.2cm}
Next, partition $[0, 1]^n$ like so: let$\{I_1,\dots,I_m\}$ be a partition of the unit interval $[0, 1]$ into $m$ pairwise disjoint subintervals each of length $1/m$. For each choice of the indices
\[
(k_1,\dots,k_n)\in \{1,\dots,m\}^n,\]
set
\[
I_{(k_1,\dots,k_n)}=I_{k_1}\times \cdots \times I_{k_n}.\] 
One now easily verifies

\vspace{0.4cm}
\textbf{Lemma 71}. \textit{$(i)$ Each set $I_{(k_1,\dots,k_n)}$ is a subcube of $[0, 1]^n$ of side-length $1/m$.}

\vspace{0.2cm}
$(ii)$ \textit{If $x, y \in I_{(k_1,\dots,k_n)}$ then $|x_i-y_i|\leq 1/m$, $i=1,\dots,n$.}

\vspace{0.2cm}
$(iii)$  \textit{The set $\big\{I_{(k_1,\dots,k_n)}: (k_1,\dots,k_n)\in \{1,\dots,m\}^n \big\}$ is a partition of $[0, 1]^n$.}

\vspace{0.4cm}
Let $\gamma \in F$. Consider the points $\phi(k\gamma), k=1,\dots,m^n+1$ of $[0, 1]^n$. There are $m^n+1$ terms of the sequence of points $(\phi(\gamma),\dots,\phi\big((m^n+1)\gamma\big))$ in $[0, 1]^n$ and there are $m^n$ subcubes in the above partition of $[0, 1]^n$. Hence the pigeon-hole principle implies that at least two terms of this sequence must be in the same subcube, say $\phi(h\gamma)$ and $\phi(l\gamma)$, with $l<h$. Lemma 71$(ii)$ implies that
\begin{equation*}
|\phi(h\gamma)_i-\phi(l\gamma)_i|\leq \frac{1}{m},\ \textrm{for all $i$.}\tag{5}\]

Let $t=h-l \leq m^n=M$. Then
\[
t\gamma=\sum_i h\gamma_i\alpha_i-\sum_i l\gamma_i\alpha_i.\]
Now write
\[
h\gamma_i=a_i+b_i,\ l\gamma_i=c_i+d_i,\ \textrm{where $a_i, c_i \in \mathbb{Z},\ 0\leq b_i<1,\ 0\leq d_i<1.$}\]
Then
\[
\omega=\sum_i(a_i-c_i)\alpha_i\in R,\]
\[
t\gamma-\omega=\sum_i(b_i-d_i)\alpha_i.\]
The above estimate of $|N(\delta)|$ for $\delta\in F$ (with $\delta=t\gamma-\omega$) implies that
\begin{equation*}
|N(t\gamma-\omega)|\leq C\left(\max_i|b_i-d_i|\right)^n.\tag{6}\]
We have from the definition of $\phi$ that
\[
\phi(h\gamma)_i=b_i,\ \phi(l\gamma)_i=d_i,\]
hence (5) implies that
\[
\left(\max_i|b_i-d_i|\right)^n\leq \left(\frac{1}{m}\right)^n,\]
and so (6) and the choice of $m>\sqrt[n]{C}$ implies that
\[
|N(t\gamma-\omega)|\leq Cm^{-n}<1.\]

$\hspace{15cm} \textrm{QED}$

Theorem 45 and Lemma 66 imply that
\begin{quote}
a number ring is a UFD if and only if it has class number 1.
\end{quote}
Hence the class number provides a precise numerical measure of by how much a number ring fails to have unique prime factorization. This highlights the importance of the
\begin{quote}
\emph{Class Number Problem}: given a number ring, calculate its class number.
\end{quote}
The Class Number Problem is one of the most important, and most difficult, problems in all of algebraic number theory, even for class number 1.
\newpage
\begin{center}
\emph{Examples}
\end{center}

(1) \emph{Cyclotomic number fields}

\vspace{0.4cm}
\textbf{Theorem 72}. \textit{$($J. Masley and H. Montgomery, $1976)$ Let $m\in \mathbb{Z},\ m\geq 3,\ \omega_m=e^{2\pi i/m}$. The ring $\mathcal{A}\cap \mathbb{Q}(\omega_m)$ of cyclotomic integers has class number $1$ if and only if $m=3, 4, 5, 7, 8, 9 ,11, 12, 13, 15, 16, 17, 19, 20, 21, 24, 25, 27, 28, 32, 33, 35, 36, 40, 44, 45, 48, 60,$ or $84$.}

\vspace{0.4cm}
\textbf{Corollary 73}. \textit{$($Kummer's Conjecture $(\approx 1857))$ If p is an odd rational prime then $\mathcal{A}\cap \mathbb{Q}(\omega_p)$ has class number $1$ if and only if $p\leq 19$.} 

\vspace{0.4cm}
(2) \emph{Quadratic number fields}

The following theorem was conjectured by Gauss in 1801(in alternative, but equivalent, language) and proved independently by A. Baker and H. M. Stark in 1966:

\vspace{0.4cm}
\textbf{Theorem 74}. \textit{Let $m\in \mathbb{Z}$ be square-free and negative. The ring $\mathcal{A}\cap \mathbb{Q}(\sqrt m)$ of quadratic integers has class number $1$ if and only if $m=-1, -2, -3, -7, -11, -19, -43, -67$, or $-163$.}

\vspace{0.4cm}
Gauss also conjectured that for infinitely many positive square-free $m \in \mathbb{Z}$, $\mathcal{A}\cap \mathbb{Q}(\sqrt m)$ has class number 1. This conjecture, sometimes called the \emph{class-number $1$ problem}, is still open, and is one of the oldest and most famous unsolved problems in algebraic number theory. For an interesting account of the current status of this problem, see H. Cohen [4], Sections 5.5 and 5.10.
\chapter{Ramification and Degree}

Given a number ring $R$ and a prime number $q\in \mathbb{Z}$, we will now study in more detail the prime-ideal factorization of the principal ideal generated by $q$ in $R$. The resulting theory will then be applied, in Chapter 16, to the computation of ideal-class groups and class numbers of quadratic fields.

Let $F$ be a fixed number field, with number ring $R=\mathcal{A}\cap F$.

\vspace{0.4cm}
\textbf{Proposition 75}. \textit{If P is a prime ideal in R then there exits a unique prime $q\in \mathbb{Z}$ such that $P\cap \mathbb{Z}=q\mathbb{Z}$. In particular q is the \textnormal{unique} rational prime contained in P.}

\vspace{0.4cm}
\emph{Proof}. The proof of Proposition 49 implies that $P\cap \mathbb{Z}\not= \{0\}$, and $P\cap \mathbb{Z}\not=\mathbb{Z}$ because $1\not \in P$. Hence $P\cap \mathbb{Z}$ is a prime ideal of $\mathbb{Z}$ and is hence generated in $\mathbb{Z}$ by a unique rational prime $q$.$\hspace{15.3cm} \textrm{QED}$

\vspace{0.2cm}
Let $P$ be a prime ideal of $R$ and let $q$ be the rational prime contained in $P$. Then the ideal $(q)=qR$ generated by $q$ in $R$ is contained in $P$, hence $P$ occurs as a factor in the prime-ideal factorization of $(q)$ in $R$. Let
\[
e=\textrm{the multiplicity of $P$ in $(q)$}.\]

\vspace{0.4cm}
\textbf{Definitions}. The integer $e$ is called the \emph{ramification index of P}. If $e>1$ then $P$ is \emph{ramified} and if $e=1$ then $P$ is \emph{unramified}.

\vspace{0.4cm}
Next we recall some basic facts and notation about quotient rings. Let $I$ be an ideal of a commutative ring  $A$. Under its addition, $A$ is an abelian group and $I$ is a subgroup of $A$. Consider the set $A/I$ of all cosets of $I$ in $A$. If $a\in A$ then $\bar{a}$ will denote the coset $a+I$ of $I$ containing $a$, and $A/I$ becomes a commutative ring when equipped with the addition and multiplication defined by $\bar{a}+\bar{b}=\overline{a+b},\ \bar{a} \bar{b}=\overline{a b}$, called the \emph{quotient ring of A by I}. The quotient map $a\rightarrow \bar{a},\ a\in A$, is a ring homomorphism of $A$ onto $A/I$, called the \emph{quotient homomorphism}. If $a, b\in A$ then we will write $a\equiv b$ mod $I$ if $a-b \in I$. Of course, $a\equiv b$ mod $I$ if and only if $\bar{a}=\bar{b}$ in $A/I$.

\vspace{0.4cm}
\emph{Notation}: if $S$ is a set then $|S|$ will denote the cardinality of $S$.

\vspace{0.4cm}
\textbf{Lemma 76}. \textit{If $I\not=\{0\}$ is an ideal in R then R/I is finite.}

\vspace{0.4cm}
\emph{Proof}. Choose $a\in I\cap\mathbb{Z},\ a>0$. Then $(a) \subseteq I$, hence there is a surjection of $R/(a)$ onto $R/I$, whence it suffices to show that $|R/(a)|$ is finite.

Let $n=[F:\mathbb{Q}]$. We will show that $R/(a)$ has $a^n$ elements. Start with an integral basis $\{\alpha_1,\dots,\alpha_n\}$ of $F$. Consider the set
\[
S=\Big\{\sum_i a_i\alpha_i: a_i \in \mathbb{Z},\ 0\leq a_i<a \Big\}.\]

\emph{Claim}: $S$ is a set of coset representatives of $R/(a)$.

If this is true then $|R/(a)|=|S|=a^n$.

\emph{Proof of the claim}. Let $\alpha=\sum_i z_i\alpha_i \in R$. 
Then there exist $m_i, r_i \in \mathbb{Z},\ 0\leq r_i<a,\ i=1,\dots,n$, such that $z_i=m_ia+r_i,\ i=1,\dots,n$. Hence 
\[
\alpha-\sum_i r_i\alpha_i=\Big(\sum_i m_i\Big)a \in (a)\ \textrm{and}\ \sum_i r_i\alpha_i \in S, \]
and so each coset of $R/(a)$ contains an element of $S$.

Let $\sum_i a_i\alpha_i, \sum_i a_i^{\prime}\alpha_i$ be elements of $S$ in the same coset. Then
\[
\sum_i (a_i-a_i^{\prime})\alpha_i=a\alpha,\ \textrm{for some}\  \alpha \in R.\]
Hence there exists $m_i \in \mathbb{Z}$ such that
\[
\sum_i (a_i-a_i^{\prime})\alpha_i=\sum_i m_ia\alpha_i,\]
and so the linear independence (over $\mathbb{Q}$) of $\{\alpha_1,\dots,\alpha_n\}$ implies that
\[
a_i-a_i^{\prime}=m_ia,\ i=1,\dots,n\]
i.e., $a$ divides $a_i-a_i^{\prime}$ in $\mathbb{Z}$. Because $|a_i-a_i^{\prime}|<a$ for all $i$, it follows that $a_i-a_i^{\prime}=0$ for all $i$. Hence each coset of $R/(a)$ contains exactly one element of $S$. $\hspace{4.5cm} \textrm{QED}$

\vspace{0.2cm}
In conjunction with the ramification index, we will now associate another parameter to a prime ideal of $R$, called the degree of the ideal. To define it, we need two lemmas, the statements of which require the following definition:
 
 \vspace{0.4cm}
 \textbf{Definition}. Let $A$ be a commutative ring, and let
 \[
 C=\{n\in \mathbb{Z}: n>0\ \textrm{and for all $a\in A$, $na=0$}\}.\]
 If $C$ is empty then $A$ \emph{has characteristic $0$}. If $C$ is not empty then the smallest element of $C$ is the \emph{characteristic of A}.
 
  \vspace{0.4cm}
N.B. It is easy to see that if the characteristic of a ring is positive then it must be a prime number.

\vspace{0.4cm}
\textbf{Lemma 77}. \textit{If $K$ is a finite field, i.e., $K$ has only a finite number of elements, then the characteristic $q$ of $K$ is positive and there exist a positive $f\in \mathbb{Z}$ such that $|K|=q^f$. Moreover, $\alpha=\alpha^{q^f}$ for all $\alpha\in K$, and $f$ is the smallest positive rational integer with this property.}

\vspace{0.4cm}
\emph{Proof}. Hungerford [8],  Proposition  V.5.6 and its proof. $\hspace{4.9cm} \textrm{QED}$

\vspace{0.4cm}
\textbf{Lemma 78}. \textit{If P is a prime ideal of R and $q$ is the unique rational prime contained in P, then $R/P$ is a finite field of characteristic $q$.}

\vspace{0.4cm}
\emph{Proof}. $P$ is a maximal ideal of $R$ (Proposition 49). A basic result of ring theory asserts that if $M$ is a maximal ideal in a commutative ring $A$ with identity then $A/M$ is a field (Hungerford [8], Theorem III.2.20$(i)$). Hence $R/P$ is a field, and is finite by Lemma 76.

To see that $R/P$ has characteristic $q$, note first that $P\cap \mathbb{Z}=q\mathbb{Z}$ (Proposition 75), hence there is a natural isomorphism of the field $\mathbb{Z}_q=\mathbb{Z}/q\mathbb{Z}$ into $R/P$ such that the identity element of $\mathbb{Z}_q$ is mapped onto the identity $\bar{1}$ of $R/P$. If we identify $\mathbb{Z}_q$ with its image under this isomorphism then we may assume that $\bar{1}\in \mathbb{Z}_q\subseteq R/P$. Now $\mathbb{Z}_q$ has characteristic $q$, hence $q\bar{1}=0$ and so $q\bar{r}=q\bar{1}\bar{r}=0$ for all $\bar{r}\in R/P$. On the other hand, if $n$ is a positive rational integer such that $n\bar{1}=0$ then $n\bar{x}=0$ for all $\bar{x}\in \mathbb{Z}_q$, and so $q\leq n$. Hence $q$ is the characteristic of $R/P$.$\hspace{11.7cm} \textrm{QED}$

\vspace{0.2cm}
Lemmas 77 and 78 imply that there is a unique positive $f\in \mathbb{Z}$ such that $R/P$ has $q^f$ elements.

\vspace{0.4cm}
 \textbf{Definition}. The integer $f$ is called the \emph{degree of P}.

\vspace{0.4cm}
Let $q \in \mathbb{Z}$ be prime and let
\[
(q)=\prod_{i=1}^g P_i^{e_i}\]
be the prime (ideal) factorization of $(q)$ in $R$, where the $P_i$'s are \emph{distinct} prime ideals. The ideals $P_1,\dots,P_g$ are precisely the prime ideals of $R$ which contain $q$, and $e_i$ is the ramification index of $P_i$. Let $f_i$ be the degree of $P_i$, and let $n=[F:\mathbb{Q}]$. The next theorem gives a remarkable and very useful relationship among the numbers $n, e_i$, and $f_i$.

\textbf{Theorem 79}. \textit{$($the ramification equation$)$.}
\[
n=\sum_i e_i f_i.\]

\vspace{0.4cm}
\emph{Proof}. We need two lemmas, with the first one coming from the general theory of commutative rings.

\vspace{0.4cm}
\textbf{Lemma 80}. \textit{$($Chinese remainder theorem for commutative rings$)$ Let $A$ be a commutative ring with identity, $I_1,\dots,I_k$ ideals of $A$ such that $I_i+I_j=A$ for $i\not=j$. Define the homomorphism} 
\[
\psi: A\rightarrow A/I_1\times \cdots \times A/I_k\]
\textit{by}
\[
\psi(\alpha)=\big(\pi_1(\alpha),\dots,\pi_k(\alpha)\big),\ \alpha \in A,\]
\textit{where $\pi_i: A\rightarrow A/I_i$ is the quotient homomorphism. Then $\psi$ is surjective with kernel $I_1\cdots I_k$, hence $\psi$ induces a natural isomorphism of $A/(I_1\cdots I_k)$ onto $ A/I_1\times \cdots \times A/I_k.$}

\vspace{0.4cm}
\emph{Proof}. Marcus [9], p. 253, and the comment after the proof.$\hspace{4.4cm} \textrm{QED}$

\vspace{0.4cm}
\textbf{Lemma 81}. \textit{If P is a prime ideal of R and if $|R/P|=q^f$ then $|R/(P^e)|=q^{ef}$.}

\vspace{0.4cm}
Assume this lemma for now. In order to prove Theorem 79, start with the prime factorization
\[
(q)=\prod_{i=1}^g P_i^{e_i}\ \textrm{with $f_i=$ the degree of $P_i, i=1,\dots,g.$}\]
Proposition 59$(ii)$ implies that $\gcd\big(P_i^{e_i}, P_j^{e_j}\big)=(1)$ for $i\not=j$, and so if we choose generators $\{\alpha_1, \alpha_2\}$ and $\{\beta_1, \beta_2\}$ for $P_i^{e_i}$ and 
$ P_j^{e_j}$, respectively, then Proposition 59$(iii)$ implies that
\[
P_i^{e_i}+P_j^{e_j}=(\alpha_1, \alpha_2, \beta_1, \beta_2)=\gcd\big(P_i^{e_i}, P_j^{e_j}\big)=(1),\ i\not=j.\]
Hence Lemma 80 implies that
\begin{equation*}
R/(q)\ \textrm{is isomorphic to $R\big/\big(P_1^{e_1}\big)\times \cdots \times R\big/\big(P_g^{e_g}\big)$}.\tag{1}\]
The proof of Lemma 76 implies that
\begin{equation*}
|R/(q)|=q^n.\tag{2}\]
Lemma 81 implies that
\begin{equation*}
\big|R/\big(P_i^{e_i}\big)\big|=q^{e_if_i},\ \textrm{for all $i$}.\tag{3}\]
Hence (1), (2), (3) imply that $q^n=q^{e_1f_1}\cdots q^{e_gf_g}$, and so equating exponents yields the ramification equation.$\hspace{11.6cm} \textrm{QED}$

\vspace{0.2cm}
\emph{Proof of Lemma} 81. We proceed by induction on $e$. The lemma is clearly true for $e=1$, hence we suppose it true for $e-1$, $e>1$. We have that $P^{e-1}/P^e$ is an ideal of $R/P^e$, hence we conclude from the third isomorphism theorem for rings (Hungerford [8], Theorem III.2.12) that
\begin{equation*}
R\big/\big(P^{e-1}\big)\ \textrm{is isomorphic to $\big(R/P^e\big) \big/\big(P^{e-1}/P^e\big)$}.\tag{4}\]
The induction hypothesis implies that
\begin{equation*}
\big|R\big/\big(P^{e-1}\big)\big|=q^{(e-1)f}.\tag{5}\]
\emph{Claim} 1: $\big|P^{e-1}/P^e\big|=q^f$.

If this is true then (4), (5) imply that
\[
q^{(e-1)f}=\frac{\big|R/P^e\big|}{\big|P^{e-1}\big/P^e\big|}=\frac{\big|R/P^e\big|}{q^f},\]
hence
\[
\big|R/P^e\big|=q^{(e-1)f}\cdot q^f=q^{ef}.\]

\vspace{0.2cm}
\emph{Proof of Claim} 1. Observe that $P^e\subsetneqq P^{e-1}$ (otherwise, we cancel $P^{e-1}$ from both sides of $P^{e-1}=P^e$ to obtain $P=(1)$, contrary to the fact that $P$ is proper), and hence choose
\[
\alpha\in P^{e-1}\setminus P^e.\]
We assert that
\begin{equation*}
(\alpha)+P^e=P^{e-1}.\tag{6}\]
In order to see this, let $Q$ be a prime factor of $(\alpha)+P^e$; then
\[
P^e\subseteq (\alpha)+P^e\subseteq Q,\]
hence $Q=P$ and $(\alpha)+P^e$ is a power $P^m$ of $P$. But then
\[
P^e\subseteq P^m=(\alpha)+P^e\subseteq P^{e-1}\ (\alpha\in P^{e-1}!),\]
hence $m=e$ or $e-1$. But $\alpha\not \in P^e$, so $m\not=e$. Hence $m=e-1$, which verifies (6).

Next map $R$ into $P^{e-1}\big/P^e$ by
\[
\phi: \gamma\rightarrow \alpha\gamma+P^e.\]
This is a homomorphism, and (6) implies that it is surjective.

\emph{Claim} 2: kernel of $\phi=P$.

If this is true then $R/P$ is isomorphic to $P^{e-1}\big/P^e$, hence
\[
\big|P^{e-1}/P^e\big|=|R/P|=q^f,\]
which verifies Claim 1.

\emph{Proof of Claim} 2. Note first that $\alpha \in P^{e-1}$ implies that $P\subseteq$ kernel of $\phi$.

Suppose that $\gamma \in$ kernel of $\phi$. Then $\alpha\gamma \in P^e$, hence if $k=$ the multiplicity of $P$ in $(\alpha\gamma)$ then
\begin{equation*}
k\geq e.\tag{7}\]
Now, let
\[
l\ (\textrm{respectively, $m$})=\textrm{the multiplicity of $P$ in $(\alpha)$(respectively, $(\gamma)$)}.\]
We want $m\geq1$, for then $\gamma \in P$ and we are done. In order to get $m\geq1$, we apply the following lemma:

\vspace{0.4cm}
\textbf{Lemma 82}. \textit{If I and J are ideals of R and if i and j are, respectively, the multiplicity of the prime ideal P in I and J, then $i+j$ is the multiplicity s of P in IJ.}

\vspace{0.4cm}
Lemma 82 applied to the product $(\alpha)(\gamma)=(\alpha\gamma)$ implies that
\[
k=l+m.\]
But $l=e-1$ because $\alpha\in P^{e-1}\setminus P^e$, hence (7) implies that
\[
e-1+m=k\geq e,\]
and so $m\geq1$.$\hspace{13.2cm} \textrm{QED}$

\vspace{0.2cm}
\emph{Proof of Lemma} 82. By definition of $i,j,s$, there exist ideals $A, B, C$ such that
\begin{equation*}
\gcd(A, P)=\gcd(B, P)=\gcd(C, P)=(1),\tag{8}\]
\[
P^sC=IJ=P^{i+j}AB.\]
If $s<i+j$ then $C=P^{i+j-s}AB\subseteq P$, and so $P=\gcd(C, P)$, contradicting (8). If $s>i+j$ then $P^{s-i-j}C=AB$, hence $AB\subseteq P$. But $P$ is prime, and so $A\subseteq P$ or $B\subseteq P$ (Proposition 48), and either one of these inclusions also contradicts (8). Hence $s=i+j$. $\hspace{2cm} \textrm{QED}$

\vspace{0.2cm}
We now define and study an important class of number fields for which the ramification equation takes an even nicer form. This will play an important role in our study of ramification of rational primes in quadratic and cyclotomic fields.

\vspace{0.4cm}
\textbf{Definition}. A number field $F$ is a \emph{Galois field} if every embedding of $F$ over $\mathbb{Q}$ maps $F$ into $F$.

\vspace{0.4cm}
\textbf{Definitions}. An automorphism $\sigma$ of a number field $F$ is a \emph{Galois automorphism over $\mathbb{Q}$} if $\sigma(q)=q$ for all $q\in \mathbb{Q}$. The set of all Galois automorphisms of $F$ over $\mathbb{Q}$ is a group under composition of automorphisms, called the \emph{Galois group of $F$ over $\mathbb{Q}$}.

\vspace{0.4cm}
Let $F$ be a number field with Galois group $G$. Because $F\subseteq \mathbf{C}$, every Galois automorphism of $F$ is an embedding of $F$ over $\mathbb{Q}$. Hence if $n$ is the degree of $F$ over $\mathbb{Q}$ then $|G|\leq n$. If $\sigma$ is an embedding of $F$ such that $\sigma(F)\subseteq F$ then $\sigma$ is an injective $\mathbb{Q}$-linear transformation of $F$ into $F$. Since $F$ is a finite-dimensional vector space over $\mathbb{Q}$, it follows that $\sigma(F)=F$, and so $\sigma$ is a Galois automorphism of $F$. Hence if $F$ is a Galois field then $G$ is the set of all embeddings of $F$ over $\mathbb{Q}$. In particular, $|G|=n$, and so whenever $F$ is Galois, the Galois group of $F$ is as large as possible.

\vspace{0.4cm}
\textbf{Proposition 83}. \textit{Let F be a number field. The following statements are equivalent:}

\noindent $(i)$ \textit{F is a Galois field};

\noindent $(ii)$ \textit{if $\theta$ is a primitive element of F then every conjugate of $\theta$ over $\mathbb{Q}$ is in F};

\noindent $(iii)$\textit{ if $\alpha\in F$ then every conjugate of $\alpha$ over $\mathbb{Q}$ is in F.}

\noindent $(iv)$\textit{ if $\alpha\in F$ then every conjugate of $\alpha$ over F is in F.}

\vspace{0.4cm}
\emph{Proof}. Exercise.$\hspace{12cm} \textrm{QED}$

\vspace{0.4cm}
\textbf{Corollary 84}. \textit{Every quadratic number field and every cyclotomic number field is a Galois field}.

\vspace{0.4cm}
\emph{Proof}. If $\mathbb{Q}(\sqrt m)$ is a quadratic field then $\sqrt m$ is obviously a primitive element and $\pm\sqrt m$, the conjugates of $\sqrt m$, are clearly in $\mathbb{Q}(\sqrt m)$. If $\omega_m=e^{2\pi i/m}$ and $\mathbb{Q}(\omega_m)$ is the corresponding cyclotomic field then $\omega_m$ is a primitive element and $\omega_m^k, 1\leq k\leq m,\ \gcd(k, m)=1$, the conjugates of $\omega_m$ (Lemma 11), are clearly contained in $\mathbb{Q}(\omega_m)$. Now apply Proposition 83.$\hspace{15.1cm} \textrm{QED}$

\vspace{0.2cm}
\emph{Example}. $\mathbb{Q}(\sqrt[3]{3})$ is \emph{not} a Galois field:  $\mathbb{Q}(\sqrt[3]{3})=\{q_1+q_2\sqrt[3]{3}+q_3\sqrt[3]{9}: (q_1, q_2, q_3)\in \mathbb{Q}^3\} $ is a subset of the real line, hence $\sqrt[3]{3} e^{2\pi i/3}$ is a conjugate of $\sqrt[3]{3}$ over $\mathbb{Q}$ that is not in  $\mathbb{Q}(\sqrt[3]{3})$.

\vspace{0.4cm}
\textbf{Proposition 85}. \textit{Let F be a number field, G be its Galois group, $R=\mathcal{A}\cap F$}.

\noindent$(i)$ \textit{If $\sigma\in G$ then $\sigma(R)=R$.}

\noindent$(ii)$ \textit{If $I\subseteq R$ and if $\sigma\in G$ then $I$ is an ideal of $R$ if and only if $\sigma(I)$ is an ideal of R.}

\noindent$(iii)$ \textit{If $\sigma\in G$ then I is a prime ideal of $R$ if and only if $\sigma(I)$ is a prime ideal of $R$.}

\vspace{0.4cm}
\emph{Proof}. Let $\alpha \in R$ and let $p$ be the minimal polynomial of $\alpha$ over $\mathbb{Q}$. If $\sigma \in G$ then $\sigma(q)=q$ for all $q\in \mathbb{Q}$, hence $0=\sigma\big(p(\alpha)\big)=p\big(\sigma(\alpha)\big)$. Because $p$ is monic and $p\in \mathbb{Z}[x]$, Lemma 22 implies that $\sigma(\alpha)\in R$. Hence $\sigma(R)\subseteq R$. As $\sigma^{-1} \in G$, it follows that  $\sigma^{-1}(R)\subseteq R$, and so $R\subseteq \sigma(R)$. This verifies $(i)$. Because $(i)$ is valid $\sigma|_R$ is an automorphism of $R$ for all $\sigma \in G$, and $(ii)$ and $(iii)$ are immediate consequences of that.$\hspace{5.8cm} \textrm{QED}$

\vspace{0.2cm}
The next theorem shows how ramification of rational primes and the ramification equation improve when $F$ is a Galois field.

\vspace{0.4cm}
\textbf{Theorem 86}. \textit{If F is a Galois field, $R=\mathcal{A}\cap F$, $q\in \mathbb{Z}$ is prime, $(q)=\prod_{i=1}^g P_i^{e_i}$ is the prime factorization of $(q)$ in R, $f_i=$ the degree of $P_i$, $i=1,\dots,g$, and $n=[F:\mathbb{Q}]$, then }
\[
e_1=\dots=e_g,\ f_1=\dots=f_g\]
\textit{and if $e$ and $f$ are the respective common values of the ramification indices and the degrees of the $P_i$'s then}
\[
n=efg.\]

\vspace{0.4cm}
\emph{Proof}. Because of Theorem 79, we need only prove that $e_1=\dots=e_g$ and $ f_1=\dots=f_g$. This requires

\vspace{0.4cm}
\textbf{Lemma 87}. \textit{If G is the Galois group of F and $i\not=j$ then there exists $\sigma\in G$ such that $\sigma(P_i)=\sigma(P_j)$, i.e., G \textnormal{acts transitively} on the set of of primes ideals of R which contain q.}

\vspace{0.4cm}
Assume this for now. For a fixed $i$, Lemma 87 implies that there exits $\sigma \in G$ such that $\sigma(P_1)=P_i$. Then $R/P_1$ and $R\big/\sigma(P_1)=R/P_i$ are isomorphic, hence 
\[
f_1=\log_q|R/P_1|=\log_q|R/P_i|=f_i,\ \textrm{for all $i$},\]
and so all the degrees $f_i$ are the same.

As a consequence of Proposition 85$(i)$, $\sigma|_R$ is an automorphism of $R$, and so it follows easily from Proposition 85$(ii)$, $(iii)$ that $\sigma|_R$ preserves multiplicities in the prime factorization of ideals in $R$. As $\sigma((q))=(q)$, it hence follows that $e_i=$ the multiplicity of $P_i$ in $(q)=$ the multiplicity of $\sigma(P_1)$ in $\sigma((q))=$ the multiplicity of $P_1$ in $(q)=e_1$ for all $i$, i.e., all indices $e_i$ are the same. $\hspace{14.4cm} \textrm{QED}$

\vspace{0.2cm}
\emph{Proof of Lemma} 87. Suppose this is false. Then there exit $i\not=j$ such that $P_i\not \in \{\sigma(P_j): \sigma \in G\}$. By reindexing if necessary, we may assume that $\sigma_1(P_j),\dots,\sigma_m(P_j),\ m\leq n$, are the distinct images of $P_j$ under the elements of $G$, and so $P_j, \sigma_1(P_j),\dots,\sigma_m(P_j)$ are distinct prime ideals, hence the greatest common divisor of any two of them is $(1)$, hence the sum of any two of them is $R$ (from the proof of Theorem 79), hence Lemma 80 implies that if $\pi:R\rightarrow R/P_i,\ \pi_k: \rightarrow R/\sigma_k(P_j) $ are the quotient homomorphisms then
\[
\psi:\gamma\rightarrow \big(\pi(\gamma), \pi_1(\gamma),\dots,\pi_m(\gamma)\big)\]
is a surjection of $R$ onto
\[
R/P_i\times R/\sigma_1(P_j)\times \cdots \times R/\sigma_m(P_j),\]
and so there exist $\gamma \in R$ such that $\psi(\gamma)=(\bar{0}, \bar{1},\dots,\bar{1})$, i.e.,
\begin{equation*}
\gamma\in P_i,\ \gamma-1\in \sigma_k(P_j),\ k=1\dots,m.\tag{9}\]

Now let $N:F\rightarrow \mathbb{Q}$ be the norm map. As $F$ is a Galois field, it follows that $G$ is the set of embeddings of $F$ over $\mathbb{Q}$, hence
\[
N(\gamma)=\prod_{\sigma\in G} \sigma(\gamma)= \gamma\left( \prod_{\sigma\not=\textrm{the identity map on $R$}} \sigma(\gamma)\right).\]
We have that $\gamma \in P_i$ by (9) and $\prod_{\sigma\not=\textrm{the identity map on $R$}} \sigma(\gamma) \in R$ by Proposition 85$(i)$, hence
\[
N(\gamma)\in P_i\cap \mathbb{Z}=q\mathbb{Z}=P_j\cap \mathbb{Z}.\]
Because $P_j$ is prime, there hence exists $\sigma \in G$ such that $\sigma(\gamma)\in P_j$, i.e., $\gamma \in \sigma^{-1}(P_j)$. But (9) implies that $\gamma-1\in \sigma^{-1}(P_j)$, hence $1\in  \sigma^{-1}(P_j)$, and this is impossible because $\sigma^{-1}(P_j)$ is a proper ideal.$\hspace{13.5cm} \textrm{QED}$

\vspace{0.4cm}
\textbf{Definitions}. Let $F$ be a Galois field, $q \in \mathbb{Z}$ a prime, $n=[F:\mathbb{Q}]$,
\[
(q)=\prod_{i=1}^g P_i^e
\]
the prime factorization of $(q)$ in $R=\mathcal{A}\cap F$. The prime $q$ is \emph{ramified in $F$} if $e>1$, \emph{completely ramified in F} if $e>1$ and $g=1$, \emph{unramified in F} if $e=1$, \emph{split in F} if $g=n$ (and hence also $e=1$), and \emph{inertial in F} if $g=e=1$.

\vspace{0.4cm}
Thus $q$ is completely ramified in $F$ if and only if $(q)$ is a power $>1$ of a single prime ideal, $q$ is unramified in $F$ if and only if $(q)$ factors into a product of distinct prime ideals, $q$ splits in $F$ if and only if $(q)$ factors into the maximum possible number $n$ of distinct prime ideals, and $q$ is inertial in $F$ if and only if $(q)$ is a prime ideal in $R$ if and only if $q$ is a prime element of $R$.

N.B. Let $S$ and $R$ be number rings with $S\subseteq R$. If $P$ is a prime ideal of $S$ then the set $PR=\big\{\sum_i x_iy_i: x_i\in P, y_i \in R\big\}$ is an ideal of $R$ and hence factors into a product of prime ideals of $R$. The distinct prime-ideal factors of this product consist precisely of the set of prime ideals of $R$ which contain $P$, and these prime ideals of $R$ are said to \emph{lie over P}. One can develop a ramification theory for this more general situation that is in complete analogy with what we have done in Chapter 13 (there we have taken $S=\mathbb{Z}$, the simplest, yet one of the most important, special cases of this general situation). For the development of ramification and degree in the general set-up alluded to here, see Marcus [9], Chapters 3 and 4.
\chapter{Ramification in Cyclotomic Number Fields}

The statement of the main theorem of this chapter requires a bit of elementary number theory. Recall

\vspace{0.4cm}
\textit{Fermat's Little Theorem. If $r\in \mathbb{Z}$ is prime and $a\in \mathbb{Z}$ then}
\[
a^r\equiv a\ \textnormal{mod}\ r.\]

\vspace{0.4cm}
Suppose that $a\in \mathbb{Z}$ and the rational prime $r$ does not divide $a$; then by Fermat's Little Theorem, 
\[
a^{r-1}\equiv 1\ \textnormal{mod}\ r.\]
Hence there is a smallest positive $f\in \mathbb{Z}$ such that $a^f\equiv 1\ \textnormal{mod}\ r$. 

\vspace{0.4cm}
\textbf{Definitions}. The number $f$ is called the \emph{order of a with respect to r}. If the order of $a$ with respect to $r$ is $r-1$ then $a$ is called a \emph{primitive root of r}.

\vspace{0.4cm}
N.B. Because the order of $a$ with respect to $r$ is the order of the element $a+r\mathbb{Z}$ in the mulitplicative group of the field $\mathbb{Z}/r\mathbb{Z}$, it follows that the order of $a$ with respect to $r$ always divides $r-1$.

Let $p$ be a fixed odd rational prime, $\omega=e^{2\pi i/p}$, $F=\mathbb{Q}(\omega),\ R=\mathcal{A}\cap F$. Corollary 84 implies that $F$ is Galois, hence, by way of Theorem 86, if $q\in \mathbb{Z}$ is prime then $(q)$ has a prime factorization of the form
\[
(q)=\prod_{i=1}^g P_i^e\]
\[
\textrm{degree of $P_i=d$, for all $i$, and $edg=p-1$}.\]

The following theorem describes precisely how the rational prime $q$ ramifies in $F$.

\vspace{0.4cm}
\textbf{Theorem 88}. \textit{$(i)$ Suppose that $q\not=p$, and let $f=$ the order of $q$ with respect to p. Then $q$ is unramified in F, each prime factor $P_i$ of $(q)$ has degree f, and the number $g$ of prime factors of $(q)$ is $(p-1)/f$.}

\textit{$(ii)$ The prime p is completely ramified in F; in fact $(1-\omega)$ is a prime ideal in R of degree $1$ and $(p)=(1-\omega)^{p-1}.$}

\vspace{0.4cm}
\textbf{Corollary 89}. \textit{$(i)$ The prime q splits in F if and only if p divides $q-1$ $($in $\mathbb{Z})$.}

\textit{$(ii)$ The prime q is inertial in F if and only if q is a primitive root of p.}

\vspace{0.4cm}
\emph{Proof}. Theorem 88 implies that $q$ splits in $F$ if and only if $f=1$ if and only if $p$ divides $q-1$ (in $\mathbb{Z})$, and $q$ is inertial in $F$  if and only if $1=g=(p-1)/f$ if and only if $f=p-1$.$\hspace{13.8cm} \textrm{QED}$

\vspace{0.2cm}
\emph{Proof of Theorem} 88. 
$(i)$ We show first that $q$ is unramified. To do this we need

\vspace{0.4cm}
\textbf{Lemma 90}. \textit{If $\alpha\in R$ then $\alpha^{q^f}-\alpha\in (q)$.}

\vspace{0.4cm}
Assume this for now. Suppose that $q$ is ramified in $F$. Then
\[
(q) \subseteq P^2\ \textrm{for some prime ideal $P$}.\]
Let $\alpha\in P\setminus P^2$. Lemma 90 implies that $\alpha^{q^f}-\alpha\in (q)$, and so
\[
\alpha^{q^f}\equiv \alpha \ \textnormal{mod}\ P^2.\]
But $q^f\geq 2$ and $\alpha \in P$, hence $\alpha^{q^f} \in P^2$, and so $\alpha\in P^2$, contrary to its choice.

We prove next that the degree $d$ of each prime factor $P_i$ of $(q)$ is $f$. Observe first that $R/P_i$ is a finite field of order $q^d$. Lemma 77 implies that
\[
\alpha^{q^d}\equiv \alpha \ \textnormal{mod}\ P_i,\ \textrm{for all $\alpha \in R$},\] 
and $d$ is the smallest positive rational integer with this property. But $(q)\subseteq P_i$, so Lemma 90 implies that
\[
\alpha^{q^f}\equiv \alpha \ \textnormal{mod}\ P_i,\ \textrm{for all $\alpha \in R$},\] 
hence
\begin{equation*}
d \leq f.\tag{1}\]
On the other hand, the proof of Proposition 20 implies that
\begin{equation*}
p=\prod_{j=1}^{p-1} (1-\omega^j),\tag{2}\]
hence, taking cosets in $R/P_i$,
\[
\bar{p}=\prod_{j=1}^{p-1} (\bar{1}-\bar{\omega}^j).\]
Now $\bar{p}=\bar{0}$ implies $p\in P_i$. But $q\in P_i$ and $q$ is the only rational prime contained in $P_i$. As $p\not=q$, we conclude that $\bar{p}\not=\bar{0}$, hence $\bar{\omega}^j\not=\bar{1},\ j=1,\dots,p-1$, and so the powers $\bar{\omega}^j,\ j=0, 1,\dots,p-1$, are all distinct. Thus $\{\bar{1}, \bar {\omega}_1,\dots,\bar{\omega}^{p-1}\}$ is a subgroup of order $p$ in the multiplicative group of $R/P_i$. The latter group has order $q^d-1$, hence Lagrange's Theorem implies that $q^d\equiv 1$ mod $p$. Hence the definition of $f$ implies that
\begin{equation*}
f \leq d,\tag{3}\]
and so (1) and (3) imply that $d=f$. Hence, by Theorem 86, $g=(p-1)/f$.

\emph{Proof of Lemma} 90. Note from Proposition 25 that $R=\mathbb{Z}[\omega]$, so if $\alpha\in R$ then we find $a_i\in \mathbb{Z}$ such that $\alpha=\sum_i a_i\omega^i$. Fermat's Little Theorem implies that
\begin{equation*}
a_i^{q}\equiv a_i \ \textnormal{mod}\ (q),\ \textrm{for all $i$}.\tag{4}\]
The quotient ring $R/(q)$ has characteristic $q$, hence from the binomial theorem, $(\bar{x}+\bar{y})^q=\bar{x}^q+\bar{y}^q$ for all $\bar{x}, \bar{y}\in R/(q)$, and so
\begin{eqnarray*}
\alpha^q&=& \Big(\sum_i a_i\omega^i\big)^q\\
&\equiv&\sum_i a_i^q\omega^{iq} \ \textnormal{mod}\ (q)\\
&\equiv&\sum_i a_i\omega^{iq} \ \textnormal{mod}\ (q),\ \textrm{by $(4)$}.
\end{eqnarray*}
Repeating this calculation $f$ times yields
\[
\alpha^{q^f}\equiv \sum _i a_i(\omega^{q^f})^i \ \textnormal{mod}\ (q).\]
But $q^f \equiv 1$ mod $p$ hence $\omega^{q^f}=\omega$ ($\omega^p=1!$), and so
\[
\alpha^{q^f}\equiv \sum _i a_i \omega^i \ \textnormal{mod}\ (q) \equiv \alpha \ \textnormal{mod}\ (q).\]

$\hspace{15cm} \textrm{QED}$

$(ii)$ Let\[
u_i=\frac{1-\omega^i}{1-\omega},\ i=1,\dots,p-1.\]

\emph{Claim}: $u_i$ is a unit in $R,\ i=1,\dots,p-1$.

If this is so, then (2) implies that
\[
p=(1-\omega)^{p-1}\prod_i u_i=(1-\omega)^{p-1}\times\ \textrm{(a unit in $R$)},\]
hence
\[
(p)=(1-\omega)^{p-1}.\]
Now let $(1-\omega)=\prod_{i=1}^h Q_i^{k_i}$ be the prime factorization of $(1-\omega)$ in $R$. Then
\[
(p)=\prod_{i=1}^h Q_i^{k_i(p-1)}.\]
Theorem 86 implies that $k_1(p-1)=\cdots=k_h(p-1)$ hence $k_1=\cdots=k_h=k$, say, and the $Q_i$'s all have the same degree $m$. Then the ramification equation for $(p)$ implies that $p-1=mhk(p-1)$, hence $m=h=k=1$, and so $(1-\omega)=Q_1$ is prime of degree 1.

\emph{Proof of the claim}. If $i\in \mathbb{Z}$ and $1\leq i<p$ then $p$ does not divide $i$ in $\mathbb{Z}$, and so there exists $j\in \mathbb{Z}$, $1\leq j<p$ such that $ij\equiv 1$ mod $p$, hence
\[
\frac{1}{u_i}=\frac{1-\omega}{1-\omega^i}=\frac{1-\omega^{ij}}{1-\omega^i}=\sum_{k=0}^{j-1} \omega^{ik}\in R.\]

$\hspace{15cm} \textrm{QED}$

Let $m\in \mathbb{Z},\ m\geq 3,\ \omega_m=e^{2\pi i/m}$. How do rational primes ramify in $\mathbb{Q}(\omega_m)$? The answer requires the following interesting and important generalization of Fermat's Little Theorem:

\vspace{0.4cm}
\textit{Euler's Theorem. Let $m, k\in \mathbb{Z},\ m\geq 2,$ $k$ and $m$ relatively prime. If $\varphi(m)$ is the value of Euler's totient at $m$ then}
\[
k^{\varphi(m)}\equiv 1\ \textnormal{mod}\ m.\]

\vspace{0.4cm}
Hence if $k$ and $m$ are relatively prime then there exist a smallest positive rational integer $f$ such that $k^f\equiv 1$ mod $m$.

\vspace{0.4cm}
\textbf{Definitions}. $f$ is called the \emph{ order of k with respect to m}. If $f=\varphi(m)$ then $k$ is a \emph{primitive root of m}.

\vspace{0.4cm}
N.B. The order of $k$ with respect to $m$ is the order of the element $k+m\mathbb{Z}$ in the multiplicative group of the ring $\mathbb{Z}/m\mathbb{Z}$, and since this group has order $\varphi(m)$, it follows that the order of $k$ with respect to $m$ always divides $\varphi(m)$ (in $\mathbb{Z}$).

\vspace{0.4cm}
\textbf{Theorem 91}. \textit{Let $q\in \mathbb{Z}$ be prime, $R_m=\mathcal{A}\cap \mathbb{Q}(\omega_m)$, }
\[
(q)=\prod_{i=1}^g P_i^e\]
\textit{the prime factorization of $(q)$ in $R_m,\ d=$ the degree of $P_i$, and let $m=q^kn$, with $q$ and n relatively prime $($in $\mathbb{Z})$. Then $e=\varphi(q^k)$,}
\[
d=\left\{\begin{array}{cc}1,\ \textrm{if $n=1$,}\\
\ \textrm{order of $q$ with respect to $n$, if $n>1$,}
\\\end{array}\right.
\]

\vspace{0.2cm}
\[
g=\left\{\begin{array}{cc}1,\ \textrm{if $n=1$,}\\
\ \textrm{$\varphi(n)$/(order of $q$ with respect to $n$), if $n>1$,}
\\\end{array}\right.
\]

\vspace{0.4cm}
\emph{Proof}. See Marcus [9], , Chapter 3, Theorem 26.$\hspace{6.2cm}\textrm{QED}$

\vspace{0.4cm}
\textbf{Corollary 92}. \textit{$(i)$ q is ramified in $\mathbb{Q}(\omega_m)$ if and only if $q|m$ in $\mathbb{Z}$.}

\textit{$(ii)$ q is completely ramified in $\mathbb{Q}(\omega_m)$ if and only if $n=1$, i.e., either m is a power of q or $n>1$ and q is a primitive root of n.}

\textit{$(iii)$ q splits in $\mathbb{Q}(\omega_m)$ if and only if $m|(q-1)$ in $\mathbb{Z}$.}

\textit{$(iv)$ q is inertial in $\mathbb{Q}(\omega_m)$ if and only if q does not divide m in $\mathbb{Z}$ and q is a primitive root of m.}

\chapter{Ramification in Quadratic Number Fields}

Let $m$ be a fixed square-free rational integer, let $F=\mathbb{Q}(\sqrt m),\ R=\mathcal{A}\cap F$. Let $d=$ the discriminant of $F$, which by Proposition 23 is $m$ (respectively, $4m$) if $m\equiv 1$ mod 4 (respectively, $m\not \equiv 1$ mod 4). Let $p\in \mathbb{Z}$ be a fixed prime.

Corollary 84 implies that $F$ is a Galois field, hence if $e$ and $f$ are the ramification index and the degree of the prime factors of $(p)$ in $R$ and $g$ is the number of prime factors then $efg=[F:\mathbb{Q}]=2$, hence there are only three possibilities:

$e=2,\ f=g=1$, i.e., $p$ is completely ramified in $F$,

$e=f=1,\ g=2$, i.e., $p$ splits in $F$, and,

$e=g=1,\ f=2$, i.e., $p$ is inertial in $F$.

The ramification of $p$ in $F$ will depend on whether $p$ is even or odd, so we treat each case separately.

\vspace{0.4cm}
\textbf{Theorem 93}. \textit{Suppose that p is odd}.

\textit{$(i)$ If p does not divide d in $\mathbb{Z}$ and there exists $a\in \mathbb{Z}$ such that $a^2\equiv m\ \textnormal{mod}\ p$, i.e., m is a quadratic residue of p, then $(p,\ a\pm \sqrt m)$ are distinct prime ideals of R and}
\[
(p)=(p,\ a+\sqrt m)(p,\ a-\sqrt m),\]
\textit{i.e., p splits in F}.

\textit{$(ii)$ If p does not divide d in $\mathbb{Z}$ and there does not exist $a\in \mathbb{Z}$ such that $a^2\equiv m\ \textnormal{mod}\ p$, i.e., m is a quadratic non-residue of p, then p is inertial in F}. 

\textit{$(iii)$ If $p|d$ in $\mathbb{Z}$ then $(p, \sqrt m)$ is a prime ideal in R and}
\[
(p)=(p, \sqrt m)^2,\]
\textit{i.e., p completely ramifies in F}.

\vspace{0.4cm}
\emph{Proof}. We first prove that if $(p)$ factors as in $(i)$ or $(iii)$ then the factors must be nonzero and proper. They are all clearly nonzero, and since $(p)$ is proper the factor of $(p)$ in $(iii)$ must be proper. Suppose that $(p)$ factors as in $(i)$ and $(1)=(p,\ a+\sqrt m)$. Then $a-\sqrt m \in (p)$ hence $p|a$ in $\mathbb{Z}$. Since  $a^2\equiv m\ \textnormal{mod}\ p$, it follows that $p|m$, and so $p|d$, in $\mathbb{Z}$, contrary to hypothesis. The same reasoning implies that $(p,\ a-\sqrt m)\not=(1)$.

Next we prove that if $(p)$ is a product of two not-necessarily-distinct, nonzero, proper ideals, then the factors must all be prime. This will show that the primality of the ideals in the conclusions of $(i)$ and $(iii)$ will follow from the validity of the other conclusions of $(i)$ and $(iii)$.

Suppose that $(p)=J_1J_2,\ \{0\}\not= J_i\not= (1)$ ideals of $R,\ i=1, 2$. Let $e$ be the ramification index of the prime factors of $(p)$ and let $g$ be the number of distinct ideals in the prime factorization of $(p)$. Theorem  86 implies that $1\leq e \leq 2$,  $1\leq g \leq 2$, hence $(p)$, and thus also $J_1$ and $J_2$ have at most two distinct prime factors. If $(p)$ has exactly one prime factor, say $Q$, then $Q$ is also the only prime factor of $J_1$ and $J_2$. Hence if $e_i=$ the multiplicity of $Q$ in $J_i,\ i=1, 2$ then (Lemma 82) $e_1+e_2=e\leq 2$, and because $J_i\not=(1)$, we have that $e_i\not=0,\ i=1, 2$. Hence $e_1=e_2=1$ and so $J_1=Q=J_2$ is prime. If $(p)$ has exactly two prime factors, say $Q_1$ and $Q_2$, then $e=1$. Now the set of prime factors of both $J_1$ and $J_2$ are contained in $\{Q_1, Q_2\}$ (otherwise $(p)$ would have at least three prime factors), and if $e_{ij}=$ the multiplicity of $Q_i$ in $J_j$, then $e_{i1}+e_{i2}=e=1,\ i=1, 2.$ Hence $\{e_{i1}, e_{i2}\}=\{0, 1\},\ i=1, 2,$ and so $\big\{J_1, J_2\big\}=\big\{Q_1, Q_2\big\}$. Hence $J_1$ and $J_2$ are both prime.

$(i)$ We have that
\begin{eqnarray*}
(*)\hspace{1cm} (p,\ a+\sqrt m)(p,\ a-\sqrt m)&=&\left(p^2,\ p(a+\sqrt m),\ p(a-\sqrt m),\ p\ \frac{a^2-m}{p}\right)\\
&=&(p)\left(p,\ a+\sqrt m,\ a-\sqrt m,\ \frac{a^2-m}{p}\right).
\end{eqnarray*}
Let $I$= the second factor of the product on the right-hand side of the second equation in $(*)$. Then $p$ and $2a=a+\sqrt m+a-\sqrt m$ are in $I$. Because $p$ is odd and $p$ does not divide $a$ in $\mathbb{Z}$, it follows that $p$ and $2a$ are relatively prime, and so there exit $x, y \in \mathbb{Z}$ such that $xp+2ay=1$. Hence $1\in I$, and so from $(*)$ it follows that
\begin{equation*}
(p)=(p,\ a+\sqrt m)(p,\ a-\sqrt m).\tag{$**$}\]
Moreover, $(p,\ a+\sqrt m)\not=(a-\sqrt m)$; otherwise, this ideal contains both $p$ and $2a$ hence $(p,\ a+\sqrt m)=(1)=(a-\sqrt m)$ and so by $(**)$, $(p)=(1)^2=(1)$, which is impossible.

$(ii)$ Let $Q$ be a prime factor of $(p)$.

\emph{Claim}: degree of $Q=2$.

If this is true then $e=g=1$, hence $(p)=Q$ is prime.

\emph{Proof of the claim}. Suppose by way of contradiction that the degree of $Q$ is 1. Then $R/Q$ has $p$ elements. The map $\mathbb{Z}/p\mathbb{Z}\rightarrow R/Q$ defined by $a+p\mathbb{Z}\rightarrow a+Q,\ a\in \mathbb{A}$, is injective, and $\mathbb{Z}/p\mathbb{Z}$ has $p$ elements. Hence this map is surjective, i.e., if $\alpha\in R$ then there exits $a\in \mathbb{Z}$ such that $a\equiv \alpha$ mod $Q$. Take $\alpha=\sqrt m$ and square both sides to obtain $a\in \mathbb{Z}$ such that $a^2\equiv m$ mod $Q$, hence
\[
a^2-m\in Q\cap \mathbb{Z}=p\mathbb{Z},\]
and so $a^2\equiv m\ \textnormal{mod}\ p$, contrary to hypothesis.

$(iii)$ We have that $p|m$ in $\mathbb{Z}$ since $p|(d=m$ or $4m)$ in $\mathbb{Z}$. Hence 
\[
(p, \sqrt m)^2=(p^2,\ p\sqrt m,\ m)=(p)(p, \sqrt m,\ m/p).\]
Let $I= (p, \sqrt m,\ m/p).$ If $p$ divides $m/p$ in $\mathbb{Z}$ then $p^2|m$ in $\mathbb{Z}$, impossible since $m$ is square-free. Hence $p$ and $m/p$ are relatively prime in $\mathbb{Z}$, and so $I=(1)$ as in the proof of $(i)$, hence $(p)=(p, \sqrt m)^2$.$\hspace{12.8cm}\textrm{QED}$

\vspace{0.2cm}
If $d$ is odd then $d=m$, hence $m\equiv 1$ mod 4, and so $m\equiv 1$ or 5 mod 8. If $d$ is even then either $m$ is even, hence $m\equiv$ 2 mod 4 or $d=4m$, hence $m\equiv$ 2 or 3 mod 4. The cases in the following theorem are thus exhaustive and mutually exclusive.

\vspace{0.4cm}
\textbf{Theorem 94}. \textit{Suppose that $p=2$}.

 \textit{$(i)$ If d is odd and $m\equiv 1\ \textnormal{mod}\ 8$ then}
 \[
 \left(2,\ \frac{1\pm \sqrt m}{2}\right)\]
\textit{are distinct prime ideals in R and}
\[
(2)= \left(2,\ \frac{1+ \sqrt m}{2}\right)\left(2,\ \frac{1-\sqrt m}{2}\right),\]
\textit{i.e., $2$ splits in F}.    

\textit{$(ii)$ If d is odd and $m\equiv 5\ \textnormal{mod}\ 8$ then $2$ is inertial in F}.

\textit{$(iii)$ If d is even and $m\equiv 2\ \textnormal{mod}\ 4$ then $(2, \sqrt m)$ is prime in R and}
\[
(2)=(2, \sqrt m)^2.\]
\textit{If d is even and $m\equiv 3\ \textnormal{mod}\ 4$ then $(2, 1+\sqrt m)$ is prime in R and}
\[
(2)=(2, 1+\sqrt m)^2.\]
\textit{Hence if d is even then $2$ is completely ramified in F}.

\vspace{0.4cm}
\emph{Proof}. As in the proof of Theorem 93, primality of all ideals appearing in the statements of $(i)$ and $(iii)$ follows from the validity of the other conclusions of those statements.

$(i)$ Note first that Proposition 23 implies that $\displaystyle{\frac{1\pm \sqrt m}{2}} \in R$. We have that
\[
\left(2,\ \frac{1+ \sqrt m}{2}\right)\left(2,\ \frac{1-\sqrt m}{2}\right)=(2)\left(2,\ \frac{1+ \sqrt m}{2},\ \frac{1-\sqrt m}{2},\ \frac{1-m}{8}\right).\]
The second factor of the product on the right-hand side is $(1)$ since it contains $\displaystyle{\frac{1+ \sqrt m}{2}}+\displaystyle{\frac{1- \sqrt m}{2}}=1$, hence
\[
\left(2,\ \frac{1+ \sqrt m}{2}\right)\left(2,\ \frac{1-\sqrt m}{2}\right)=(2).\]
 
 \vspace{0.2cm}
\noindent Moreover, the ideals on the left-hand side of this equation are distinct; otherwise this ideal contains 1, hence $(1)=(1)^2=(2)$, a contradiction.

$(ii)$ Let $Q$ be a prime factor of $(2)$. We need only prove that $Q$ has degree 2, for then $(2)=Q$ is prime.

Suppose that the degree of $Q$ is 1. Then as in the proof of Theorem 93$(ii)$, there is an $a \in \mathbb{Z}$ such that $a\equiv \displaystyle{\frac{1+ \sqrt m}{2}}$ mod $Q$. Now $\displaystyle{\frac{1+ \sqrt m}{2}}$ is a root of the polynomial
\[
x^2-x+\frac{1-m}{4},\]
hence 
\[
a^2-a+\frac{1-m}{4} \equiv 0\ \textnormal{mod}\ Q,\]
hence
\[
a^2-a+\frac{1-m}{4} \in \mathbb{Q}\cap \mathbb{Z}=2\mathbb{Z},\]
and so
\[
a^2-a+\frac{1-m}{4}\ \textrm{is even}.\]
But $a^2-a=a(a-1)$ is even, hence $\displaystyle{\frac{1-m}{4}}$ is even, i.e., $m\equiv 1$ mod 8, contrary to hypothesis.

$(iii)$ This follows by arguments very similar to the ones used in the proof of Theorem 93$(iii)$, and so we leave the details as an exercise.$\hspace{7cm}\textrm{QED}$

\chapter{Computing the Ideal-Class Group in Quadratic Fields}

Ramification theory for quadratic number fields, when combined with some additional mathematical technology, can be used effectively to compute ideal-class groups and class numbers for those fields. We illustrate how things go with three examples. But first, the additional technology that is required.

Let $R=\mathcal{A}\cap F$ be a fixed number ring with $n=[F:\mathbb{Q}]$.

\vspace{0.4cm}
\textbf{Definition}. If $I$ is an ideal of $R$ then the \emph{norm $N(I)$ of I} is defined to be the cardinality of $R/I$.

\vspace{0.4cm}
N.B. Lemma 76 implies that $N(I)<+\infty$ for all nonzero ideals $I$ of $R$.

\vspace{0.4cm}
\emph{Notation}. We let $E$ denote the set of real numbers.

\vspace{0.4cm}
\textbf{Definition}. If $\sigma$ is an embedding of $F$ over $\mathbb{Q}$ then $\sigma$ is \emph{real} (respectively, \emph{complex}) if $\sigma(F)\subseteq E$ (respectively, $\sigma(F)\nsubseteq E$).

\vspace{0.4cm}
Let
\[
r=\textrm{the number of real embeddings of $F$},\]
\[
c=\textrm{the number of complex embeddings of $F$.}\]

\noindent We have $n=r+c$, and we claim that $c$ is even. In order to see that, let $\theta$ be a primitive element of $F$ over $\mathbb{Q}$, $\theta^{(1)},\dots,\theta^{(n)}$ the (distinct) conjugates of $\theta$ over $\mathbb{Q}$. Proposition 14 implies that if $\{\sigma_1,\dots,\sigma_n\}$ are the embeddings of $F$ over $\mathbb{Q}$ then
\[
\sigma_i(F)=\{p(\theta^{(i)}): p\in \mathbb{Q}_n[x]\},\ \textrm{for each $i$},\]
and because each $p\in \mathbb{Q}_n[x]$ has real coefficients it follows that
\[
\sigma_i\ \textrm{is real if and only if $\theta^{(i)}$ is real, for each $i$}.\]
Hence $r$ is the number of real conjugates of $\theta$ and $c$ is the number of non-real conjugates. But $\theta^{(1)},\dots,\theta^{(n)}$ are the roots of the minimal polynomial of $\theta$ over $\mathbb{Q}$, which polynomial has real coefficients, and so its non-real roots occur in (complex) conjugate pairs. Hence $c$ is even. 

Let $s=c/2$. The next lemma contains the additional mathematical tools that we need for the calculation of ideal-class groups; we will first use it in those calculations and then we will prove (most of) it. 

\vspace{0.4cm}
\textbf{Lemma 95}. \textit{$(i)$ If I and J are ideals of R, then}
\[
N(IJ)=N(I)N(J).\]

\textit{$(ii)$ If $0\not= \alpha \in R$ then the norm of the principal ideal generated by $\alpha$ is $|N(\alpha)|$.}

\textit{$(iii)$ If d is the discriminant of F then in each ideal class of R there is an ideal I such that}
\[
N(I)\leq \lambda=\frac{n!}{n^n}\left(\frac{4}{\pi}\right)^s\sqrt{|d|}.\]

\vspace{0.4cm}
The constant $\displaystyle{\frac{n!}{n^n}\left(\frac{4}{\pi}\right)^s\sqrt{|d|}}$ is called \emph{Minkowski's constant}, and arises in the study of the geometry of numbers. We will get a glimpse of this very important subject when we study the structure of the group of units of a number ring in Chapter 17. 

\vspace{0.4cm}
\emph{Example} 1

Let $F=\mathbb{Q}(\sqrt 2\ ),\ R=\mathcal{A} \cap F=\mathbb{Z}+\sqrt 2\ \mathbb{Z}.$ Then $n=2,\ s=0,\ d=8$ (from the calculation of the discriminant on p. 38), hence the value of Minkowski's constant $\lambda$ in Lemma 95$(iii)$ is
\[
\frac{2!}{2^2}\sqrt 8<2,\]
and so by Lemma 95$(iii)$, every ideal class of $R$ contains an ideal $I$ with $N(I)\leq 1$ hence $|R/I|=N(I)=1$, hence $I=(1)$. Conclusion: $R$ has only one ideal class, the principal class, and so $R$ has class number 1.

\vspace{0.4cm}
\emph{Example} 2

Let $F=\mathbb{Q}(\sqrt {-5}\ ),\ R=\mathcal{A} \cap F=\mathbb{Z}+\sqrt {-5}\ \mathbb{Z}.$ Then $n=2,\ s=1,\ d=-20$, so $\lambda=\displaystyle{\frac{4\sqrt 5}{\pi}}<3$, hence every ideal class of $R$ contains an ideal $I$ such that $N(I)$ is either 1 or 2.

If $N(I)=1$ then $I=(1)$. Suppose that $N(I)=2$. Then the additive group of $R/I$ has order 2, and so
\[
2(\alpha+I)=I,\ \textrm{for all $\alpha \in R$,}\]
and taking $\alpha=1$, we obtain $2 \in I$. Hence all of the prime factors of $I$ must contain 2, so we factor the ideal $(2)$ by way of Theorem 94$(iii)$ as
\[
(2)=(2, 1+\sqrt{-5})^2.\]
It follows that $I$ must be a power $J^k$ of $J=(2, 1+\sqrt{-5})$. The ramification equation of $(2)$ implies that $J$ has degree 1, hence $N(J)=2$. But then Lemma 95$(i)$ implies that $2=N(I)=N(J)^k=2^k$, hence $k=1$ and so $I=J$.
\begin{quote}
Conclusion: there are at most \emph{two} ideal classes of $R$, namely $[(1)]$ and $[J]$.
\end{quote}

\emph{Claim}: $J$ is not principal.

Hence the ideal-class group of $R$ is $\{[(1)], [J]\}$ and $R$ has class number 2.

\emph{Proof of the claim}. Suppose there exits $\alpha\in R$ such that $J=(\alpha)$. Lemma 95$(ii)$ implies that
\[
|N(\alpha)|=N(J)=2,\]
hence $N(\alpha)=2$ (all nonzero elements of $R$ have positive norm). But there exist $a, b\in \mathbb{Z}$ such that $\alpha=a+b\sqrt{-5}$, hence
\[
a^2+5b^2=N(\alpha)=2,\]
and this is clearly impossible.

\vspace{0.4cm}
\emph{Example} 3

Let $F=\mathbb{Q}(\sqrt {-23}\ ),\ R=\mathcal{A} \cap F=\mathbb{Z}+\displaystyle{\left(\frac{1+\sqrt {-23}}{2}\right)}\ \mathbb{Z}.$ Then $n=2,\ s=1,\ d=-23$ hence $\lambda=\displaystyle{\frac{2\sqrt {23}}{\pi}}<4$, and so every ideal class contains an ideal with norm 1, 2, or 3. As in example 2, every ideal of norm 2 (respectively, 3) must have all of its prime factors containing 2 (respectively, 3), and so factoring via Theorems 93 and 94, we obtain
\[
(2)=\left(2,\ \frac{1+\sqrt {-23}}{2}\right) \left(2,\ \frac{1-\sqrt {-23}}{2}\right)=I_1I_2,\]
\[
(3)=(3,\ 1+\sqrt {-23}) (3,\ 1-\sqrt {-23})=I_3I_4,\]  
hence the ideals of norm 2 are $I_1,\ I_2$ and the ideals of norm 3 are $I_3,\ I_4$. 

The proof of Proposition 23 implies that the elements of $R$ are of the form $a+b\sqrt{-23}$, where $a, b\in \mathbb{Z}$ or $\frac{1}{2}(a+b\sqrt{-23})$, where $a$ and $b$ are \emph{odd} elements of $\mathbb{Z}$. Hence the norm of an element of $R$ is either $a^2+23b^2$ or $\frac{1}{4}(a^2+23b^2)$ for $a, b\in \mathbb{Z}$, neither of which can be 2 or 3. Hence $I_1,\ I_2,\ I_3,$ and $I_4$ are all not principal.
\begin{quote}
Conclusion: in order to calculate the ideal-class group of $R$, we must determine the inequivalent ideals among $I_1,\ I_2,\ I_3,$ and $I_4$.
\end{quote}

We first look at $I_1$ and $I_4$. $I_1\sim I_4$ if and only if $[I_1][I_4]^{-1}=[(1)]$. But $I_3I_4=(3)\sim (1)$, and so $[I_4]^{-1}=[I_3]$, hence we need to see if $I_1I_3$ is principal. 

Lemma 95$(i)$ implies that
\begin{equation*}
N(I_1I_3)=N(I_1)N(I_3)=2\cdot3=6.\tag{1}\]

\emph{Claim}: an ideal $I\not= \{0\}$ of $R$ is principal if and only if there exits $\alpha \in R$ such that $N(\alpha)=N(I)$ and there is a generating set $S$ of $I$ such that $s/\alpha\in R$, for all $s\in S$.

The necessity of this is clear. For the sufficiency, let $\alpha\in R$ satisfy the stated conditions. Then $J=(1/\alpha)I$ is an ideal of $R$ and Lemma 95 $(i), (ii)$ imply that
\[
N(I)=N\big((\alpha)J\big)=N(I)N(J).\]
hence $N(J)=1$ and so $J=(1)$, whence $I=(\alpha)$.

So in light of (1), we must look for elements of $R$ of norm 6. If $a, b\in \mathbb{Z}$ then $a^2+23b^2\not= 6$; on the other hand, 
\[
6=\frac{a^2+23b^2}{4}\]
if and only if $a=\pm 1$ and $b=\pm 1$. Hence there are exactly two principal ideals of norm 6: $\left(\frac{1}{2}(1\pm \sqrt{-23}\ )\right)$. Let $\alpha=\displaystyle{\frac{1+\sqrt {-23}}{2}}$. We have that
\[
I_1I_3=\left(6,\ 2+2\sqrt{-23},\ \frac{3}{2}\big(1+\sqrt {-23}\ \big),\ \frac{1}{2}\big(1+\sqrt {-23}\ \big)^2\right).\]
Now divide each of these generators by $\alpha$: you always get an element of $R$. Hence by the claim, $I_1I_3=(\alpha)$, and so $I_1\sim I_4$.

We assert next that $I_2\sim I_3$. One way to see this is to use the above reasoning again, but there is a more elegant way that we will now explain.

Let $F$ be a Galois field, $A=\mathcal{A}\cap F$. Each of the embeddings of $F$ over $\mathbb{Q}$ is a Galois automorphism of $F$, hence restricts to an automorphism of $A$. Hence 
\begin{quote}
if $\sigma$ is an embedding of $F$ over $\mathbb{Q}$ and $I$ and $J$ are ideals of $A$, then $I\sim J$ if and only if $\sigma(I)\sim \sigma(J)$.
\end{quote}
Moreover, these $\sigma$'s are obtained precisely from the maps of $F$ induced by the conjugates over $\mathbb{Q}$ of a primitive element of $F$ over $\mathbb{Q}$.

If we apply the above observations to the Galois field $F=\mathbb{Q}(\sqrt{-23}\ )$ then the map $q_1+q_2\sqrt{-23}\rightarrow q_1-q_2\sqrt{-23},\ q_1, q_2\in\mathbb{Q}$, of $F$ induced by the conjugate $-\sqrt{-23}$ of $\sqrt{-23}$ restricts to an automorphism $\sigma$ of $R$ which fixes each element of $\mathbb{Z}$. Hence $\sigma(I_1)=I_2$ and $\sigma(I_4)=I_3$ and so we conclude from the equivalence $I_1\sim I_4$ that $I_2\sim I_3$.

We now claim that $I_1$ is \emph{not} equivalent to $I_2$. Otherwise, $[(1)]=[(2)]=[I_1I_2]=[I_1^2]$, hence there exits $\alpha\in R$ such that $I_1^2=(\alpha)$, and so $N(\alpha)=N(I_1^2)=4$, whence $\alpha=\pm 2$. But then $(2)I_1=I_1^2I_2=(2)I_2$, hence $I_1=I_2$, which contradicts the fact that these ideals are distinct.

It follows that the ideal-class group of $R$ is $\{[(1)], [I_1], [I_2]\}$, and $R$ has class number 3. Since the ideal-class group is of prime order, it's cyclic, and since the order is 3, both ideal classes $[I_1]$ and $[I_2]$ are generators of the group.

\emph{Proof of Lemma} 95.

$(i)$ If $I=\{0\}$ or $(1)$ then $N(I)=+\infty$ or 1, and $IJ=\{0\}$ or $J$, hence the conclusion of Lemma 95$(i)$ clearly holds in this case, and so we may assume that $\{0\}\not= I\not= (1),\ \{0\}\not= J\not= (1).$

If $\mathcal{S}$ (respectively, $\mathcal{T}$) is the set of prime ideals containing $I$ (respectively, $J$) then $I$ and $J$ have prime factorizations
\[
I=\prod_{P\in \mathcal{S}} P^{e(P)},\ J=\prod_{Q\in \mathcal{T}} Q^{e(Q)}.\]
If $P\in \mathcal{S} \cap \mathcal{T}$ then we let $e_1(P)$ (respectively, $e_2(P)$) denote the multiplicity of $P$ in $I$ (respectively, $J$). Hence $IJ$ has prime factorization
\[
IJ=\prod_{P\in \mathcal{S}\setminus \mathcal{T}} P^{e(P)} \prod_{P\in \mathcal{S}\cap \mathcal{T}} P^{e_1(P)+e_2(P)} \prod_{Q\in \mathcal{T}\setminus \mathcal{S}} Q^{e(Q)},\]
where an empty product, if it occurs, is replaced by $R$. Because all of the prime factors in this product are distinct, the definition of the norm of an ideal and Lemma 80 imply that
\begin{equation*}
N(IJ)=\prod_{P\in \mathcal{S}\setminus \mathcal{T}} \big|R\big/P^{e(P)}\big| \prod_{P\in \mathcal{S}\cap \mathcal{T}} \big|R\big/ P^{e_1(P)+e_2(P)}\big| \prod_{Q\in \mathcal{T}\setminus \mathcal{S}} \big|R\big/Q^{e(Q)}\big|,\tag{2}\]

If $P\in \mathcal{S} \cap \mathcal{T}$, $f$ is the degree of $P$, and $q$ is the rational prime contained in $P$ then Lemma 81 implies that
\begin{eqnarray*}
\big|R\big/ P^{e_1(P)+e_2(P)}\big|&=&q^{f(e_1(P)+e_2(P))}\\
&=&q^{fe_1(P)}q^{fe_2(P)}\\
&=&\big|R\big/ P^{e_1(P)}\big| \big|R\big/P^{e_2(P)}\big|.
\end{eqnarray*}
Upon substitution of this into (2) and using Lemma 80 again, we find that
\begin{eqnarray*}
N(IJ)&=&\prod_{P\in \mathcal{S}} \big|R\big/P^{e(P)}\big| \prod_{Q\in \mathcal{T}} \big|R\big/ Q^{e(Q)}\big|\\
&=&\Big|R\Big/\Big(\prod_{P\in \mathcal{S}} P^{e(P)}\Big)\Big| \Big|R\Big/\Big(\prod_{Q\in \mathcal{T}} Q^{e(Q)}\Big)\Big|\\
&=&N(I)N(J).
\end{eqnarray*}

$(ii)$ Since this conclusion has been used only for quadratic number fields, we will verify it when $F$ is a Galois field, where the argument simplifies considerably. We will then indicate a way to prove the full-strength version.

So assume $F$ is Galois, with $n=[F:\mathbb{Q}]$. Let $G$ be the Galois group of $F$, and let $0\not= \alpha \in R$. Then $N(\alpha)=\prod_{\sigma \in G} \sigma(\alpha)$, hence
\[
(|N(\alpha)|)= \prod_{\sigma \in G} (\sigma(\alpha)).\]
Lemma 95$(i)$and the proof of Lemma 76 imply that
\begin{equation*}
|N(\alpha)|^n=N\big((|N(\alpha)|)\big)=\prod_{\sigma\in G} N\big((\sigma(\alpha))\big).\tag{3}\]
Since each $\sigma \in G$ is an automorphism of $R$, it follows that
\[
N\big((\sigma(\alpha))\big)=\big|R/(\sigma(\alpha))\big|=|R/(\alpha)|=N\big((\alpha)\big),\ \textrm{for all $\sigma \in G$}.\]
Because $|G|=n$, it hence follows from (3) that
\[
|N(\alpha)|^n=N\big((\alpha)\big)^n,\]
i.e., $|N(\alpha)|=N\big((\alpha)\big)$.

A proof of the general case of $(ii)$ can be easily deduced from the following useful formula for the norm of a nonzero ideal $I$ in an arbitrary number field $F$: if $\{\alpha_1,\dots,\alpha_n\}$ is a basis of $I$ and $d$ is the discriminant of $F$, then
\begin{equation*}
N(I)=\left|\frac{\Delta(\alpha_1,\dots,\alpha_n)}{d}\right|^{1/2}.\tag{4}\]
For an elementary and fairly elegant proof of this formula, see H. Pollard [10], Theorem 9.10.

Now for the proof of $(ii)$ for arbitrary $F$. Let $0\not=\alpha\in R$ and let $\{\omega_1,\dots,\omega_n\}$ be an integral basis of $F$. Then $\{\alpha\omega_1,\dots,\alpha\omega_n\}$ is a basis of $(\alpha)$ and 
\begin{eqnarray*}
\Delta(\alpha\omega_1,\dots,\alpha\omega_n)&=&\big(\det [\alpha^{(j)}\omega_i^{(j)}]\big)^2\\
&=&\Big(\prod_j \alpha^{(j)}\Big)^2\Big(\det [\omega_i^{(j)}]\Big)^2\\
&=&N(\alpha)^2\Delta(\omega_1,\dots,\omega_n)\\
&=&d\cdot N(\alpha)^2.
\end{eqnarray*}
Hence from (4), we conclude that $N(\alpha)^2=N\big((\alpha)\big)^2$, i.e., $|N(\alpha)|=N\big((\alpha)\big)$.

$(iii)$ 
The proof of this is fairly long and intricate, so our discussion will focus on the main ideas involved; for the rest of the details, see Marcus [9], Chapter 5.

The argument makes use of geometric estimates which employ some interesting  volume calculations of certain convex subsets of $E^n$.

We begin with the geometry.

\vspace{0.4cm}
\textbf{Definitions}. An \emph{n-dimensional vector lattice in $E^n$} is a set $L$ of the form
\[
\mathbb{Z} e_1+\cdots+\mathbb{Z} e_n,\]
where $\{e_1,\dots,e_n\}$ is a vector basis of $E^n$. The \emph{co-volume of L} is the $n$-dimensional volume of the set
\[
\Big\{\sum_1^n x_ie_i: x_i\in E,\ 0\leq x_i<1,\ \textrm{for all $i$}.\Big\}.\] 

\vspace{0.4cm}
\emph{Fact} 1. The co-volume of $L$ is the absolute value of the determinant of the $n\times n$ matrix whose rows are $e_1,\dots,e_n$.

\vspace{0.4cm}
This is a well-known fact from linear algebra: see R. C. Buck [3], Theorem 8.3.3.

\vspace{0.4cm}
Any $n$-dimensional vector lattice, when endowed with the vector addition from $E^n$, is an abelian group, hence if $M$ and $L$ are $n$-dimensional vector lattices with $M\subseteq L$ then $M$ is a subgroup of $L$.

\vspace{0.4cm}
\emph{Fact} 2. The index $[L:M]=|L/M|$ of $M$ in $L$ is finite, and
\[
\textrm{the co-volume of $M=[L:M]\times$(the co-volume of $L$).}\]
(Marcus [9], Chapter 5, exercise 3).

\vspace{0.4cm}
Recall that we have set
\[
r=\textrm{the number of real embeddings of $F$ over $\mathbb{Q}$,}\]
\[
s=\frac{1}{2}\times\textrm{the number of complex embeddings of $F$ over $\mathbb{Q}$,}\]
and so 
\[
n=r+2s.\]
We define a ``norm"  $\rho$ on $E^n$ like so: if $x=(x_1,\dots,x_n)\in E^n$ then
\[
\rho(x)=\prod_{i=1}^r x_i\times \prod_{i=1}^s (x_{r+2i-1}^2+x_{r+2i}^2).\]

\vspace{0.4cm}
\emph{Fact} 3. If $L$ is an $n$-dimensional vector lattice then there is a point $0\not= x\in L$ such that
\begin{equation*}
|\rho(x)|\leq \frac{n!}{n^n}\left(\frac{8}{\pi}\right)^s\times(\textrm{co-volume of $L$).}\tag{5}\]

\vspace{0.4cm}

This estimate lies fairly deep and is the main technical step in the proof of Lemma 95$(iii)$. Its derivation depends on the following geometric property of compact, convex, centrally symmetric subsets $C$ of $E^n$ (a subset $S$ of $E^n$ is \emph{centrally symmetric} if it is symmetric about 0, i.e., if $s\in S$ then $-s\in S$): if the ($n$-dimensional) volume of $C$ is positive and if $|\rho(x)|\leq 1$ for all $x\in C$ then every $n$-dimensional lattice $L$ contains a point $x\not=0$ such that
\begin{equation*}
|\rho(x)|\leq  \frac{2^n}{\textrm{volume of $C$}}\times (\textrm{co-volume of $L$)}\tag{6}\] 
(Marcus [9], Corollary of the Lemma, p. 138).

To obtain (5), one considers the subset $C$ of $E^n$ defined by
\[
C=\Big\{(x_1,\dots,x_n)\in E^n: \sum_{i=1}^r |x_i|+2\sum_{i=1}^s \sqrt{x_{r+2i-1}^2+x_{r+2i}^2}\ \leq n\Big\}.\]
It is not difficult to show that $C$ is compact, convex, and centrally symmetric. We check that $|\rho(x)\leq 1$ for all $x\in C$ by using the geometric-mean-arithmetic-mean  inequality, which asserts that the geometric mean $\displaystyle{\Big({\prod_{i=i}^n y_i}\Big)^{1/n}}$ of the sequence $(y_1,\dots,y_n)$ of nonnegative real numbers does not exceed its arithmetic mean $\displaystyle{\frac{1}{n}\sum_{i=1}^n y_i}$. When we take $x=(x_1,\dots,x_n)\in C$ and apply this inequality to the sequence of numbers
\[
|x_1|,\dots,|x_r|, \sqrt{x_{r+1}^2+x_{r+2}^2}, \sqrt{x_{r+1}^2+x_{r+2}^2},\dots,\sqrt{x_{n-1}^2+x_n^2}, \sqrt{x_{n-1}^2+x_n^2},\]
the geometric mean is $\sqrt[n]{|\rho(x)|}$ and its arithmetic mean does not exceed 1. Next one calculates that the volume of $C$ is
\[
\frac{n^n}{n!}\cdot2^r\cdot\Big(\frac{\pi}{2}\Big)^s.\]
We can then insert $C$ into (6) to deduce (5). For the rest of the details here, see Marcus [9], Chapter 5, proof of Theorem 37.

Although the mapping $\rho$ at first glance seems rather contrived, it in fact arises naturally as follows. Let $\sigma_1,\dots,\sigma_r$ denote the real embeddings of $F$ over $\mathbb{Q}$ and $\tau_1,\ \bar{\tau_1},\dots,\tau_s,\ \bar{\tau_s}$ the complex embeddings, with the latter indexed so that if $\alpha\in F$ then
\[
\bar{\tau_i}(\alpha)=\textrm{the complex conjugate of $\tau_i(\alpha)$}.\]
(It follows from Proposition 14 that the complex embeddings of $F$ occur in complex-conjugate pairs.) Then it can be shown that the map $\psi: F\rightarrow E^n$ defined by
\[
\psi(\alpha)=(\sigma_1(\alpha),\dots,\sigma_r(\alpha), \textrm{Re($\tau_1(\alpha)$), Im($\tau_1(\alpha)$),\dots,Re($\tau_s(\alpha)$), Im($\tau_s(\alpha)$))},\ \alpha\in F,\]
where Re and Im denote the real and imaginary parts of a complex number, is an injective homomorphism of the additive group of $F$ into $E^n$. It is then straightforward to check that $\rho$ is simply the map on $E^n$ induced by $\psi$ and the field norm $N: F\rightarrow \mathbb{Q}$, i.e., we have

\vspace{0.4cm}
\emph{Fact} 4.
\[
\rho\big(\psi(\alpha)\big)=N(\alpha),\ \alpha\in F.\]

\vspace{0.4cm}

The mapping $\psi$ will be used to tie the ideals of $R$ to the geometry and group structure of $E^n$ by means of

\vspace{0.4cm}
\emph{Fact} 5. $\psi(R)$ is an $n$-dimensional lattice in $E^n$ and if $\{0\}\not=I$ is an ideal of $R$ then $\psi(I)$ is an $n$-dimensional sublattice of $\psi(R)$ and 

\vspace{0.2cm}
$(a)\  N(I)=[\psi(R):\psi(I)].$

\vspace{0.2cm}
\noindent Moreover, if $d$ is the discriminant of $F$ then

\vspace{0.2cm}

$(b)$ the co-volume of $\psi(R)$ is $2^{-s}\sqrt{|d|}\ .$

\vspace{0.4cm}
In order to verify Fact 5, we first assert that if $\{\alpha_1,\dots,\alpha_n\}$ is a linearly independent subset of $F$ over $\mathbb{Q}$ then $\{\psi(\alpha_1),\dots,\psi(\alpha_n)\}$ is a vector basis of $E^n$. This can be shown by considering the determinant of the $n\times n$ matrix $M$ whose rows are $\psi(\alpha_1),\dots, \psi(\alpha_n)$. Performing elementary column operations on $M$ converts this determinant into $\displaystyle{\frac{1}{(2\sqrt{-1}\ )^s}} \times$ the determinant of the $n\times n$ matrix $M^{\prime}$ whose $i$-th row is
\[
(\sigma_1(\alpha_i),\dots,\sigma_r(\alpha_i), \tau_1(\alpha_i), \bar{\tau}_1(\alpha_i),\dots,\tau_s(\alpha_i), \bar{\tau}_s(\alpha_i)),\ i=1,\dots,n.\]
Hence the definition of discriminant and Corollary 18 imply that
\[
(\det M^{\prime})^2=\Delta(\alpha_1,\dots,\alpha_n)\not= 0,\]
and so
\[
\det M=\frac{1}{(2\sqrt{-1}\ )^s}\det M^{\prime}\not= 0.\]
It follows that $\{\psi(\alpha_1),\dots,\psi(\alpha_n)\}$ is  linearly independent over $E$, and is hence a vector basis of $E^n$. 

Observe next that in light of the definition of $\psi$,
\begin{equation*}
\psi(q\alpha)=q\psi(\alpha),\ \textrm{for all $q\in \mathbb{Q}$ and for all $\alpha \in F$.}\tag{6}\]
Hence if $\{\omega_1,\dots,\omega_n\}$ and  $\{\beta_1,\dots,\beta_n\}$ are, respectively, an integral basis of $F$ and a basis of $I$ then  $\{\psi(\omega_1),\dots,\psi(\omega_n)\}$ and $\{\psi(\beta_1),\dots,\psi(\beta_n)\}$ are vector bases of $E^n$ and, as a consequence of (6),
\[
\psi(R)=\sum_1^n \mathbb{Z} \psi(\omega_i),\ \psi(I)=\sum_1^n \mathbb{Z} \psi(\beta_i),\]
i.e., $\psi(I)$ and $\psi(R)$ are $n$-dimensional vector lattices. When we recall that $\psi$ is a group isomorphism on the additive group of $R$, it hence follows that
\[
[\psi(R):\psi(I)]=\big|\psi(R)/\psi(I)\big|=|R/I|=N(I).\]
This proves $(a)$.

In order to verify $(b)$, observe that if $M$ is the $n \times n$ matrix with rows $\psi(\omega_1),\dots,\psi(\omega_n),$ 
then Fact 1 and a previous calculation imply that 
\begin{eqnarray*}
\textrm{the co-volume of $\psi(R)$}&=&|\det M|\\
&=&\left|\frac{1}{(2\sqrt{-1}\ )^s}\right|\sqrt{|\Delta(\omega_1,\dots,\omega_n)|}\\
&=&2^{-s}\sqrt{|d|}.
\end{eqnarray*}

We now have all the ingredients required for a proof of Lemma 95$(iii)$. From Facts 2 and 5, we deduce that for any nonzero ideal $I$ of $R$,
\begin{equation*}
\textrm{co-volume of $\psi(I)$}=N(I)\times \textrm{co-volume of $\psi(R)$}=2^{-s}\sqrt{|d|}\cdot N(I).\tag{7}\]
Apply Fact 3 to the lattice $L=\psi(I)$ and then use Fact 4 and (7) to find $0\not=\alpha\in I$ such that
\begin{eqnarray*}
(8)\hspace{2cm}  |N(\alpha)|&=&\big|\rho\big(\psi(\alpha)\big)\big|\\
&\leq&\frac{n!}{n^n}\left(\frac{8}{\pi}\right)^s\times \textrm{co-volume of $\psi(I)$}\\
&=&\frac{n!}{n^n}\left(\frac{4}{\pi}\right)^s \sqrt{|d|}\cdot N(I).
\end{eqnarray*}

Now let $X$ be an ideal class of $R$. Let $I\in X^{-1}$, and choose $0\not=\alpha\in I$ which satisfies estimate (8). Because $(\alpha)\subseteq I$, there exits an ideal $J\not=\{0\}$ such that $(\alpha)=IJ$ (Proposition 57), hence
\[
[(1)]=[I][J]=X^{-1}[J],\]
and so $[J]=X$, i.e., $J\in X$. By virtue of Lemma 95$(i)$, $(ii)$ and (8), it follows that
\[
N(I)N(J)=|N(\alpha)|\leq \frac{n!}{n^n}\left(\frac{4}{\pi}\right)^s \sqrt{|d|}\cdot N(I),\]
and because $N(I)>0$, we obtain
\[
N(J)\leq \frac{n!}{n^n}\left(\frac{4}{\pi}\right)^s \sqrt{|d|}.\]
$\hspace{15.5cm}\textrm{QED}$
\chapter{Structure of the Group of Units in a Number Ring}

Let $F$ be a number field, $R=\mathcal{A}\cap F$. If $I$ and $J$ are nonzero ideals of $R$ then $I\sim J$ if and only if there exits $0\not=\alpha\in F$ such that $J=\alpha I$. Hence if we set
\[
M(I)=\{0\not=\alpha\in F: \alpha I\subseteq R\},\]
then the ideal class containing $I$ is given by
\[
[I]=\{\alpha I: \alpha\in M(I)\}.\]

Now consider the set $F^{*}$ of all nonzero elements of $F$. When $F^{*}$ is endowed with the multiplication of $F$, it becomes an abelian group, the \emph{multiplicative group of F}. Let $U(R)$ denote the group of units of $R$. $U(R)$ is a subgroup of $F^{*}$; let $C$ denote a set of coset representatives of the quotient group $F^{*}/U(R)$, i.e., each coset of $F^{*}/U(R)$ contains exactly one element from $C$. It is easy to see that
\[
[I]=\{\alpha I: \alpha\in C\cap M(I)\},\]
so that in order to calculate the ideal class containing $I$, you need only look for elements \emph{in C} that multiply $I$ into $R$. Hence, at least in principle, the structure of $U(R)$ is of considerable interest in computing ideal-class groups and/or the class number of $R$: the more that you know about the structure of $U(R)$, the more that you know about $C$, and this information may (or may not!) be useful in finding ideal classes, or at least counting them.

Another important situation in which $U(R)$ plays a crucial role is concerned with the algebra of polynomials over $R$. If $p(x)$ is a nonzero polynomial in $R[x]$ then the \emph{division algorithm holds for $p(x)$} if for all nonzero polynomials $f(x)\in R[x]$, there exist unique polynomials $q(x), r(x) \in R[x]$ such that $f(x)=q(x)p(x)+r(x)$ and either $r(x)\equiv 0$ or the degree of $r(x)$ is less that the degree of $p(x)$. The division algorithm in $R[x]$ is a very useful tool for studying the algebraic structure of $R$, and it holds for $p(x)$ if and only if the leading coefficient of $p(x)$ is a unit. Hence in order to use the division algorithm in $R[x]$, one needs to keep track of the units in $R$.

Finally, in Chapter 18, we will indicate how the structure of $U(R)$ determines some important and interesting features of the distribution of ideals in $R$.
\newpage
\begin{center}
\emph{Examples}
\end{center}

(1) \emph{Quadratic number fields}

Let $m$ be a square-free rational integer, $F=\mathbb{Q}(\sqrt m\ ),\ R=\mathcal{A}\cap F$. Let $N: F\rightarrow \mathbb{Q}$ be the norm map. Proposition 38 implies that
\[
\alpha \in R\ \textrm{is a unit if and only if $N(\alpha)=\pm 1$.}\]

In order to determine the units of $R$, we must distinguish, according to Proposition 23, two cases:

\emph{Case} 1. $m\not \equiv 1$ mod 4.

Here $R=\mathbb{Z}+\sqrt m\ \mathbb{Z}$, hence $\alpha=x+y\sqrt m\in R$ is a unit if and only if 
\begin{equation*}
x^2-my^2=N(\alpha)=\pm 1.\tag{1}\]

\emph{Case} 2: $m \equiv 1$ mod 4.  

Here, in addition to $\mathbb{Z}+\sqrt m\ \mathbb{Z}$, $R$ also has elements of the form
\[
\frac{x+y\sqrt m}{2},\ \textrm{$x$ and $y$ both \emph{odd} rational integers,}\]
and so in addition to the units coming from solutions $(x, y)\in \mathbb{Z}\times \mathbb{Z}$ of (1), there are also units of the form $\displaystyle{\frac{x+y\sqrt m}{2}}$ where
\begin{equation*}
x^2-my^2=\pm 4,\ \textrm{$x$ and $y$ both odd rational integers.}\tag{2}\]

Suppose now that $m<0$, i.e., $F$ is an \emph{imaginary quadratic number field}. Then only the $+$ sign can occur in (1) and (2). If $m<-1$ then  $x^2-my^2=1$ has only the solutions $x=\pm 1,\ y=0$ (otherwise $x^2-my^2>1$), and if $m<-4$ then  $x^2-my^2=4$ has no solutions with $x$ and $y$ both odd.
\[
\textrm{Conclusion: if $m<-4$ then $U(R)=\{-1, 1\}$.}\]

Next, consider $m=-1, -2, -3$.

\vspace{0.2cm}
$m=-1$.

\vspace{0.2cm}
We are in Case 1, hence $R=\mathbb{Z}+\sqrt {-1}\ \mathbb{Z}=$ the Gaussian integers, and so units $x+y\sqrt{-1}$ are determined by the solutions of $x^2+y^2=1$, i.e., $(x, y)=(1, 0), (0, 1), (-1, 0), (0, -1).$ Hence 
\[
U(R)=\{ \pm 1, \pm\sqrt{-1}\}.\]

\vspace{0.2cm}
$m=-2$.

\vspace{0.2cm}
This is also in Case 1, so units $x+y\sqrt{-2}$ are determined by the solutions of $x^2+2y^2=1$, i.e., $(x, y)=(1, 0), (-1, 0),$ and so
\[
U(R)=\{-1, 1\}.\]

\vspace{0.2cm}
$m=-3$.

\vspace{0.2cm}
This is in Case 2, so units other than $\pm 1$ come from solutions of $x^2+3y^2=4$, with $x$ and $y$ both odd. Hence $(x, y)=(1, 1), (1, -1), (-1, 1), (-1, -1)$, 
and so
\[
U(R)=\left\{\pm 1, \frac{1\pm \sqrt{-3}}{2}, \frac{-1\pm \sqrt{-3}}{2}\right\}.\]

We have proved

\vspace{0.4cm}
\textbf{Proposition 96}. \textit{ The imaginary quadratic number field $\mathbb{Q}(\sqrt m )$ has only the units $\pm 1$, unless $m=-1$, in which case there are the additional units $\pm \sqrt{-1}$, or unless $m=-3$, in which case there are the additional units $\displaystyle{\frac{1\pm \sqrt{-3}}{2}}$ and $\displaystyle{\frac{-1\pm \sqrt{-3}}{2}}$.}

\vspace{0.4cm}
Suppose next that $m>0$, i.e., $F$ is a \emph{real quadratic number field}. Here the group of units has a much more complicated structure. The units are now determined by solutions $(x, y)\in \mathbb{Z}\times \mathbb{Z}$ of the equations
\[
x^2-my^2=\pm 1,\ x^2-my^2=\pm 4\]
for a fixed positive square-free $m\in \mathbb{Z}$, called \emph{Pell's equations}. These equations always have infinitely many solutions, and very efficient algorithms exist for finding them; for further details on this topic , see H. Cohen [5], section 6.3.5 and H. Cohen [4], sections 5.7 and 5.8. We will be content to discuss the following classical example, which illustrates the typical structure of the group of units in real quadratic fields.

\vspace{0.4cm}
\textbf{Proposition 97}. \textit{If $F=\mathbb{Q}(\sqrt 2 )$ then $\lambda=1+\sqrt 2$ is a unit in $R=\mathcal{A} \cap F$ and $U(R)=\{\pm \lambda^n: n\in \mathbb{Z}\}$.}

\vspace{0.4cm}
\emph{Proof}. For this we will need

\vspace{0.4cm}
\textbf{Lemma 98}. \textit{ $R$ has no unit between $1$ and $\lambda$.}

\vspace{0.4cm}
Assume this lemma for now. $1+\sqrt 2$ is a unit because it has norm $-1$, and so $\pm\lambda^n$ is a unit for each $n\in \mathbb{Z}$.

Let $\varepsilon$ be a unit. Because all elements of $\mathbb{Q}(\sqrt 2 )$ are real, $\varepsilon>0$ or $\varepsilon<0$. Suppose that $\varepsilon>0$. As $\lambda>1$, it follows that there is an $n\in \mathbb{Z}$ such that $\lambda^n\leq \varepsilon<\lambda^{n+1}$. If $\lambda^n< \varepsilon<\lambda^{n+1}$ then $1<\varepsilon\lambda^{-n}<\lambda$.
Since $\varepsilon\lambda^{-n}$ is a unit, this contradicts Lemma 98. Hence $\varepsilon=\lambda^{n}$. If $\varepsilon<0$ then apply this argument to $-\varepsilon$, also a unit, to get $\varepsilon=-\lambda^n$ for some  $n\in \mathbb{Z}$.$\hspace{15.2cm}\textrm{QED}$

\vspace{0.2cm}
\emph{Proof of Lemma} 98. Suppose that $\varepsilon=x+y\sqrt 2$ is a unit such that
\begin{equation*}
1<\varepsilon<1+\sqrt 2.\tag{3}\]
Since $\pm 1=N(\varepsilon)=x^2-2y^2$, we have that
\[
x-y\sqrt 2=\pm\frac{1}{x+y\sqrt 2},\]
and so
\begin{equation*}
-1<x-y\sqrt 2<1.\tag{4}\]
Add (3)and (4) to obtain $0<2x<2+\sqrt 2$, hence $0<x<1.8$. But $x\in \mathbb{Z}$, and so $x=1$. Then $1<1+y\sqrt 2<1+\sqrt 2$, hence $0<y<1$, which is impossible because $y$ is a rational integer.$\hspace{14.2cm}\textrm{QED}$

\vspace{0.2cm}
The structure of the group of units for a real quadratic field as illustrated by Proposition 97 persists in general. We will eventually see that in every real quadratic field there is a unit $\lambda>1$, called the \emph{fundamental unit}, such that the group of units is $\{\pm \lambda^n: n\in \mathbb{Z}\}$. Hence the group of units is infinite and there are infinitely many positive (respectively, negative) units of arbitrarily large (respectively, arbitrarily small) absolute value. Moreover, if $\mathbb{Z}/2\mathbb{Z}=\{\bar{0}, \bar{1}\}$ is the cyclic group of order 2, then the map $\lambda^n\rightarrow (\bar{0}, n),\ -\lambda^n\rightarrow (\bar{1}, n),\ n\in \mathbb{Z}$, is a group isomorphism of the group of units onto $\mathbb{Z}/2\mathbb{Z}\oplus \mathbb{Z}$.

\vspace{0.4cm}
(2) \emph{Cyclotomic number fields}

\vspace{0.4cm}
Let $p$ be an odd rational prime, $\omega=e^{2\pi i/p}$. The structure of the group of units of the cyclotomic field $\mathbb{Q}(\omega)$ is more complicated than the structure of the group of units in a quadratic field, hence we will simply state the results and forgo the proofs.

In order to do that, we need to recall some basic group theoretic facts and notation.The rational integers $\mathbb{Z}$ equipped with its addition is a countably infinite cyclic group and is the unique such group up to isomorphism. If $m\in \mathbb{Z},\ m\geq 1$, then the set $\mathbb{Z}_m=\mathbb{Z}/m\mathbb{Z}$ of integers modulo $m$ equipped with its addition is a cyclic group of order $m$ and is the unique such group up to isomorphism. We will denote by $\mathbb{Z}^m$ the abelian group defined by the direct sum of $m$ copies of $\mathbb{Z}$, i.e., $\mathbb{Z}^m=\{(a_1,\dots,a_m): a_i\in \mathbb{Z},\ \textrm{for all $i$}\}$, with the group operation defined coordinate-wise by $(a_i)+(b_i)=(a_i+b_i)$. $\mathbb{Z}^m$ is free abelian of rank $m$, and is, up to isomorphism, the unique free abelian group of rank $m$.

\vspace{0.4cm}
\textbf{Proposition 99}. \textit{Let $\varepsilon=e^{\pi i/p}=\omega^{(p+1)/2}$ and let $r=\displaystyle{\frac{p-3}{2}}$.}

\textit{$(i)$ If $p=3$ then the group of units $U(\mathbb{Q}(\omega))$ is $\{\varepsilon^k: k=0, 1, 2, 3, 4, 5\}$, hence the map $\varepsilon^k\rightarrow k$ is an isomorphism of $U(\mathbb{Q}(\omega))$ onto $\mathbb{Z}_6$.}

\textit{$(ii)$ If $p\geq 5$ then there exits real, positive units $u_1,\dots, u_r$ in $\mathbb{Q}(\omega)$ such that every unit in $\mathbb{Q}(\omega)$ can be written uniquely as}
\[
u=\varepsilon^k\prod_{i=1}^r u_i^{a_i},\]
\textit{where}
\[
k\in\{0, 1, \dots,2p-1\}\ \textrm{and}\ (a_1,\dots,a_r)\in\mathbb{Z}^r.\]
\textit{Hence the map $u\rightarrow (k, a_1,\dots,a_r)$ is an isomorphism of $U(\mathbb{Q}(\omega))$ onto $\mathbb{Z}_{2p}\oplus \mathbb{Z}^r$.}

\vspace{0.4cm}
\emph{Proof}. See Z. Borevich and I. Shafarevich [2].$\hspace{7.2cm}\textrm{QED}$

\vspace{0.2cm}
N.B. If $p=3$ then $\mathbb{Q}(\omega)=\mathbb{Q}(\sqrt{-3})$, and the 6 sixth roots of unity are 
\[
\pm 1,\ \displaystyle{\frac{1\pm \sqrt{-3}}{2}},\  \displaystyle{\frac{-1\pm \sqrt{-3}}{2}},\] 
and so in this case Proposition 99 and Proposition 96 assert the same thing.

It can be shown that $\{\varepsilon^k: k=0, 1,\dots,2p-1\}$ is the set of roots of unity that are contained in $\mathbb{Q}(\omega)$ (Borevich and Shafarevich [2], Lemma 3.1.3). Hence Proposition 99 implies that $U(\mathbb{Q}(\omega))$ is the direct sum of the finite cyclic group consisting of the roots of unity contained in  $\mathbb{Q}(\omega)$ and a free abelian group of finite rank generated by \emph{fundamental units} $u_1,\dots, u_r$. As we mentioned before, the group of units of a quadratic field also has the exact same structure. We will now show that this is no accident, i.e., we will prove that the the group of units in \emph{any} number ring has the same structure.

\vspace{0.4cm}
We return to the general set-up of a number field $F$ and its number ring $R=\mathcal{A}\cap F$, with group of units $U(R)$ of $R$. Let $n=[F:\mathbb{Q}]$.

\vspace{0.4cm}
\textbf{Proposition 100}. \textit{The set $U_0(R)$ of all roots of unity in $R$ is a finite cyclic subgroup of $U(R)$.}

\vspace{0.4cm}
\emph{Proof}. $U_0(R)$ is clearly the set of all elements of $R$ of multiplicative finite order and is hence closed under multiplication. If $\alpha\in U_0(R)$ then $\alpha^k=1$ for some positive $k\in \mathbb{Z}$, hence $1/\alpha=\alpha^{k-1}\in R$, and so $U_0(R)\subseteq U(R)$ and $U_0(R)$ is closed under reciprocation, whence $U_0(R)$ is a subgroup of $U(R)$. To prove that $U_0(R)$ is cyclic, we need only prove that it is finite, since any finite subgroup of the multiplicative group of a field is cyclic (proof: exercise). Finiteness of $U_0(R)$ is a consequence of 

\vspace{0.4cm}
\textbf{Lemma 101}. \textit{There are at most finitely many elements of $R$ which, together with their conjugates, do not exceed a fixed constant in absolute value.}

\vspace{0.4cm}Since a root of unity and all of its conjugates have absolute value 1, Lemma 101 implies that there are only finitely many roots of unity in $R$.$\hspace{6cm} \textrm{QED}$

\vspace{0.2cm} \emph{Proof of Lemma } 101. Let $C\geq 1$ be a fixed constant. Consider the set $\mathcal{S}$ of all polynomials in $\mathbb{Z}[x]$ of degree $\leq n$, and all of whose coefficients have absolute value $\leq (2C)^n$. $\mathcal{S}$ is a finite set, hence the set 
\[
\mathcal{R}=\bigcup_{p\in \mathcal{S}}\ \textrm{set of all roots of $p$}\]
is a finite set.

Now let $\alpha\in R$, let $m$ be the degree of $\alpha$ over $\mathbb{Q}$, and suppose that
\begin{equation*}
|\alpha^{(i)}|\leq C,\ i=1,\dots,m.\tag{$*$}\]
If $p$ is the minimal polynomial of $\alpha$ over $\mathbb{Q}$ then $p\in \mathbb{Z}[x]$ and the degree $m$ of $p$ is $\leq n$. The proof of Proposition 6 implies that each coefficient of $p$ is $\pm$ an elementary symmetric function in $m$ variables of the $\alpha^{(i)}$'s hence $(*)$ implies that each of these coefficients has absolute value not exceeding $(2C)^m\leq (2C)^n$. Hence $p\in \mathcal{S}$ and so $\alpha \in \mathcal{R}$. Because $\mathcal{R}$ is finite, there are only finitely many $\alpha \in R$ what satisfy $(*)$.$\hspace{9.1cm}\textrm{QED}$

\vspace{0.2cm}
N.B. It is always true that $\{-1, 1\}\subseteq U_0(R)$ and it may be the case (Proposition 96) that $U_0(R)=\{-1, 1\}=U(R)$.

The next theorem, one of the fundamental results of algebraic number theory, describes precisely the structure of $U(R)$. In order to state it, let
\[
r_1=\textrm{the number of real embeddings of $F$,}\]
\[
c=\textrm{the number of complex embeddings of $F$ (which we showed on p. 90 is even),}\]
\[
r_2=\frac{c}{2},\]
\[
r=r_1+r_2-1.\]

\vspace{0.2cm}
\noindent We have $n=r_1+2r_2$, and either $r_1$ or $r_2$ (but not both) may be 0.

Suppose that $r=0$, i.e., $r_1+r_2=1$, hence $(r_1, r_2)=(0, 1)$ or $(1, 0)$. If $r_1=0$ and $r_2=1$ then $n=2$, hence $F$ is an imaginary quadratic field, and so $U(R)$ is determined by Proposition 96. If $r_1=1$ and $r_2=0$ then $n=1$, and so $F=\mathbb{Q}$ and $R=\mathbb{Z}$, hence $U(R)=\{-1, 1\}$. In particular, if $r=0$ then $U(R)=U_0(R)$.  Hence we may assume that $r\geq 1$.

\vspace{0.4cm}
\textbf{Theorem 102}. \textit{$($Dedekind-Dirichlet Units Theorem$)$ If $r\geq 1$, m is the order of $U_0(R)$, and $\rho_0$ is a generator of $U_0(R)$, i.e.,}
\[
U_0(R)=\{\rho_0^k: k=0,1,\dots,m-1\}\]
\textit{then there exist r units $\varepsilon_1,\dots,\varepsilon_r$ in $R$ such that every unit u in R can be written uniquely in the form}
\[
u=\rho_0^k\prod_{i=1}^r \varepsilon_i^{a_i},\]
\textit{for some}
\[
k\in\{0,1,\dots,m-1,\}\ \textrm{and}\ (a_1,\dots,a_r)\in \mathbb{Z}^r.\]
\textit{Hence the map $u\rightarrow (k, a_1,\dots,a_r)$ is an isomorphism of $U(R)$ onto $\mathbb{Z}_m\oplus \mathbb{Z}^r$.}

\vspace{0.4cm}
\textbf{Definitions}. Any set $\varepsilon_1,\dots,\varepsilon_r$ of units of $R$ which satisfy the conclusion of Theorem 102 are called \emph{fundamental units of R} and $r=r_1+r_2-1$ is the \emph{unital rank of R}.

\vspace{0.4cm}
N.B. If $\mathbb{Q}(\sqrt m\ )$ is a real quadratic field then $r=1$ and $U_0(\mathbb{Q}(\sqrt m\ ))=\{-1, 1\}\ (\pm 1$ are the only roots of unity that are real!). Hence Theorem 102 implies that $\mathcal{A}\cap \mathbb{Q}(\sqrt m \ )$ has a fundamental unit $\lambda$, uniquely determined by also requiring that $\lambda>1$, such that $U(\mathbb{Q}(\sqrt m \ ))=\{\pm\lambda^n: n\in \mathbb{Z}\}$. This justifies our previous assertion that the structure of the group of units in the example in Proposition 97 typifies the general situation for real quadratic fields.

N.B. The Fundamental Theorem of Ideal Theory (Theorem 50), the theorem on the finiteness of the class number (Theorem 69), and the Dedekind-Dirichlet Units Theorem (Theorem 102) are the three pillars on which stands the entire theory of algebraic numbers.

\vspace{0.4cm}
\emph{Proof of Theorem} 102. This proof is long and technically intricate, so we will skip some of the more complicated details in the reasoning. We will follow the classical argument given in E. Hecke [7], Chapter V, section 34. For a more modern approach, see Marcus [9], Chapter 5. First, recall the following definitions from abelian group theory:

\vspace{0.4cm}
\textbf{Definitions}. Let $A$ be an abelian group, written multiplicatively. The \emph{torsion subgroup $A_t$ of A} is the subgroup of $A$ consisting of all elements of $A$ of finite order. An element of $A\setminus A_t$ is said to be \emph{torsion-free}.

A (finite) set $\{a_1,\dots,a_k\}$ of elements of $A$ is \emph{independent} if $a_1^{n_1}\cdots a_k^{n_k}=1$ for some $(n_1,\dots,n_k)\in \mathbb{Z}^k$ implies that $n_1=\cdots=n_k=0$.

\vspace{0.4cm}
N.B. Every element of an independent subset of $A$ must be torsion-free.

Let $\theta$ be a primitive element of $F$ over $\mathbb{Q}$. The proof of Theorem 102 makes use of a certain indexing of the conjugates $\theta^{(1)},\dots,\theta^{(n)}$ of $\theta$, described like so: let $\theta^{(1)},\dots,\theta^{(r_1)}$ be any fixed indexing of the real conjugates of $\theta$, and then index the $2r_2$ non-real conjugates so that
\[
\theta^{(p+r_2)}=\textrm{the complex conjugate of $\theta^{(p)},\ p=r_1+1,\dots,r_1+r_2$.}\]
Note that if $\alpha\in F$ then this indexing also induces the same indexing of the conjugates $\alpha^{(1)},\dots,\alpha^{(n)}$ of $\alpha$ over $F$.

We will divide the proof of Theorem 102 into a series of four assertions and their verifications.

\vspace{0.4cm}
\emph{Assertion} 1: $U(R)$ contains at most $r$ independent units.

In order to see this, let $\{u_1,\dots,u_k\}$ be a fixed subset of $U(R)$. We first verify

\emph{Claim} 1. $\{u_1,\dots,u_k\}$ is independent if and only if the only element $(a_1\dots,a_k)\in \mathbb{Z}^k$ for which
\begin{equation*}
\sum_{i=1}^k a_i\log |u_i^{(l)}|=0,\ l=1,\dots,n,\tag{5}\]
is $a_1=\cdots =a_k=0$.

$(\Rightarrow)$ Suppose that for some $(a_1\dots,a_k)\in \mathbb{Z}^k$, equations (5) hold. Then $|(u_1^{a_1}\cdots u_k^{a_k})^{(l)}|=1$, for all $ l=1,\dots,n$, i.e., $u_1^{a_1}\cdots u_k^{a_k}$ is an element of $R$ which, together with all of its conjugates, has absolute value 1. Hence Lemma 101 implies that the positive rational integral powers of $u_1^{a_1}\cdots u_k^{a_k}$ cannot all be distinct, and so there is a positive $m\in \mathbb{Z}$ such that  $u_1^{ma_1}\cdots u_k^{ma_k}=1$. We are assuming that 
$\{u_1,\dots,u_k\}$ is independent, and so it follows that $a_1=\cdots =a_k=0$.

$(\Leftarrow)$ Suppose that for some $(a_1\dots,a_k)\in \mathbb{Z}^k$, $u_1^{a_1}\cdots u_k^{a_k}=1$. Then $(u_1^{a_1}\cdots u_k^{a_k})^{(l)}=1$, for all $l=1,\dots,n$, and so upon taking the absolute value and then the logarithm of these equations, we deduce that (5) holds. Hence from the hypothesis of this implication  it follows that $a_1=\cdots =a_k=0$.

We next verify

\emph{Claim} 2. If
\begin{equation*}
\sum_{i=1}^k a_i\log |u_i^{(l)}|=0,\ l=1,\dots,r,\tag{6}\]
then the equations in (5) for $l=r+1,\dots,n$ automatically hold.

In order to see this, note first that by our indexing of the conjugates of the elements of $F$ it follows that
\begin{equation*}
|u_i^{(p+r_2)}|=|u_i^{(p)}|, \ p=r_1+1,\dots,r_1+r_2.\tag{7}\]
Because $u_i$ is a unit, 
\[
1=|N(u_i)|=\prod_{l=1}^n |u_i^{(l)}|=\prod_{l=1}^{r_1}  |u_i^{(l)}| \prod_{l=r_1+1}^{r_1+r_2}  |u_i^{(l)}|^2,\]
hence
\[
0=\sum_{l=1}^{r_1} \log  |u_i^{(l)}|+2\sum_{l=r_1+1}^{r_1+r_2}  \log  |u_i^{(l)}|,\]
i.e.,
\[
2\log |u_i^{(r_1+r_2)}|=-\sum_{l=1}^{r_1} \log  |u_i^{(l)}|-2\sum_{l=r_1+1}^r  \log  |u_i^{(l)}|.\]
Now multiply this equation by $a_i$, sum from $i=1,\dots,k$, and invoke (6): we obtain
\[
2\sum_{i=1}^k a_i\log |u_i^{(r_1+r_2)}|=-\sum_{l=1}^{r_1}\left(\sum_{i=1}^k a_i\log  |u_i^{(l)}|\right)-2\sum_{l=r_1+1}^r \left(\sum_{i=1}^k a_i\log  |u_i^{(l)}|\right)=0.\]
Hence (5) is true for $l=r_1+r_2$, and so from (7) we deduce (5) for the remaining values of $l$.

N.B. The verification of Claim 2 is the reason why the special indexing of the conjugates over $F$ was introduced.

We deduce from Claims 1 and 2 that

\vspace{0.4cm}
\begin{quote} $\{u_1,\dots,u_k\}$ is independent if and only if the only solution of the equations (6) for $(a_1\dots,a_k)\in \mathbb{Z}^k$ is $a_1=\cdots =a_k=0$.
\end{quote}

\vspace{0.4cm}
One now proves that

\vspace{0.4cm}
\begin{quote} if (6) has a nonzero solution $(a_1\dots,a_k)\in E^k$ then  (6) has a nonzero solution $(a_1\dots,a_k)\in \mathbb{Z}^k$.
\end{quote}

\vspace{0.4cm}
\noindent (Hecke [7], section 34, Lemma (b)). Hence it follows that  $\{u_1,\dots,u_k\}$ is independent if and only if the vectors $(\log |u_i^{(1)}|,\dots,\log |u_i^{(r)}|),\ i=1,\dots,k$ are linearly independent in $E^r$. Consequently, $U(R)$ contains at most $r$ independent units, and Assertion 1 is verified.

\vspace{0.4cm}
\textbf{Definition}. The set of vectors
\[
\big\{(\log |u^{(1)}|,\dots,\log |u^{(n)}|): u \in U(R)\big\}\]
in $E^n$ is called the \emph{logarithm space of R}.

 Our proof of Assertion 1 shows that there are in fact only at most $r=r_1+r_2-1$ free parameters in the logarithm space. We will eventually prove that  there are exactly $r$ free parameters in the logarithm space.

Let $\{u_1,\dots,u_k\}$ be a set of independent units of $R$ of maximum cardinality. Assertion 1 implies that $k\leq r$.

\vspace{0.4cm}
\emph{Assertion} 2. $U(R)$ is finitely generated.

In order to prove this, use the maximality of $k$ and the argument in Hecke [7], section 34, Lemma (c) to verify that there is a positive $M\in \mathbb{Z}$ with the following property: if $u\in U(R)$ then there exits $(g_1,\dots,g_k)\in \mathbb{Z}^k$ such that
\begin{equation*}
\log |u^{(i)}|=\sum_{j=1}^k \frac{g_j}{M} \log |u_j^{(i)}|,\ i=1,\dots,n.\tag{8}\]

Now let $u\in U(R)$. Then (8) implies that 
\[
\big|\big(u_1^{g_1}\cdots u_k^{g_k}u^{-M}\big)^{(i)}\big|=1,\ \textrm{for all $i=1,\dots,n$.}\]
This says that 
\[
\alpha=u_1^{g_1}\cdots u_k^{g_k}u^{-M} \in R,\]
together with all of its conjugates, has absolute value 1, and so by an argument in the proof of Assertion 1, Lemma 101 implies that $\alpha \in U_0(R)$, hence $\alpha=\rho_0^t$ for some $t\in \mathbb{Z}$. Hence
\[
u=u_1^{g_1/M}\dots u_k^{g_k/M}\rho_0^{-t/M}.\]
If we now let $H$ denote the subgroup of the multiplicative group of \textbf{C} generated by the $M$-th roots of $u_1,\dots,u_k$ and $\rho_0$, then $H$ is a finitely generated abelian group and we have just shown that $U(R)$ is a subgroup of $H$. Because every subgroup of a finitely generated abelian group is finitely generated (Hungerford [8], Corollary II.1.7), it follows that $U(R)$ is finitely generated.

\vspace{0.4cm}
\emph{Assertion} 3. $U(R)$ has a set of $k$ fundamental units, i.e., there is a set $\{\varepsilon_1,\dots,\varepsilon_k\}$ of independent units of $R$ which satisfy the conclusion of Theorem 102 with $r$ replaced by $k$.

To see that this is true, we use Assertion 2 and some structure theory for finitely generated abelian groups, to wit, if $A$ is such a group and $A_t$ is its torsion subgroup then there exits a unique free abelian subgroup $B$ of $A$ of finite rank such that $A$ is the internal direct product $A_t\times B$ of $A_t$ and $B$ (Hungerford [8], section II.2). We note next that $U_0(R)$ is the torsion subgroup of the finitely generated abelian group $U(R)$, and so we find a free abelian, finite-rank subgroup $B$ of $U(R)$ such that
\[
U(R)=U_0(R)\times B.\]
Now a basis of $B$ is a finite set of independent elements of $U(R)$, and so by the maximality of $k$,
\[
\textrm{ rank of $B$ $\leq k$.}\]
On the other hand, if $m$ is the order of $U_0(R)$, then $u_1^m,\dots,u_k^m$ are $k$ independent elements of $B$, and hence generate a free subgroup of $B$ of rank $k$. Since a free subgroup of a free abelian group has rank that does not exceed the rank of the group (Hungerford [8], remark before the proof of Theorem II.1.6), we obtain
\[
\textrm{$k \leq$ rank of $B$.}\]
Hence the rank of $B$ is $k$, and so if $\{\varepsilon_1,\dots,\varepsilon_k\}$ is a basis of $B$ then every unit $u$ of $R$ is uniquely represented in the form
\[
u=\rho_0^l\prod_{i=1}^k \varepsilon_i^{a_i},\]
\textrm{for some}
\[
l\in\{0,1,\dots,m-1,\}\ \textrm{and}\ (a_1,\dots,a_k)\in \mathbb{Z}^k,\]
i.e., $\varepsilon_1,\dots,\varepsilon_k$ are fundamental units of $R$.

It follows from the proof of Assertion 3 that Theorem 102 will be a consequence of 

\vspace{0.4cm}
\emph{Assertion} 4. $k=r$.

The verification of this assertion, the heart of the proof of Theorem 102, requires that we produce $r$ independent units in $R$. In order to do that, we need

\vspace{0.4cm}
\textbf{Lemma 103}. \textit{For each $r$-tuple of real numbers $(c_1,\dots,c_r)\not=0$, there exits $u\in U(R)$ such that}
\[
\sum_{i=1}^r c_i\log |u^{(i)}|\not=0.\]

\vspace{0.4cm}
Assume Lemma 103 for now; we use it to inductively construct $r$ units $u_1,\dots,u_r$ of $R$ like so: Lemma 103 implies that there is a unit $u_1$ such that
\[
\log |u_1^{(1)}|\not=0.\]

Suppose that the units $u_1,\dots,u_t$, with $t<r$, have been constructed so that   

\begin{equation*}
\det \left[ \begin{array}{ccc}
\log |u_1^{(1)}|  & \dots &  \log |u_t^{(1)}|\\
\vdots & \dots & \vdots\\
\log |u_1^{(t)}| & \dots & \log |u_t^{(t)}|   \\
 \end{array} \right]\not=0.\tag{9}\]
Consider the matrix

\[
 \left[ \begin{array}{ccc}
\log |u_1^{(1)}|  & \dots &  \log |u_t^{(1)}|\\
\vdots & \dots & \vdots\\
\log |u_1^{(t)}| & \dots & \log |u_t^{(t)}|   \\
\log |u_1^{(t+1)}| & \dots & \log |u_t^{(t+1)}|   \\
 \end{array} \right].\]

\noindent For $i=1,\dots,t+1$, let
\[
M_i=\textrm{$t\times t$ matrix formed from this matrix by deletion of the $i$-th row,}\]

\noindent and set
\[
c_i=(-1)^{t+1+i}\det M_i.\]

\noindent For $i=t+2,\dots,r$ (if any such $i$ exist), set $c_i=0$. Because $c_{t+1}\not=0$ (by (9)), it follows that $(c_1,\dots c_r)\not=0$, hence Lemma 103 implies that there is a unit $u_{t+1}$ such that
\[
\sum_{i=1}^{t+1} c_i\log |u_{t+1}^{(i)}|\not=0.\]
By construction of the $c_i$'s, this sum is the cofactor expansion down the last column of

\[
\det \left[ \begin{array}{ccc}
\log |u_1^{(1)}|  & \dots &  \log |u_{t+1}^{(1)}|\\
\vdots & \dots & \vdots\\
\log |u_1^{(t+1)}| & \dots & \log |u_{t+1}^{(t+1)}|   \\
 \end{array} \right],\]

\vspace{0.4cm}
\noindent hence this determinant is nonzero. 

This construction yields $r$ units $u_1,\dots,u_r$ such that

\begin{equation*}
\det \left[ \begin{array}{ccc}
\log |u_1^{(1)}|  & \dots &  \log |u_r^{(1)}|\\
\vdots & \dots & \vdots\\
\log |u_1^{(r)}| & \dots & \log |u_r^{(r)}|   \\
 \end{array} \right]\not=0.\tag{10}\]

\vspace{0.4cm}
\noindent We claim that $\{u_1,\dots,u_r\}$ is independent: suppose that for some $(a_1,\dots,a_r)\in \mathbb{Z}^r$,
\[
u_1^{a_1}\cdots u_r^{a_r}=1.\]
As before, upon taking the conjugates over $F$ of this equation and then taking the logarithm of the absolute value of the equations resulting from that, we obtain
\[
\sum_{j=1}^r a_j\log |u_j^{(i)}|=0,\ i=1,\dots,r.\]
But the coefficient matrix of this system of linear equations in the $a_i$'s has nonzero determinant by (10), hence $a_1=\cdots=a_r=0$.

\vspace{0.4cm}
\emph{Proof of Lemma} 103. The main idea of this argument is based on the following classical theorem of H. Minkowski in the geometry of numbers:

\vspace{0.4cm}
\textbf{Theorem 104}. \textit{Suppose that}
\[
L_i(x_1,\dots,x_n)=\sum_{j=1}^n a_{ij} x_j\]
\textit{are $n$ linear homogeneous forms defined on $E^n$ with real coefficients $a_{ij}$ such that}
\[
D=\det [a_{ij}]\not=0.\]
\textit{If $(\delta_1,\dots,\delta_n)$ is an $n$-tuple of positive real numbers such that}
\begin{equation*}
\prod_{i=1}^n \delta_i \geq |D|,\tag{$*$}\]
\textit{then there exists $0\not=(z_1,\dots,z_n)\in \mathbb{Z}^n$ such that}
\[
\big|L_i(z_1,\dots,z_n)\big|\leq \delta_i,\ \textrm{for $i=1,\dots,n$.}\]

\vspace{0.4cm}
To derive Lemma 103 from Theorem 104, we start with an integral basis $\{\omega_1,\dots,\omega_n\}$ of $F$ and define the $n$ homogeneous linear forms
\begin{equation*}
L_i(x)=\sum_{j=1}^n \omega_j^{(i)}x_j,\ i=1,\dots,n.\tag{11}\]
These forms do not necessarily have real coefficients, but notice that if $L_i$ has non-real coefficients then the form obtained by replacing all coefficients by their complex conjugates is also on the list (11). If in addition to the hypothesis $(*)$ on the $\delta_i$'s in Theorem 104 we add $\delta_{i+r_2}=\delta_i$, for $i=r_1+1,\dots,r_1+r_2$, one can then deduce from Theorem 104 that the conclusion of Theorem 104 also holds for the forms (11) (see Hecke [7], Chapter V, Theorem 95).

We wish to apply Theorem 104 (as modified so as to be applicable to linear forms with non-real coefficients) to the linear forms (11) for certain choices of the parameters $\delta_1,\dots,\delta_n$. In order to do that, note first that if $d$ is the discriminant of $F$ then
\[
\big|\det[\omega_j^{(i)}]\big|=\sqrt{|d|},\]
and so this determinant is nonzero. Let $D=\max\{2, \sqrt{|d|}\}$ and let $\delta_1,\dots,\delta_n$ be positive real numbers such that 
\begin{equation*}
\prod_{i=1}^n \delta_i=D,\tag{12}\]
\begin{equation*}
\delta_{i+r_2}=\delta_i,\ \textrm{ for $i=r_1+1,\dots,r_1+r_2$.}\tag{13}\]
Theorem 104 then implies that there is $0\not= x=(z_1,\dots,z_n)\in \mathbb{Z}^n$ such that
\[
\big|L_i(z_1,\dots,z_n)\big|\leq \delta_i,\ \textrm{for $i=1,\dots,n$.}\]

Now, let $\alpha=\sum_i z_i \omega_i$. Then $0\not=\alpha\in R$ and $L_i(x)=\alpha^{(i)}$, and so
\begin{equation*}
|\alpha^{(i)}|\leq \delta_i,\ \textrm{for all $i$},\tag{14}\]
\begin{equation*}
1\leq \big|N(\alpha)\big|=\big|\prod_i \alpha^{(i)}\big|\leq \prod_i \delta_i=D.\tag{15}\]
Hence
\begin{equation*}
|\alpha^{(i)}|\geq \frac{1}{\displaystyle{\prod_{j\not= i}|\alpha^{(j)}|}}\geq \frac{1}{\displaystyle{\prod_{j\not= i}\delta_i}}=\frac{\delta_i}{\displaystyle{\prod_{j}\delta_i}}=\frac{\delta_i}{D}\ ,\ \textrm{for all $i$.}\tag{16}\]

Next, let $0\not=(c_1,\dots,c_n)\in E^n$ and define the map $L: R\setminus \{0\}\rightarrow E$ by
\[
L(\gamma)=\sum_{m=1}^r c_m \log |\gamma^{(m)}|.\]
Lemma 103 will be proven provided that we can find a unit $u$ of $R$ such that $L(u)\not=0$. We hence proceed to find such a unit. 

Begin by deducing from (14) and (16) that
\begin{eqnarray*}
(17)\hspace{1.5cm} \left|L(\alpha)-\sum_{m=1}^r c_m \log \delta_m\right|&=&\left|\sum_{m=1}^r c_m \log\frac{\delta_m}{|\alpha^{(m)}|}\right|\\
&\leq&\sum_{m=1}^r |c_m| \log\frac{\delta_m}{|\alpha^{(m)}|}\\
&\leq&(\log D) \sum_{m=1}^r |c_m| \\
&<&C,
\end{eqnarray*}
where $C$ is a fixed constant chosen independently of $\alpha$ and the $\delta_i$'s. 

We now exploit the fact that $\delta_1,\dots,\delta_r$ can be chosen arbitrarily to make certain specific choices of $\delta_1,\dots,\delta_n$ satisfying (12) and (13). Let $k$ be a fixed positive rational integer. First, choose $\delta_{1k},\dots,\delta_{rk}$ positive with
\begin{equation*}
\sum_{m=1}^r c_m \log \delta_{mk}=2Ck.\tag{18}\]
Then choose
\[
\delta_{i+r_2, k}=\delta_{ik}\ \textrm{for $i=r_1+1,\dots,r=r_1+r_2-1$.}\]
It remains to choose $\delta_{r_1+r_2, k}$ and $\delta_{nk}$. Let both of these be
\[
\sqrt{\frac{D}{\prod_{i}\delta_{ik}}}\ .\]
Then $\delta_{1k},\dots,\delta_{nk}$ satisfy (12) and (13), hence (15), (17), and (18) imply that there exist $0\not= \alpha_k\in R$ such that, for all $k$,
\begin{equation*}
\big|L(\alpha_k)-2Ck\big|<C,\tag{19}\]
\begin{equation*}
\big|N(\alpha_k)\big|\leq D.\tag{20}\]
But (19) implies that
\[
C(2k-1)<L(\alpha_k)<C(2k+1),\ \textrm{for all $k$,}\]
hence the sequence $\alpha_1,\alpha_2,\dots$ of nonzero elements of $R$ satisfies
\begin{equation*}
L(\alpha_1)< L(\alpha_2)< L(\alpha_3)<\dots,\tag{21}\]
\begin{equation*}
\textrm{the sequence $\big|N(\alpha_1)\big|, \big|N(\alpha_2)\big|, \big|N(\alpha_3)\big|, \dots$ is bounded.}\tag{22}\]

\emph{Claim}: the principal ideals $(\alpha_1), (\alpha_2),\dots$ cannot all be distinct.

If this is true then $(\alpha_s)=(\alpha_t)$ for at least two distinct $s$ and $t$, hence there is a unit $u$ of $R$ such that $\alpha_s=u\alpha_t$. Then (21) implies that
\[
L(\alpha_t)\not=L(\alpha_s)=L(u\alpha_t),\]
and so
\[
L(u)=L(u\alpha_t)-L(\alpha_t)\not=0.\]

\emph{Proof of the claim}. Suppose this is false. Because the numbers $\big|N(\alpha_1)\big|, \big|N(\alpha_2)\big|,\dots$ are positive rational integers, (22) implies that there are infinitely many $\alpha_s$ such that the ideals $(\alpha_s)$ are all distinct and $\big|N(\alpha_s)\big|$ have the same value $z\in \mathbb{Z}$. For each such $s$,
\[
\pm \alpha_s^{(2)}\cdots\alpha_s^{(n)}=\frac{z}{\alpha_s}\in \mathcal{A}\cap F=R,\]
i.e., $z\in (\alpha_s)$ for all such $s$, which is impossible since $z$ is contained in only finitely many ideals (see the proof of the claim in the proof of Lemma 51).$\hspace{5cm}\textrm{QED}$

\vspace{0.2cm} 
N.B. The special value $r_1+r_2-1$ for $r$ is crucial for the validity of Lemma 103; indeed, the lemma is false for any value of the parameter $r$ in its statement that is larger than  $r_1+r_2-1$.

\vspace{0.2cm} 
\emph{Proof of Theorem} 104. This argument uses the following ingenious geometric idea. With $L_i$ and $\delta_i$  linear forms and parameters which satisfy the hypotheses of Theorem 104, let $\Pi$ denote the parallelotope in $E^n$ defined by
\[
\Pi=\Big\{x\in E^n:\big|L_i(x)\big|\leq \frac{\delta_i}{2},\ \textrm{for all $i$}\Big\}.\]
Then
\[
2 \Pi=\{2x: x\in \Pi\}=\big\{x\in E^n:\big|L_i(x)\big|\leq \delta_i,\ \textrm{for all $i$}\big\},\]
\[
\textrm{volume of $\Pi=|D|^{-1}\prod_{i=1}^n \delta_i$.}\]

Now consider the translates
\[
\Pi_{(a_1,\dots,a_n)}=\big\{(a_1,\dots,a_n)+x: x\in \Pi\big\}\]
of $\Pi$ by all points $(a_1,\dots,a_n)$ of the integer lattice $\mathbb{Z}^n$ of $E^n$. Note that hypothesis $(*)$ in Theorem 104 implies that the volume of $\Pi$ is at least 1. We will show that it follows that at least two translates $\Pi_{(a_1,\dots,a_n)}$ and $\Pi_{(a_1^{\prime},\dots,a_n^{\prime})}$ have a point in common. A simple calculation then verifies that 
\[
0\not=(a_1-a_1^{\prime},\dots,a_n-a_n^{\prime})\in \mathbb{Z}^n\cap (2\Pi)=\mathbb{Z}^n\cap \big\{x\in E^n:\big|L_i(x)\big|\leq \delta_i,\ \textrm{for all $i$}\big\},\] 
which is the conclusion of Theorem 104. 

We must now prove that if the volume of $\Pi$ is at least 1 then at least two of the translates $\Pi_{(a_1,\dots,a_n)}$ are not disjoint. In fact, we do this under the assumption that  the volume of $\Pi$ is greater than 1, i.e., we suppose that
\begin{equation*}
\prod_{i=1}^n \delta_i>|D|\tag{23}\]
(this is, in fact, all that we need for the proof of Lemma 103).

Suppose on the contrary that all of the translated parallelotopes are pairwise disjoint. The strategy of the argument is to then show that this disjointness, together with the fact, from (23), that each parallelotope has volume $>$ 1, prevents the correct number of parallelotopes from fitting inside $n$-cubes centered at the origin in $E^n$ . It was Minkowski's brilliant insight to realize that such ``closest-packing" arguments in geometry can be used to solve important problems in number theory.

For each positive real number $T$, let
\[
S(T)=\{x\in E^n: |x_i|\leq T\}|\]
denote the $n$-cube of side-length $2T$ centered at the origin, and for each positive $k\in \mathbb{Z}$, consider the set of all  parallelotopes $\Pi_{(a_1,\dots,a_n)}$ such that
\begin{equation*}
|a_i|\leq k,\ \textrm{for all $i$}.\tag{24}\]
Note that
\begin{equation*}
\textrm{there are $(2k+1)^n$ of these parallelotopes.}\tag{25}\]

Next, we find the smallest $n$-cube that contains all of these parallelotopes. Let 
\[
c=\textrm{maximum of the absolute value of all coordinates of all points in $\Pi$}.\]
If $\textbf{a}=(a_1,\dots,a_n)\in\mathbb{Z}^n$ then
\[
\Pi_{\textbf{a}}=\{x\in E^n: x-\textbf{a}\in \Pi\},\]
and so if $\textbf{a}$ satisfies (24) then the definition of $c$ implies that
\[
|x_i|\leq |x_i-a_i|+|a_i|\leq c+k,\ \textrm{for all $x\in \Pi_{ \textbf{a}}$,}
\]
hence
\[
\Pi_{\textbf{a}}\subseteq S(c+k),\ \textrm{for all $\textbf{a}$ satisfying (24).}\]
As the $\Pi_{\textbf{a}}$'s  are all pairwise disjoint, we have that
\begin{equation*}
\textrm{sum of the volumes of the $\Pi_{\textbf{a}}$ such that $\textbf{a}$ satisfies (24)}\tag{26}\]
 \[
 \textrm{$\leq$ volume of $S(c+k)=2^n(c+k)^n$.}\]
But
\[
\textrm{volume of $\Pi_{\textbf{a}}=$ volume of $\Pi$, for all $\textbf{a}$,}\]
hence it follows from (25) that
\begin{equation*}
\textrm{sum on the left-hand side of (26)$=(2k+1)^n \times$ volume of $\Pi$,}\tag{27}\]
and so (26) and (27) imply that
\[
\textrm{$(2k+1)^n \times$ volume of $\Pi \leq 2^n(c+k)^n$,} \]
i.e.,
\[
\textrm{ volume of $\Pi \leq \left(\frac{k+c}{k+\frac{1}{2}}\right)^n$.}\]

\vspace{0.2cm}
\noindent Because $k\in \mathbb{Z}$ here is arbitrary, we let $k\rightarrow +\infty$ and conclude that
\[
|D|^{-1}\prod_{i=1}^n \delta_i=\textrm{volume of $\Pi \leq 1$,}\]
which contradicts (23). Hence Theorem 104 is true when strict inequality holds in $(*)$.

Now assume that equality holds in $(*)$. We are to prove: $2\Pi$ contains a nonzero element of $\mathbb{Z}^n$. For $m=1, 2, 3,\dots$, let
\[
\Pi_m=2\left(1+\frac{1}{m}\right)\Pi=\left\{x\in E^n: \big|L_i(x)\big|\leq \left(1+\frac{1}{m}\right)\delta_i,\ \textrm{for all $i$}\right\}.\]
What we just proved implies that there exits $0\not=a_m \in \mathbb{Z}^n\cap \Pi_m$, for all $m$. Now $\Pi_m\subseteq 4\Pi$ and $4\Pi$ is a bounded set in $E^n$, and so the sequence $a_1,\dots,a_m,\dots$ is bounded. A Cantor diagonalization argument on the coordinates of the $a_m$'s then shows that $a_1,\dots,a_m,\dots$ must have terms with the same nonzero value $\bar{a}\in \mathbb{Z}^n$ for infinitely many $m=m_1,m_2,\dots$. Hence $\bar{a}\in \Pi_{m_i}$, for all $i$, i.e.,
\begin{equation*}
\frac{1}{1+\displaystyle{\frac{1}{m_i}}}\ \bar{a}\in 2\Pi,\ \textrm{for all $i$.}\tag{28}\]
But $2\Pi$ is a closed subset of $E^n$, and so (28) implies that\[
\bar{a}=\lim_i\frac{1}{1+\displaystyle{\frac{1}{m_i}}}\ \bar{a}\in 2\Pi.\]
$\hspace{15.7cm}\textrm{QED}$

\newpage
\begin{center}
\emph{Scholium on the Logarithm Space}
\end{center}

\vspace{0.4cm}
If $F$ is a number field, $R=\mathcal{A}\cap F,$ and $n=[F:\mathbb{Q}],$ recall that the logarithm space of $R$, which we will denote by Log $R$, is the set of vectors in $E^n$ defined by
\[
\textrm{Log $R$}=\{(\log |u^{(1)}|,\dots,\log |u^{(n)}|): u \in U(R)\}.\]
If we consider $E^n$ as an abelian group with respect to its vector-space addition, then Log $R$ is a subgroup of $E^n$ and the map
\begin{equation*}
u\rightarrow (\log |u^{(1)}|,\dots,\log |u^{(n)}|)\tag{29}\]
is an epimorphism of $U(R)$ onto Log $R$. It follows from an argument in the proof of Assertion 1 above that the kernel of this epimorphism is $U_0(R)$. Hence as abelian groups,
\[
\textrm{$U(R)\big/U_0(R)$ is canonically isomorphic to Log $R$.}\]
Moreover, it follows from Theorem 102 that  if $r$ is the unital rank of $R$ then $U(R)$ is canonically isomorphic to the direct product $U_0(R)\times \big(U(R)\big/U_0(R)\big)$ and  $U(R)\big/U_0(R)$ is free abelian of rank $r$. If $\{u_1U_0(R),\dots,u_rU_0(R)\}$ is a basis of $U(R)\big/U_0(R)$ (as a free abelian group), or, equivalently, $\{u_1,\dots,u_r\}$ is a set of fundamental units of $U(R)$, and if we set
\[
e_i=(\log |u_i^{(1)}|,\dots,\log |u_i^{(n)}|),\ i=1,\dots,r,\]
then another argument from the proof of Assertion 1 implies that $\{e_1,\dots,e_r\}$ is linearly independent in $E^n$. An application of the epimorphism (29) also allows us to conclude that
\[
\textrm{Log $R=\sum_{i=1}^r \mathbb{Z} e_i$,}\]
i.e., 
\vspace{0.2cm}
\[
\textrm{Log $R$ is an $r$-dimensional vector sublattice of $E^n$ with basis $\{e_1,\dots,e_r\}$.}\]

\vspace{0.2cm}
We collect all of these facts together in

\vspace{0.4cm}
\textbf{Theorem 105}. \textit{$($Structure of the Group of Units$)$ $U(R)$ is canonically isomorphic to the direct product 
\[
U_0(R)\times \big(U(R)\big/U_0(R)\big),\] 
$ U(R)\big/U_0(R)$ is canonically isomorphic to the logarithm space \textnormal{Log} $R$ of $R$, \textnormal{Log} $R$ is an $r$-dimensional  vector sublattice of $E^n$, and if $\{u_1,\dots,u_r\}$ is a set of fundamental units of $U(R)$ then $\big\{(\log |u_1^{(1)}|,\dots,\log |u_1^{(n)}|),\dots,(\log |u_r^{(1)}|,\dots,\log |u_r^{(n)}|)\big\}$ is a basis of \textnormal{Log} $R$.}

\chapter{The Regulator of a Number Field and the Distribution of Ideals}

We conclude these lecture notes with an illustration of how the structure of the group of units determines important and interesting algebraic properties of a number ring.

Let $\varepsilon_1,\dots,\varepsilon_r$ be fundamental units in $R=\mathcal{A}\cap F$ and consider the determinant

\begin{equation*}
\det \left[ \begin{array}{ccc}
\log |\varepsilon_1^{(1)}|  & \dots &  \log |\varepsilon_r^{(1)}|\\
\vdots & \dots & \vdots\\
\log |\varepsilon_1^{(r)}| & \dots & \log |\varepsilon_r^{(r)}|   \\
 \end{array} \right]\tag{1}\]

\vspace{0.5cm} 
\noindent that played an important role in the proof of Theorem 102. We claim that the absolute value of this determinant does not depend on the set of fundamental units used to define it. In order to see this, let $\eta_1,\dots,\eta_r$ be another set of fundamental units. Because $\{\varepsilon_1,\dots,\varepsilon_r\}$ and $\{\eta_1,\dots,\eta_r\}$ are contained in the unique rank-$r$, free-abelian factor of $U(R)$, it follows that there exist $a_{ij}\in \mathbb{Z}$ and $b_{jk}\in\mathbb{Z}$ such that
\begin{equation*}
\eta_i=\prod_{j=1}^r\varepsilon_j^{a_{ij}},\ i=1,\dots,r,\tag{2}\]
\begin{equation*}
\varepsilon_j=\prod_{k=1}^r\eta_k^{b_{jk}},\ j=1,\dots,r.\tag{3}\]
Substitution of (3) into (2) implies that
\begin{equation*}
\eta_i=\prod_{k=1}^r\eta_k^{(\sum_{j=1}^r a_{ij}b_{jk})},\ i=1,\dots,r.\tag{4}\]
Independence of $\{\eta_1,\dots,\eta_r\}$ and (4) imply that
\[
c_{ik}=\sum_{j=1}^r a_{ij}b_{jk}=\left\{\begin{array}{cc}0,\ \textrm{if $i\not=k$,}\\
1,\ \textrm{if $i=k$.}\\\end{array}\right.\]
Hence
\[
1=\det[c_{ij}]=\det[a_{ij}][b_{ij}]=\det[a_{ij}]\det[b_{ij}],\]
and since both determinants on the right side of this equation are in $\mathbb{Z}$, we conclude that
\begin{equation*}
\det [a_{ij}]=\pm1.\tag{5}\]
But (2) implies that
\[
\log\big|\eta_i^{(j)}\big|=\sum_{k=1}^r a_{ik}\log \big|\varepsilon_k^{(j)}\big|,\]
i.e.,
\begin{equation*}
\textrm{transpose of $ \big[\log\big|\eta_j^{(i)}\big|\big]=[a_{ij}]\times$ transpose of $\big[\log\big|\varepsilon_j^{(i)}\big|\big]$.}\tag{6}\]
Hence (5) and (6) imply that
\[
\big|\det \big[\log\big|\eta_j^{(i)}\big|\big]\big|=\big|\det \big[\log\big|\varepsilon_j^{(i)}\big|\big]\big|,\] 
as we claimed. Thus, if 
\[
s=\textrm{$\frac{1}{2}\times$ (number of complex embeddings of $F$),}\]
\[
D=\textrm{absolute value of the determinant $(1)$,}\]
then
\[
\rho=\left\{\begin{array}{cc}D,\ \textrm{if $s=0,\ 1$,}\\
2^{s-1}D,\ \textrm{if $s\geq 2$,}\\\end{array}\right.\]

\vspace{0.2cm}
\noindent is a real-valued and positive invariant of $F$. If $r=0$, i.e., if $F=\mathbb{Q}$ or $F$ is an imaginary quadratic field, then we take $\rho$ to be 1.

\vspace{0.4cm}
\textbf{Definition}. The parameter $\rho$ is called the \emph{regulator of F}.

\vspace{0.4cm}
In order to see what $\rho$ regulates, we turn to the ideals of $R$. If $C$ is an ideal class of $R$ and $t\geq 0$, let
\[
\mathcal{I}_C(t)=\{I\in C: N(I)\leq t\}.\]
We claim that this is a finite set. To verify this, let $J$ be a fixed ideal in $C^{-1}$. Let $0\not= \alpha\in J$. Then there is a unique ideal $I$ such that $(\alpha)=IJ$, hence $I\in C\   ([I]=C[IJ]=C[(1)]=C)$, and 
\begin{equation*}
|N(\alpha)|=N(I)N(J).\tag{7}\]
Moreover, the map $(\alpha)\rightarrow I$ is a bijection of the set of all nonzero principal ideals contained in $J$ onto $C$ (prove!). Also, (7) implies that
\[
N(I)\leq t\ \textrm{if and only if $|N(\alpha)|\leq tN(J)$.}\]
Hence there is a bijection of $\mathcal{I}_C(t)$ onto the set
\[
\mathcal{J}=\big\{\{0\}\not=(\alpha)\subseteq J: |N(\alpha)|\leq tN(J)\big\}.\]
Now the argument in the verification of the claim in the proof of Lemma 103 shows that there is only a finite number of principal ideals of $R$ whose norms do not exceed a fixed constant. Hence $\mathcal{J}$, and so also $\mathcal{I}_C(t)$, is a finite set.

For $k=1, 2, 3,\dots$, set
\[
Z_C(k)=|\mathcal{I}_C(k)|,\]
i.e., $Z_C(k)$ is the number of ideals in the ideal class $C$ whose norms do not exceed $k$. By what we just showed, $Z_C(1),\ Z_C(2),\ Z_C(3),\dots $ is a nondecreasing sequence of positive rational integers. The following remarkable theorem, proved by Dirichlet for quadratic number fields and by Dedekind for all number fields, gives the sharp asymptotic behavior of this sequence.

\vspace{0.4cm}
\textbf{Theorem 106} \textit{$($The Ideal-Class Distribution Theorem$)$. If $C$ is an ideal class of R,}
 \begin{eqnarray*}
 d&=&\ \textnormal{discriminant of $F$},\\
r&=& \ \textnormal{unital rank of $R$},\\
 \rho&=&\ \textnormal{regulator of $F$},\\
 s&=&\ \textnormal{$\frac{1}{2}\times $(number of complex embeddings of $F$ over $\mathbb{Q}$)},\\
 w&=&\ \textnormal{order of the group of roots of unity in $R$},
 \end{eqnarray*}
\textit{then}
\[
 \lim_{k \rightarrow \infty} \frac{Z_C(k)}{k}=\frac{2^{r+1} \pi^s \rho}{w \sqrt {|d|}}\ . \]

\vspace{0.2cm}
\noindent \textit{Moreover, if $\sigma$ denotes this limit and $n=[F:\mathbb{Q}]$, then there exists a constant $M$, depending only on $C$ and $n$, such that}
\begin{equation*}
\left|\frac{Z_C(k)}{k}-\sigma\right|\leq Mk^{-(1/n)},\ k=1, 2, 3,\dots.\tag{8}\]

\vspace{0.4cm}
\emph{Proof}. See Marcus [9], Chapter 6, Theorems 39 and 40 or Hecke [7], Chapter VI, Theorem 121.$\hspace{14.9cm}\textrm{QED}$

Now let
\[
\textrm{$\mathcal{I}(k)=$ the set of all ideals of $R$ whose norms do not exceed $k$,}\]
\[
Z(k)=\big|\mathcal{I}(k)\big|,\ k=1, 2, 3,\dots.\]
If $h$ is the class number of $R$ and $C_1,\dots,C_h$ are the ideal classes of $R$, then $\mathcal{I}(k)$ is the pairwise disjoint union
\[
\bigcup_{i=1}^h\  \mathcal{I}_{C_i}(k),\]
and so
\[
Z(k)=\sum_{i=1}^hZ_{C_i}(k),\ k=1, 2, 3,\dots.\]

\vspace{0.2cm}
\noindent We hence deduce as an immediate consequence of Theorem 106

\vspace{0.4cm}
\textbf{Theorem 107} \textit{$($The Ideal Distribution Theorem$)$. If $h$ is the class number of $R$ then}
\[
\lim_{k \rightarrow \infty} \frac{Z(k)}{k}=\sigma h.\]
\textit{Moreover,}
\[
\left|\frac{Z(k)}{k}-\sigma h\right|\leq Mhk^{-(1/n)},\ k=1, 2, 3,\dots,\]
\textit{where $M$ is any constant for which the error estimate $(8)$ in Theorem $106$ is valid for all ideal classes of $R$. }

\vspace{0.4cm}
 Thus the structural parameters of the group of units, together with the discriminant of $F$ and the class number of $R$, determine the sharp asymptotic distribution of the number of ideals of $R$ and the number of ideals in each ideal class whose norms do not exceed a series of fixed constants. As we have seen throughout these notes, the discriminant of $F$, the class number of $R$, the unital rank of $R$, the regulator of $F$, the number of complex embeddings of $F$, and the order of the group of roots of unity in $R$ are fundamental parameters associated with $F$ which govern many aspects of the arithmetic and algebraic structure of $F$ and $R$. Theorems 106 and 107 are two remarkable examples of how all of those parameters work in concert to do that.  

Let $F=\mathbb{Q}(\sqrt m\ )$ be a real quadratic field. Then $r=1$ and $s=0$, hence if $\lambda>1$ is the fundamental unit of $F$ then the regulator of $F$ is $\log \lambda$. Also,  $d=m$ or $4m$ if, respectively,  $m\equiv 1$ mod 4 or  $m\not \equiv 1$ mod 4, and $w=2$. If $h$ is the class number of $\mathcal{A}\cap F$ then Theorem 107 implies that
\begin{equation*}
\lim_{k \rightarrow \infty} \frac{Z(k)}{k}=\left\{\begin{array}{cc}\displaystyle{\frac{2\log \lambda}{\sqrt m}}h,\ \textrm{if  $m\equiv 1$ mod 4,}\\
\vspace{0.01cm}\\
\displaystyle{\frac{\log \lambda}{\sqrt m}h},\ \textrm{if $m\not \equiv 1 $ mod 4.}\\\end{array}\right.\tag{9}\]
Now the fundamental unit $\lambda$ can be readily computed by solving the appropriate Pell equation, and so if the limit on the left-hand side of (9) can be estimated accurately enough, then the class number $h$ can be found. This strategy for calculation of the class number is what motivated Dirichlet to prove Theorems 106 and 107 for quadratic fields.

\vspace{2cm}
\begin{center}
\textit{FINIS}
\end{center}

\begin{theindex}
\item algebraic integer, 30
\item algebraic number, 12
\item algebraic number field, 15
\item algebraic number over $F$, 8
\subitem conjugates of an, 10
\subitem degree of an, 9
\item associates (in an integral domain), 4
\item characteristic of a ring, 72
\item class number, 66
\item class-number problem, 69
\item class-number 1 problem, 70
\item completely ramified, 79
\item complex embedding, 90
\item complex number field, 8
\subitem extension of a, 14
\subitem simple extension of a, 14
\item conjugates of $\alpha$ over $K$, $\alpha \in K$, 24
\item cyclotomic number field, 16

\item Dedekind domain, 62
\item degree of a prime ideal, 73
\item discriminant of $(\alpha_1,\dots,\alpha_n)$, 24
\item discriminant of $F$, 39
\item embedding of $K$ over $F$, 21

\item Euler's Theorem, 84
\item Euler's totient function, 16
\item Fermat's Last Theorem, 4
\item Fermat's Little Theorem, 81
\item field of fractions, 62
\item field polynomial, 35
\item FLT for primes, 4
\item fundamental unit(s), 103, 104, 106
\item Galois automorphism, 77
\item Galois field, 77
\item Galois group, 77
\item Gaussian integers, 34
\item Gauss' Lemma, 31
\item greatest common divisor of ideals, 58
\item group of units, 4
\item ideal, 48
\subitem basis of an, 49
\subitem generators of an, 48
\subitem maximal, 50
\subitem prime, 50
\subitem principal, 48
\item ideal class, 65
\item ideal-class group, 66
\item ideal-class product, 65
\item ideal product, 49
\item independent subset of an abelian group, 107
\item inertial, 79
\item integral basis, 37
\item integral domain, 4
\item integrally closed, 62
\item irreducible polynomial, 9
\item logarithm space, 119
\item minimal polynomial, 8
\item Minkowski's constant, 91
\item multiplicative group of a field, 100
\item multiplicity of a prime ideal, 55
\item norm of $\alpha$ over $F$, 22
\item norm of an ideal, 90

\item notation
\subitem $\mathcal{A}$, 30
\subitem $A[x],\ A$ a commutative ring, 8
\subitem \textbf{C}, 6
\subitem $\det M$, 24
\subitem $E$, 90
\subitem $F_n[x]$, 21
\subitem $S^n$, $S$ a set, 5
\subitem $|S|$,  $S$ a set, 71
\subitem $\mathbb{Q}$, 7
\subitem $\mathbb{Z}$, 5
\item numbered corollaries, lemmas, propositions, and theorems
\subitem Proposition 1, 8
\subitem Corollary 2, 9
\subitem Proposition 3, 9
\subitem Proposition 4, 9
\subitem Theorem 5, 10
\subitem Proposition 6,10
\subitem Corollary 7, 11
\subitem Lemma 8, 11
\subitem Theorem 9, 12
\subitem Proposition 10, 14
\subitem Lemma 11, 16
\subitem Proposition 12, 18
\subitem Lemma 13, 19
\subitem Proposition 14, 21
\subitem Proposition 15, 22
\subitem Proposition 16, 25
\subitem Proposition 17, 25
\subitem Corollary 18, 25
\subitem Proposition 19, 26
\subitem Proposition 20, 28
\subitem Theorem 21, 30
\subitem Lemma 22, 30
\subitem Proposition 23, 32
\subitem Corollary 24, 34
\subitem Proposition 25, 34
\subitem Proposition 26, 34
\subitem Lemma 27, 35
\subitem Lemma 28, 36
\subitem Proposition 29, 37
\subitem Corollary 30, 37
\subitem Theorem 31, 37
\subitem Proposition 32, 38
\subitem Proposition 33, 39
\subitem Lemma 34, 39
\subitem Proposition 35 (Eisenstein's Irreducibility Criterion), 42
\subitem Proposition 36, 43
\subitem Proposition 37, 43
\subitem Proposition 38, 44
\subitem Theorem 39, 44
\subitem Proposition 40, 45
\subitem Proposition 41, 46
\subitem Lemma 42, 46
\subitem Lemma 43, 46
\subitem Proposition 44, 48
\subitem Theorem 45, 48
\subitem Proposition 46, 49
\subitem Proposition 47, 50
\subitem Proposition 48, 50
\subitem Proposition 49, 50
\subitem Theorem 50 (Fundamental Theorem of Ideal Theory), 52
\subitem Lemma 51, 52
\subitem Lemma 52, 53
\subitem Lemma 53, 54
\subitem Proposition 54, 54
\subitem Proposition 55, 54
\subitem Proposition 56, 55
\subitem Proposition 57, 57
\subitem Corollary 58, 58
\subitem Proposition 59, 58
\subitem Corollary 60, 59
\subitem Lemma 61, 59
\subitem Theorem 62 (Ideal Generation Theorem), 60
\subitem Proposition 63, 62
\subitem Theorem 64  (Fundamental Theorem of Ideal Theory for Dedekind Domains), 63
\subitem Proposition 65, 65
\subitem Lemma 66, 65
\subitem Proposition 67, 65
\subitem Proposition 68, 66
\subitem Theorem 69 (Finiteness of the Class Number), 66
\subitem Lemma 70, 66
\subitem Lemma 71, 68
\subitem Theorem 72, 70
\subitem Corollary 73 (Kummer's Conjecture), 70
\subitem Theorem 74, 70
\subitem Proposition 75, 71
\subitem Lemma 76, 72
\subitem Lemma 77, 73
\subitem Lemma 78, 73
\subitem Theorem 79 (Ramification Equation), 74
\subitem Lemma 80 (Chinese Remainder Theorem for Commutative Rings), 74
\subitem Lemma 81, 74
\subitem Lemma 82, 76
\subitem Proposition 83, 77
\subitem Corollary 84, 77
\subitem Proposition 85, 77
\subitem Theorem 86, 78
\subitem Lemma 87, 78
\subitem Theorem 88, 81
\subitem Corollary 89, 82
\subitem Lemma 90, 82
\subitem Theorem 91, 84
\subitem Corollary 92, 85
\subitem Theorem 93, 86
\subitem Theorem 94, 88
\subitem Lemma 95, 91
\subitem Proposition 96, 102
\subitem Proposition 97, 102
\subitem Lemma 98, 102
\subitem Proposition 99, 104
\subitem Proposition 100, 104
\subitem Lemma 101, 105
\subitem Theorem 102 (Dedekind-Dirichlet Units Theorem), 106
\subitem Lemma 103, 110
\subitem Theorem 104, 112
\subitem Theorem 105 (Structure of the Group of Units), 118
\subitem Theorem 106 (Ideal-Class Distribution Theorem), 121
\subitem Theorem 107 (Ideal Distribution Theorem), 122
\item number ring, 31
\item order of $a$ modulo $m$, 81, 84

\item Pell's equation, 102
\item prime (in an integral domain), 4

\item primitive element, 18
\item Primitive Element Theorem, 18

\item primitive root, 81, 84
\item principal class, 65
\item principal-ideal domain (PID), 48
\item quadratic number field, 15
\item quotient of ideals, 58
\item ramified, 71
\item ramification index, 71
\item rational integer, 31
\item real embedding.90
\item regulator, 120
\item ring of integers, 31

\item split, 79
\item symmetric polynomial, 10
\subitem elementary, 10
\item torsion-free, 106
\item torsion subgroup, 106
\item trace of $\alpha$ over $F$, 22
\item transcendental number, 12
\item unique factorization domain (UFD), 5

\item unit, 4
\item unital rank, 106
\item unramified, 71
\item vector lattice, 96
\subitem co-volume of a, 96

\end{theindex}

\end{document}